%% file: finalv7_PV.tex
\documentclass[english,11pt]{elsarticle}

\usepackage[T1]{fontenc}
\usepackage[latin9]{inputenc} 
\usepackage{float}
\usepackage{amsmath}
\usepackage{amsthm}
\usepackage{amssymb}
\usepackage{amsfonts}
\usepackage{bbm}
\usepackage{graphicx}
\usepackage{epstopdf}
\usepackage[pdftex,colorlinks,linkcolor={blue},citecolor={blue}]{hyperref}
\usepackage{tikz, xcolor}
\usepackage{tkz-euclide}
\usetkzobj{all}
\usetikzlibrary{shapes,arrows,calc,intersections,through,backgrounds}

 \usepackage{pgfplots}
  \pgfplotsset{compat=newest}
  \usetikzlibrary{plotmarks}
  \usetikzlibrary{arrows.meta}
  \usepgfplotslibrary{patchplots}
  \usepackage{grffile}
  \usepackage{amsmath}
\long\def\comment#1{}

\input{format}

\newcommand{\blue}[1]{\textcolor{blue}{#1}}

\begin{document}
\begin{frontmatter}
 \title{\textbf{Regulating TNCs: Should Uber and Lyft Set Their Own Rules?}}

\author[1staddress]{Sen Li}
\ead{lisen1990@berkeley.edu}
\author[1staddress]{Hamidreza Tavafoghi}
\ead{tavaf@berkeley.edu}
\author[1staddress,2ndaddress]{Kameshwar Poolla}
\ead{poolla@berkeley.edu}
\author[2ndaddress]{Pravin Varaiya}
\ead{varaiya@berkeley.edu}

\address[1staddress]{Department of Mechanical Engineering, University of California, Berkeley }
\address[2ndaddress]{Department of Electrical Engineering and Computer Science, University of California, Berkeley}


\begin{abstract}
We evaluate the impact of three proposed regulations of transportation network companies (TNCs) like Uber, Lyft and Didi:  (1) a minimum wage for drivers,  (2) a cap on the number of drivers or vehicles, and (3) a per-trip congestion tax. The impact  is assessed using a queuing theoretic equilibrium model which incorporates the stochastic dynamics of the app-based ride-hailing matching platform, the ride prices and driver wages established by the platform, and the incentives of passengers and drivers. We show that a  floor placed under driver earnings  pushes the ride-hailing platform to hire {\em more} drivers and offer {\em more} rides, at the same time that  passengers  enjoy  {\em faster} rides and {\em lower} total cost, while platform rents are reduced. Contrary to standard economic theory,  enforcing a minimum wage  for drivers  benefits both  drivers and passengers, and promotes the  efficiency of the entire system. This surprising outcome holds for almost all model parameters, and  it occurs  because the wage floors curbs TNC labor market power. In contrast to a wage floor, imposing a cap on the number of vehicles  hurts drivers, because the platform  reaps all the benefits of limiting  supply.  
The congestion tax has the expected  impact: fares increase, wages and platform revenue decrease.   
We also construct variants of the model to briefly discuss platform subsidy, platform competition, and autonomous vehicles.

\end{abstract}


\end{frontmatter}

\section{Introduction}

In December 2018, New York  became the first US city to adopt a minimum wage for drivers working for app-based transportation network companies (TNCs) like Uber and Lyft. The New York City Taxi and Limousine Commission (NYTLC) established a ``minimum per-trip payment formula'' that gives an estimated gross hourly driver earnings before expenses of at least \$27.86 per hour and a net income of 
\$17.22 per hour after expenses, equivalent to the minimum wage of \$15 per hour because, as ``independent contractors,'' drivers pay additional payroll taxes and get no paid time off \blue{\cite{ban2018gan}.} 
The NYTLC formula for non-wheelchair accessible vehicles is
\begin{equation}\label{utilization}
\mbox{Driver pay per trip} =
\Bigg( \frac{\$0.631 \times \text{Trip Miles}}{\text{Company Utilization Rate}}\Bigg)
+ \Bigg( \frac{\$0.287 \times \text{Trip Minutes}}{\text{Company Utilization Rate}}\Bigg)
+ \text{Shared Ride Bonus}
\end{equation}
amounting to \$23 for a 30-min, 7.5-mile ride.\footnote{The utilization rate is calculated as the total amount of time drivers spend transporting passengers on trips dispatched by the base divided by the total amount of time drivers are available to accept dispatches from the base \cite{ban2018gan}. Wheelchair accessible vehicles receive a higher rate.} New York City's \$15/hour minimum wage for large employers, which went into effect on December 31, 2018 doesn't apply to drivers who work for ride-hailing apps.

The Commission imposed this wage floor based on  testimony  on driver expenses, meetings with stakeholders, and on the report of  labor economists J.A. Parrott and M. Reich \blue{\cite{parrott2018earning}} which showed that
median driver earnings had declined almost \$3.00 per hour from \$25.78 in September  2016 to \$ 22.90 in October 2017, a decrease of 11.17\%.   The TNCs imposed the \$3.00 per hour wage cut  
during a period when the number of  drivers in the largest four TNCs (Uber, Lyft, Gett/Juno, and Via) had grown by 80,000 \cite{ban2018gan}. Uber  would be the largest for-profit private employer in New York City if its drivers were classified as employees rather than independent contractors \cite{parrott2018earning}. The ingenious wage formula \eqref{utilization} encourages TNCs to increase driver pay through higher utilization, instead of trying to restrict the number of drivers through regulation.   Lyft, however, opposed the regulation saying that because of its larger size, Uber's higher utilization rate gave it an unfair advantage \cite{Lyft_lawsuit}.  Lyft's complaint was overruled \cite{Lyft_lawsuit_dropped}.  

The subminimum wage of drivers working for TNCs also prompted the Seattle City Council in April 2018  to pass a unanimous resolution to explore setting a minimum base rate  of \$2.40 per mile for TNCs compared with the prevailing rate of \$1.35 per mile and the rate of \$2.70 per mile charged by taxis \blue{\cite{stranger.2018}}.  The resolution also asked TNCs to voluntarily hand over anonymous data on hours, trips, fares and compensation.  Unlike NYTLC, however, no other US city has access to TNC data to estimate what their drivers are paid or the TNC impact on traffic. For example, the California Public Utilities Commission which regulates TNCs will not share TNC data with San Francisco County Transportation Authority \cite{TNC_datasharing}.
TNC regulation ``follows an elite political process dominated by concentrated actors and government decision makers largely acting ex officio (committee heads, regulators, and judges)'' \cite{Uber_Collier}.

In December 2018, Uber lost its case at the U.K. Court of Appeal  against the  October 2016 ruling that its drivers should be classified as workers entitled to rights such as minimum wage and paid holidays.  The Court ruled against Uber's claim that its drivers were just self-employed contractors who use its app in exchange for a share of their fares at the level dictated by Uber \cite{Uber_London}. The case can be used to challenge the self-employed status of millions of gig-economy workers who  work for companies like Airbnb and Deliveroo on a freelance basis without fixed contracts.
New York and London are the largest Uber markets in the US and EU.
The California state assembly recently passed  bill AB5 that would make hundreds of thousands of independent contractors including TNC drivers become employees. The bill now goes to the senate \cite{vox_gig}. Uber and Lyft are aggressively campaigning against AB5. In its SEC filing, Uber states ``If, as a result of legislation or judicial decisions, we are required to classify Drivers as employees $\ldots $ we would incur significant additional expenses [that would] require us to fundamentally change our business model, and consequently have an adverse effect on our business and financial condition \cite[p.28]{Uber_Sec}.''\footnote{For a thoughtful discussion of labor-market trends in the gig-economy see \cite{Tirole_gig}.}

As of January 1, 2019, all trips by \textit{for-hire} vehicles that cross 96th street in NYC will pay a congestion surcharge of \$2.75 per TNC trip, \$2.50 per taxi trip, and \$0.75 per pool trip. Further,
NYC will also charge a toll on \textit{every} vehicle that enters the busiest areas, currently defined as south of 61st street.  This `cordon' price will raise about \$1B per year (assuming a \$11.52 toll) for the Metropolitan Transportation Authority.

Uber's reaction to these  adverse decisions was predictable.
Responding to the NYTLC ruling Uber's director of public affairs stated, ``legislation to increase driver earnings will lead to higher than necessary fare increases for riders while missing an opportunity to deal with congestion in Manhattan's central business district'' \blue{\cite{Uber_USA}.}\footnote{Lyft echoed the Uber response 
stating, ``These rules would be a step backward for New Yorkers, and we urge the TLC to reconsider them \cite{Uber_USA}.''}  Uber  challenged the Seattle resolution: its general manager for  Seattle said, ``we are generally unclear how nearly doubling per-mile rider rates would not result in an increased cost for riders''\blue{\blue{\cite{stranger.2018}}.}  Uber also declared it would fight the U.K. Appeal Court's decision in the Supreme Court \blue{\cite{Uber_London}.} Contradicting Uber's claims, this study shows that raising driver wages will \textit{increase} the number of drivers and riders at the same time that passengers enjoy \textit{faster} rides and \textit{lower} total cost, while platform rents are reduced. 

The aforementioned regulations are part of the political response to the public anxiety over the disruption of the urban transportation system caused by the rapid growth of TNCs. Worldwide, the
monthly number of Uber users  is forecast to reach 100 million in 2018, up from 75 million in 2017. In New York, the four largest TNCs Uber, Lyft,  Juno and Via combined dispatched nearly 600,000  rides per day in the first quarter of 2018, increasing their annual trip totals by over 100 percent in 2016 and by 71 percent in 2017. About 80,000 vehicles are affiliated with these four companies \cite{parrott2018earning}.  In San Francisco, 5,700 TNC vehicles operate in peak times.  They daily make over 170,000 vehicle trips, approximately 12 times the number of taxi trips, and 15 percent of all intra-San Francisco  trips, comprising  
at least 9 percent of all San Francisco person trips \cite{castiglione2016tncs}.  This explosive growth of TNCs has raised two  public concerns.  

As noted earlier, one concern is with the working conditions of TNC drivers.  The TNC  business model  places much of the economic risk associated with the app sector on  drivers, who are classified as independent contractors. Furthermore, the model relies on having many idle cars and drivers, resulting in low driver pay per hour and high TNC platform rents.\footnote{TNC expenditures comprise a fixed initial cost for setting up the platform and a small variable cost as the company grows.  Thus the  average cost per trip falls and its profit margin increases as the TNC grows.} TNCs need idle drivers to reduce passenger waiting time. Uber's annual revenue from passenger fares in New York City amounts to about \$2
billion, of which it keeps about \$375 million in commissions and fees, for a markup  estimated at six times its variable operating cost or 600 percent \cite{parrott2018earning}. 
One common opinion is that ``Uber's driver-partners are attracted to the flexible schedules that driving on the Uber platform affords \ldots because the nature of the work, the flexibility, and the compensation appeals to them compared with other available options \cite{Hall_Kreuger}.''
  In fact, more than 60 percent of  New York City drivers  work full-time and provide 80 percent of all rides; their work hours are not flexible \cite{parrott2018earning}.

The second concern is with the negative impact of   TNCs on a city's traffic congestion and its public transit ridership.
A detailed 2017 report  \cite{schaller2017empty}  examined the impact of TNC growth on traffic conditions in Manhattan's CBD. The analysis shows that, from 2013 to 2017, TNC trips increased  15 percent, VMT increased  36 percent, traffic speed declined 15 percent, the number of vehicles increased  59 percent, and the number of unoccupied vehicles increased  81 percent. The report suggested reducing the unoccupied time of TNC vehicles as a means of congestion control. Responding to the increased congestion, the New York City Council in 2018 passed  a regulation freezing the number of TNC vehicles on the road for one year. Supporters of the cap, including Mayor Bill de Blasio,  said the regulation will protect drivers, fairly regulate the industry and reduce congestion \cite{cnbc_NY}.  However, our analysis shows that imposing a cap hurts drivers, because the TNC retains as profit the benefits of limiting supply.

Another detailed report \cite{castiglione2016tncs} by San Francisco Transportation Authority  provides information on the size, location, and time-of-day characteristics of TNC activities in San Francisco. A follow-up report  \cite{castiglione2018tncs} identifies the impact of TNC activities on road congestion in San Francisco. It shows that after subtracting the impact of employment growth, population change and network capacity change, TNCs contributed 51 percent of the increase in vehicle hours of delay, 47 percent of increase in VMT, and 55 percent of the average speed decline  between 2010 and 2016.  Moreover, ``TNC trips are concentrated in the densest and most congested parts of San Francisco including the downtown and northeastern core of the city. At peak periods, TNCs are estimated to comprise 25 percent of vehicle trips in South of Market.'' The report cites studies showing that ``between 43 percent and 61 percent of TNC trips substitute for transit, walk, or bike travel or would not have been made at all.''

This paper evaluates three TNC regulations: a minimum driver wage,  a cap on the  number of drivers or vehicles, and a per-trip congestion tax. We analyze the  impacts of these regulations on several aspects of the  app-based 
ride-hailing market, including  ride prices and driver wages established by the platform,  the incentives of passengers and drivers, vehicle occupancy, and platform rent or profit. We use a model to  determine 
the arrival of passengers, number of drivers, ride prices and platform commissions, conditioned on the imposed regulation.  The model employs 
 a queuing theoretic model with dynamic matching of passengers and drivers,  an equilibrium model  that predicts the long-term average arrivals of passengers and drivers, and an optimization model of platform decision-making. We summarize the key results.
\begin{itemize}
\item Imposing a minimum wage will motivate TNCs to hire {\em more} drivers and offer {\em more} rides, and passengers to enjoy {\em faster}  rides and {\em lower} total cost, while TNC rent or profit shrinks. It indicates that raising the minimum wage   will benefit both  drivers and passengers, while platform rent will decline. This counter-intuitive result holds for almost all model parameters, and it occurs because the  wage floor curbs TNC labor market power. 
\item Contrary to  common belief,   a cap on the number of drivers will  hurt  driver earnings.   This is because when fewer drivers are permitted, the platform will hire cheaper labor by reducing  driver pay. Thus, the benefit of limiting the driver supply is retained by the platform.    
\item Imposing a congestion surcharge has a predictable impact:  the numbers of passengers and drivers  and  the platform revenue reduce as the congestion surcharge increases. Our numerical study shows that a congestion surcharge of $\$2.75$/trip significantly reduces the platform profit in NYC. This suggests that the business model of TNC is  vulnerable to the adverse effect of congestion policies.
\end{itemize}

We also present variants of our model to analyze platform subsidy, platform competition and autonomous vehicles.

{\bf Related Work}: There are several studies of ride-hailing platforms. A recurrent concern is to evaluate decisions that maximize platform profit, with particular attention to static vs. dynamic pricing.  
A queuing  model is proposed in \cite{banerjee2015pricing} to study the profit maximizing prices of ride-hailing platforms. It shows that the throughput and profit under dynamic pricing strategy can not exceed that under the optimal static pricing strategy that is agnostic to stochastic dynamics of demands. On the other hand,  dynamic pricing is  more robust to fluctuations in system parameters compared to static pricing. Hence, the platform can use dynamic pricing to realize the benefits of optimal static pricing without perfect knowledge of system parameters. 

A similar question is studied in \cite{cachon2017role}, with a  focus on the self-scheduling capacity of for-hire drivers. It is shown that the additional flexibility of drivers is beneficial to   platforms, consumers and drivers.  It also suggests that when some periods have predictably higher demand than others (e.g., a rainy evening),  with static pricing it is hard to find service at peak demand times, so  surge pricing is likely to benefit all stakeholders.  In the same vein, 
\cite{bai2018coordinating}  suggests dynamic pricing for  the platform to maximize the profit across different time periods when the underlying operating characteristics  change significantly.   It is shown in \cite{taylor2018demand} that platform pricing can be more complicated when there is uncertainty in passenger's valuation or driver's opportunity cost. 
A general economic equilibrium model is developed in \cite{nycdoc} to evaluate the impacts of ride-hailing services on   deadhead miles and traffic congestion. 
Ride-hailing platforms are also examined as a special kind of two-sided platforms. See \cite{rysman2009economics} and \cite{rochet2006two} for a  summary of literature on two-sided platforms, and \cite{weyl2010price} for a general theory of monopoly pricing in multi-sided platforms.

The literature on regulation of the app-based ride-hailing marketplace is relatively limited. A ride-hailing platform that manages a group of self-scheduling drivers to serve  time-varying demand is studied in \cite{gurvich2016operations}. The study shows that under a wage floor, the platform starts to limit agent flexibility because it limits the number of agents that can work in some time intervals.

The  work closest to ours is by Parrott and Reich \cite{parrott2018earning}. The authors use TNC administrative data collected by the New York City Taxi and Limousine Commission (NYTLC) to examine the likely impact of the NYTLC's proposed regulations \cite{ban2018gan}. By numerical simulation, they show that the proposed policy will increase driver earnings by 22.5 percent, while  passengers will only experience moderate increase of trip fare (less than 5 percent) and waiting times (12 to 15 seconds). However, our analysis shows that both the trip cost  and the waiting time will decrease. This is because in our model we assume that the passengers are sensitive to the pickup time of the ride-hailing services, which is not captured in \cite{parrott2018earning}.

\section{TNC Environment}\label{sec2}

\begin{figure}[t!]%
\centering
\includegraphics[scale=0.5]{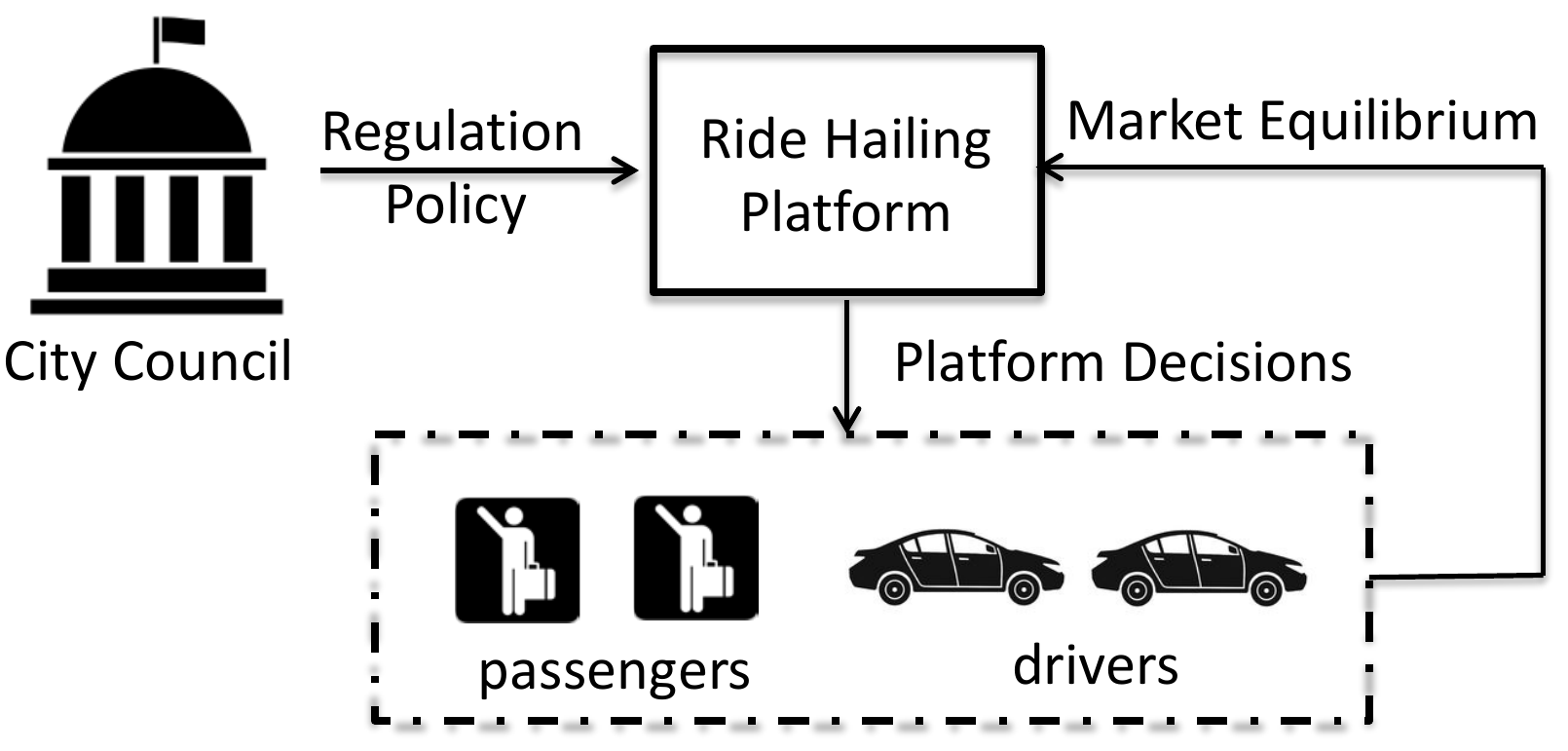}
\caption{The TNC system  includes the city council, platform, passengers and drivers.}
\label{systemdiagram}
\end{figure}

This section describes the TNC enviroment. Agents of the transportation system are comprised of the city council, the app-based ride-hailing platform (TNC), a group of passengers and for-hire vehicle drivers (see Figure \ref{systemdiagram}).  The city council sets legislation to regulate TNC operations. Examples of regulations include minimum driver wage, maximum number of vehicles and regional licensing.\footnote{Unlike TNCs, taxicabs are heavily regulated.}  The regulations are enforced by auditing the operational data of TNCs. (See \cite{ban2018gan} for details of enforcement in New York.) The platform responds to the regulations by  setting profit-maximizing ride fares and driver commissions (or equivalently, wages). These fares and wages are called `platform decisions' in Figure~ \ref{systemdiagram}. The platform decisions  influence the choices of passengers and drivers. For instance, passengers  have diverse ride choices including TNC, public transit, walking,  and biking. They select an option based on the cost and convenience of each choice. Drivers also have alternative job opportunities, such as delivering food, grocery, packages, and mail. They take the job with the highest expected wage. The  choices of passengers and drivers form a market equilibrium, which determines the TNC profit or rent. The equilibrium is affected by regulations.

The objective of the paper is  to understand how regulations impact the ride-hailing transportation system.  We  consider  three  regulations: (a) a floor under driver wage; (b) a cap on total number of drivers; and 
(c) a per-trip congestion tax. We  analyze their impact from various perspectives of the ride-hailing system, including  ride fares, commission rate, passenger pickup time, driver earnings, platform rent, number of riders, number of for-hire vehicles, and vehicle occupancy rate. 

The rest of this paper is organized as follows. In Section \ref{lowerlevel} we introduce the market equilibrium model of the response of passengers and drivers to a  platform decision. In Section \ref{higherlevel} we predict  TNC  decisions in the absence of regulation. In Section \ref{sec5} we  examine TNC decisions operating with  a floor under driver wage rate.  In Section \ref{sec6} we consider TNC decisions when there is a cap on the number of drivers.In Section \ref{sec7} we study the impact of congestion surchage.
Platform competition and other model variations are discussed in Section~\ref{extension}. Conclusions are offered in Section \ref{section_conclusion}. Several proofs are deferred to the appendix.

\section{Market Equilibrium Model} 
\label{lowerlevel}

We now develop the market equilibrium model of the  decisions of  drivers and passengers in  response to the platform decision. The model is used to predict the average arrival rates of passengers and number of drivers.

\subsection{Matching Passengers and Drivers: M/G/N Queue}
We use a continuous-time queuing process to model the matching of passengers and drivers. Consider $N$ TNC drivers or vehicles,  each modeled as a server. A vehicle is `busy' if there is a passenger on board, or a passenger is assigned and the vehicle is on its way to picking her up.  Otherwise, it is considered `idle'.  We assume that the arrival process of passengers is Poisson with rate $\lambda>0$.  Newly arrived passengers immediately join the queue and wait until an idle vehicle is dispatched by the platform. Hence this is an M/G/N queue, and the expected number of idle servers (vehicles) is $N_{I}=N-\lambda/\mu$, where $\mu^{-1}$ is the average trip duration.

\begin{remark}
To ensure stability of the queue the model requires  $N_I>0$, i.e., $N>\lambda/\mu$. This is consistent with the TNC business model that ``relies upon very short wait times for passengers requesting rides, which in turn depends on a large supply of available but idle drivers and vehicles''  \citep{parrott2018earning}. For instance,   New York has an average of 5089 TNC vehicle  \citep{schaller2017empty}, 187 passenger per minute, and a trip takes 16.3 minutes, i.e., $\mu=1/16.3 \text{ min}^{-1}$. This gives  $N_I=N-\lambda/u=2041$. Cities with limited supply of drivers (as in US suburbs and in cities like Singapore) require a distinct model \cite{banerjee2015pricing}, \cite{zha2016economic}.  
\end{remark}

\subsection{Passenger and Driver Incentives}
The passenger arrival rate $\lambda$ and the number of drivers $N$ are endogenously determined in the market equilibrium.

{\bf Passenger Incentives}: 
Passengers  choose their rides from available options like app-based TNCs, public transit, walking, or biking, by comparing their prices and waiting times. We model the cost of the app-based ride-hailing service as 
\begin{equation}
\label{costoftravel}
c=\alpha t_w+\beta p_f,
\end{equation}
where $t_w$ is the average waiting time (from sending a request to being picked up), $p_f$ is the per trip  fare of the ride-hailing service \footnote{Most app-based ride-hailing platforms charge passengers based on the  formula: $\text{total cost} = \text{base fare} + \text{price/mile} \times \text{trip miles} +  \text{price/time} \times \text{trip time}$. $p_f$ represents the sum of these three costs.}, and $\alpha$ and $\beta$ specify the passenger's trade-off between convenience and money. We refer to $c$ as the total cost of a TNC trip, including the trip fare plus the money value of the trip time.

In the ride-hailing service, a  ride  is initiated when a passenger sends  a  request to the platform, and is completed when the passenger is dropped off at the destination. We divide a ride into three periods: (1) from the ride request being received to a vehicle  being assigned; (2) from a vehicle  being assigned  to passenger pickup; (3) from passenger pickup to passenger drop-off.  Let $t_m$, $t_p$ and $t_o$ be the average duration of these periods. Here $t_m$ is the average waiting time in the M/G/N queue. Assuming that the platform matches the passenger to the nearest idle vehicle,  $t_p$ depends on the  distance of the nearest idle vehicle to the passenger. Typically,  $t_m$ ranges from a couple of seconds to a half minute, and $t_p$ is between three to six minutes. The sum of $t_m$ and $t_p$ is the passenger waiting time, denoted as $t_w=t_m+t_p$. Clearly, the passenger waiting time depends on the average number of idle vehicles $N_I$. We denote $t_w$ as a function of $N_I$, i.e., $t_w(N_I)$, and impose the following assumption
\begin{assumption}
\label{assump_pickuptimefunc}
The function $t_w: \mathbb{R}_+\rightarrow \mathbb{R}_+$ is convex, decreasing and twice differentiable. 
\end{assumption}
This assumption says that passenger waiting time decreases as $N_I$ increases, and the marginal benefit of recruiting extra vehicles to reduce waiting time diminishes as $N_I$ increases. It is a standing assumption throughout the paper. 

In some special cases, the waiting time function $t_w(\cdot)$ can be derived analytically. Let $d(x)$ denote the  distance of a passenger requesting a ride at location $x$ in a city to the nearest idle vehicle. Let $N_{I_0}$ be the average number of idle TNC vehicles before regulatory intervention (more precisely specified later). We have the following proposition.

\begin{proposition}\label{pickup_theory}
Consider a city with an arbitrary geometry . Assume that (1) the platform  matches each arriving passenger to the nearest idle vehicle available; (2) idle vehicles are uniformly and independently distributed across the city; (3) location of  passengers requesting a ride is uniformly distributed across the city,  independently of the position of idle vehicles. Then,
	\begin{align}\label{pickup-exact}
	\hspace*{-7pt}\mathbb{E}_x\{d(x)|N_I \text{\hspace*{-2pt} idle\hspace*{-2pt} vehicles}\}&\hspace*{-2pt}=\hspace*{-2pt}
	\frac{\sqrt{N_{I_0}}}{\sqrt{N_I}}\mathbb{E}_x\{d(x)|N_{I_0} \text{\hspace*{-2pt} idle\hspace*{-2pt} vehicles}\}\hspace*{-2pt}\left(\hspace*{-2pt}1\hspace*{-2pt}+\hspace*{-2pt}\mathcal{O}\hspace*{-2pt}\left(\hspace*{-2pt}\max\{\hspace*{-1pt}N_{I_0}^{-1}\hspace*{-2pt},\hspace*{-1pt}N_{I_0}^{\frac{1}{8}}N_I^{-\frac{9}{8}}\hspace*{-2pt},\hspace*{-2pt}N_I^{-1}\hspace*{-1pt}\}\hspace*{-2pt}\right)\hspace*{-3pt}\right)\hspace*{-2pt}.
	\end{align} 
\end{proposition}

The result implies that the average pickup time $t_p$ is (approximately) inversely proportional to the square root of the number of idle vehicles since $t_p=\frac{\mathbb{E}_x\{d(x)|N_I \text{\hspace*{-2pt} idle\hspace*{-2pt} vehicles}\}}{v}$,  with $v$ being the average traffic speed. The result recalls Mohring's ``square root law''  \cite{mohring1987values} and can be explained intuitively as follows. Consider a square city of unit size with $N_I$ idle vehicles  located in a grid with each idle car equally distant from its four closest neighbors to its left, right, top, or bottom, then the shortest distance between idle cars is equal to $\frac{1}{\sqrt{N_I}}$. The exact proof of Proposition \ref{pickup_theory} for a city of general shape and when the locations of idle vehicles and passengers are random  is deferred to Appendix A.

\begin{remark}
Proposition \ref{pickup_theory} has a few limitations. First, the platform may wait to accumulate idle vehicles and waiting passengers before matching \cite{ke2019optimizing}. This can potentially benefit the passenger/driver as they receive a closer match after waiting for a few more seconds. We do not capture this in Proposition \ref{pickup_theory}.  Second, we assume that both passengers and vehicles are uniformly and independently distributed across the city.  In practice, passengers/drivers may strategically choose their locations to wait for the next vehicle/customer. This is also not considered  in Proposition \ref{pickup_theory}. Nevertheless, we emphasize that our analysis does not require any specific form of function $t_w(\cdot)$. The qualitative results of this paper hold as long as Assumption \ref{assump_pickuptimefunc} is satisfied, and the result of Proposition \ref{pickup_theory} is only used to generate the numerical results.
\end{remark}

In  (\ref{pickup-exact}) we select $N_{I_0}$ as a reference so that   
$\frac{\mathbb{E}_x\{d(x)|N_{I_0} \text{\hspace*{-2pt} idle\hspace*{-2pt} vehicles}\}}{v}$ can be  computed 
 from available TNC  data for a city. For instance, on average, Manhattan has $5089$ TNC vehicles  on the road. 
 Every minute there are $187$ new TNC trips. Each trip takes around $16.3$ minutes \cite{schaller2017empty}, and 
 passengers wait $5$ minutes for pickup \cite{Uber_S1}. In this case  
 we have $N=5089$, $\lambda=187$ trips/min, 
 $ \mu=1/16.3$ min$^{-1}$, and the average pickup time is $\simeq 5$ min. Taking $N_{I_0}=N-\lambda/\mu\simeq 2041$, then $\frac{\mathbb{E}_x \{d(x)|N_{I_0} \text{\hspace*{-2pt} idle\hspace*{-2pt} vehicles}\}}{v}\simeq 5$min, and the pickup time function (\ref{pickup-exact}) becomes:
\begin{align*}
 \mathbb{E}_x\{d(x)|N_I \text{\hspace*{-2pt} idle\hspace*{-2pt} vehicles}\}
 = \frac{226}{\sqrt{N_I}}\left(1+\mathcal{O}\left(\max\{2041^{-1},2.59N_I^{-\frac{9}{8} },N_I^{-1}\}\right)\right),
\end{align*}


The estimate in Proposition \ref{pickup_theory} has approximation error $\mathcal{O}(N_{I_0}^{-1}+(\hspace*{-1pt}1\hspace*{-1pt}+\hspace*{-1pt}(\frac{N_{I_0}}{N_I})^{\hspace*{-1pt}\frac{1}{8}})N_I^{-1})$  for large $N_I$ and $N_{I_0}$. Medium to big sized cities usually have a few thousands of TNC vehicles\footnote{\noindent San Francisco has around  $6000$ active TNC vehicles on average \cite{castiglione2016tncs}, and Manhattan has more than $10,000$ TNC vehicles during peak hours \cite{schaller2017empty}.}, so \eqref{pickup-exact} is a  good approximation to the average pickup time for practical parameter values. In summary, we have:

\begin{corollary}
Assume that all the conditions in Proposition \ref{pickup_theory} hold, and the ride confirmation time $t_m$ is negligible compared to $t_p$, i.e., $t_w=t_m+t_p \simeq t_p$.  Then
	\begin{align}\label{pickup}
 t_w\simeq t_p=\frac{1}{v}\mathbb{E}_x\{d(x)|N_I \text{\hspace*{-2pt} idle\hspace*{-2pt} vehicles}\}\simeq
	\frac{M}{\sqrt{N_I}},
	\end{align}
	where $M=\frac{1}{v}\sqrt{N_{I_0}}\mathbb{E}_x\{d(x)|N_{I_0}^{-1} \text{\hspace*{-2pt} idle\hspace*{-2pt} vehicles}\}$ and $N_I=N-\lambda/\mu$.
\end{corollary}

The ride-hailing platform has a distinctive  supply-side {\em network externality}. As the number of drivers increases, so do the number and spatial density of idle drivers which, in turn, reduces pickup time and increases service quality. This enables larger platforms to offer the same service quality at a lower cost. For instance, consider a small platform and a large platform with the same vehicle occupancy. Assume the small platform is half the size of the large platform in terms of number of vehicles and passengers. Based on (\ref{pickup}), the waiting time for the small platform,  $t_w^{s}$, is $\sqrt{2}$ times that of the large platform, $t_w^{l}$.  Let $t_w^{s}=6$ min, then we have $t_w^{l}=4.2$ min.  The monetary value of this difference is $\alpha*(t_w^s-t_w^l)=\$5.7$ (See Section 4.2 for the value of $\alpha$). This indicates that the smaller platform has to lower the rider fare by $\$5.7$ to attract the same number of passengers of the larger platform.

\begin{remark}
The ride confirmation time $t_m$ is orders of magnitude smaller than the pickup time $t_p$ in large cities in the US. For instance,  New York city has an average of 5089 TNC vehicles, 187 passengers per minute, and each trip takes 16.3 minutes \cite{schaller2017empty}. If  passenger arrivals  are Poisson, then $t_m$ as the average waiting time in the M/G/N queue is sub-second (virtually 0). 
In areas with limited supply of drivers, $t_m$ could be significant and can not be neglected. We can add $t_m$ to the travel cost. We conjecture that in this case if $t_m+t_w$ satisfies Assumption \ref{assump_pickuptimefunc} the qualitative results of the paper still hold.
\end{remark}

Passengers have a reservation cost that summarizes their other travel options: if the TNC travel cost $c$ is greater than the reservation cost, the passenger abandons the TNC for an alternative transport mode.  We assume that the reservation costs of passengers are heterogeneous, and let $F_p(c)$ be the cumulative distribution of reservation costs. The  passenger arrival rate then is  given by
\begin{equation}
\label{demandfunction}
\lambda=\lambda_0\left[ 1-F_p\bigg(\alpha t_w(N-\lambda/\mu)+\beta p_f\bigg) \right]\;\mbox{rides/min},
\end{equation}
where  $\lambda_0$ is the arrival rates of potential passengers total travel demands in the city. Note that  the trip time $t_o$ does not depend on $\lambda$ or $N$, so we absorb it into $F_p$ as a constant. According to (\ref{demandfunction}), the passengers that use the app-based ride-hailing service are all potential passengers except those that leave the platform because its  cost  is greater than their reservation cost.  

{\bf Driver Incentives}: Drivers are sensitive to earnings and respond to the offered wage by joining or leaving the platform.  The average hourly wage of TNC drivers is
\begin{equation}
\label{driverwage}
w= \dfrac{\lambda {p_d}}{N},
\end{equation}
wherein $p_d$ is the per trip payment the driver receives from the platform. The platform keeps the difference between $p_f$ and $p_d$ as its commission or profit.  In 2018 Uber collected \$41B from passengers of which drivers received 78\% corresponding to a 22\% commission rate \cite{Uber_Sec}, and Lyft collected \$8B from passengers and received 26.8\% as commission  \cite{Lyft2019IPO}. 

The average hourly wage  (\ref{driverwage}) is derived as follows.  The total platform payment to all drivers sums  to  $ \lambda {p_d}\; \; \mbox{\$/min}$. Therefore the average hourly wage per driver is   $ \lambda {p_d}\times 60 \; \mbox{\$/hr}$ divided by $N$, where the constant 60 captures the time period of one hour. 

Each driver has a reservation wage. He joins TNC if the platform wage is greater than his reservation wage.   We assume that the reservation wages of drivers are heterogeneous, and denote $F_d(c)$ is the cumulative distribution of reservation wages across the population of drivers.  Hence  
\begin{equation}
\label{supplyfunction}
N=N_0 F_d\left( \dfrac{\lambda p_d}{N}   \right).
\end{equation}
Here $N_0$ is the number of potential drivers (all drivers seeking a job). For ease of notation we drop the constant factor $60$ from the hourly wage formula and absorb it in the function $F_d$ in (\ref{supplyfunction}).  According to (\ref{supplyfunction}), the number of TNC drivers is the  number of potential drivers multiplied by the proportion that joins the platform since their reservation wage is smaller than $w$. 

\begin{remark}
In practice, both supply and demand of a ride-hailing system  vary  within a day. This can be approximated  in a quasi-static analysis by varying $\lambda_0$ and $N_0$ for  peak and off-peak hours.
\end{remark}

\section{TNC decisions in absence of regulation}
\label{higherlevel}
The objective of the app-based ride-hailing platform adapts over time.  In the initial phase it   subsidizes  passengers and drivers to grow the business.  Eventually it shifts to  maximizing the profit. Here we focus on profit maximization assuming that the platform is  unregulated.  Platform subsidy and competition are discussed in Section \ref{extension}.

\subsection{Pricing without Regulation}
The platform rent or profit is 
\begin{equation}
\label{platformrevenue}
\Pi=\lambda (p_f-p_d).
\end{equation}
In a certain period (e.g., each minute),  $\lambda$ trips are completed. Since the platform pockets $p_f-p_d$ from each trip, the total rent in this period is (\ref{platformrevenue}). 

In the absence of regulation, the platform maximizes its rent subject to the   market equilibrium conditions:
\begin{equation}
\label{optimalpricing}
 \hspace{-4cm} \max_{p_f\geq 0, p_d\geq 0} \lambda(p_f-p_d)
\end{equation}
\begin{subnumcases}{\label{constraint_optimapricing}}
\lambda=\lambda_0\left[ 1-F_p\bigg(\alpha t_w(N-\lambda/\mu)+\beta p_f\bigg) \right] \label{demand_constraint}\\
N=N_0 F_d\left( \dfrac{\lambda p_d}{N}   \right) \label{supply_constraint}
\end{subnumcases}
We have the following result on the existence of solution to (\ref{constraint_optimapricing}):
\begin{proposition}
\label{feasibility}
If $F_p(\alpha t_w(N_0))<1$ and $N_0>\lambda_0/\mu$, there exist strictly positive $\lambda, N, p_f$ and $p_d$ that constitute a market equilibrium satisfying (\ref{constraint_optimapricing}).
\end{proposition}

The proof can be found in Appendix B. The assumption $F_p(\alpha t_w(N_0))<1$ means that when the platform recruits all potential drivers $N_0$ and sets the ride price at $0$ ($p_f=0$), there will be a positive number of passengers. This assumption rules out the situation in which  passenger reservation costs are so low and driver reservation wages  are  so high that supply and demand curves do not intersect.

Since (\ref{optimalpricing}) is not a convex problem,  it is not straightforward to determine its  solution. One  approach is  via numerical computation as in \cite{bai2018coordinating}. This is suitable for small problems. Instead, we 
proceed analytically.  We view $\lambda = \lambda(p_f,p_d)$ as a function of $p_f$ and $p_d$ determined implicitly by (\ref{constraint_optimapricing}). The first order necessary conditions of (\ref{optimalpricing}) then simplify to
\begin{subnumcases}{\label{1srordercondition}}
\dfrac{\partial \lambda}{\partial p_f}(p_f-p_d)+\lambda=0 \label{condition1_1}\\
\dfrac{\partial \lambda}{\partial p_d}(p_f-p_d)-\lambda=0 \label{condition2_1}
\end{subnumcases}
in which (\ref{condition1_1}) is equivalent to $\dfrac{\partial \Pi}{\partial p_f}=0$ and (\ref{condition2_1}) is  $\dfrac{\partial \Pi}{\partial p_d}=0$. For non-convex problems, these conditions are   only necessary. However, they are  sufficient in the following case.
\begin{proposition}
\label{1storderconditionsufficient}
Assume that (a) the waiting time function $t_w$ satisfies (\ref{pickup}); (b) the reservation cost and the reservation wage are uniformly distributed as $F_p(c)=\min\{e_p c, 1\}$ and $F_d(w)=\min\{e_d w,1\}$, with $e_p\in \mathbb{R}$ and $e_d \in\mathbb{R}$; (c) the profit maximizing problem (\ref{optimalpricing}) has at least one solution at which the objective value is positive. Then the following equations have a unique solution $(p_f,p_d,\lambda,N)$\footnote{Since $p_f = p_d =\lambda=N=0$ is always a solution, throughout the paper by `solution' we refer to strictly positive solutions unless otherwise stated.}, which is the globally optimal solution to (\ref{optimalpricing}).
\begin{subnumcases}{\label{1srordercondition_proposition1}}
\dfrac{\partial \lambda}{\partial p_f}(p_f-p_d)+\lambda=0  \label{proposition1_constraint1}\\
\dfrac{\partial \lambda}{\partial p_d}(p_f-p_d)-\lambda=0 \label{proposition1_constraint2}\\
\lambda=\lambda_0\left[ 1-F_p\left(\dfrac{\alpha M}{\sqrt{N-\lambda / \mu}}+\beta p_f\right) \right] \label{proposition1_constraint3}\\
N=N_0 F_d\left( \dfrac{\lambda p_d}{N}   \right) \label{proposition1_constraint4}
\end{subnumcases} 
\end{proposition}
The proof of Proposition \ref{1storderconditionsufficient} is deferred to Appendix C. It asserts that (\ref{optimalpricing}) can be effectively computed by finding the unique solution to (\ref{1srordercondition_proposition1}). Note that if the assumptions of Proposition \ref{1storderconditionsufficient}  are not satisfied, we can still solve  (\ref{optimalpricing}) by brute-force computation.

\subsection{Numerical Example}
\label{parametersetup_section}
We present a numerical example and calculate the platform's profit-maximizing decision (\ref{optimalpricing}). To apply Proposition \ref{1storderconditionsufficient}, we assume that the waiting time function $t_w$ satisfies (\ref{pickup}), and that the reservation cost of passengers and the reservation wage of drivers are both  uniformly distributed. Below we specify the model parameters used in the simulation.

{\bf Parameters:} We take the TNC data for the Manhattan Central Business District (CBD) in New York city. It records all trips that started from or ended in Manhattan CBD on a regular weekday. Let $L$ denote the average TNC trip distance. We obtain the following estimates based on \cite{schaller2017empty}:
\begin{equation}
\label{NYsolutions}
N= 5089, \,\lambda=187 \text{ ride}/\text{min}, \,  L = 2.4 \text{ mile}, \, t_o = 16.3 \text{ min}, \, p_f=\$17/\text{trip}, \, p_d=\$10.2/\text{trip}.
\end{equation}
Note that $t_o$ denotes the average TNC trip duration. 

Our estimation proceeds as follows. 
Based on \cite{schaller2017empty}, each day TNC vehicles make  202,262 trips over 91,608 vehicle hours and 802,135 miles. On average, the vehicle are occupied $60\%$  of the time \cite{schaller2017empty}.  Since there are virtually no rides between $1\text{AM}-7\text{AM}$, we  divide the  daily numbers by 18 (hours) to get $N=91,608/18=5089$ and $\lambda=202,262/18/60=187$ ride/min. The average trip length is
\[ \dfrac{\text{total mileage}}{\text{number of trips}}\times \text{occupancy}=\dfrac{802,135}{202,262}*0.6=2.4 \text{ mile}.\]
The average trip duration is \[\dfrac{\text{vehicle hours}}{\text{number of trips}}\times \text{occupancy}=\dfrac{91,608}{202,262}*0.6*60=16.3 \text{ min}.\]
We estimate that  a 16-min, 2.4-mile ride  In Manhattan costs \$17. TNC drivers in New York earn an average of \$22.6 per hour before expenses \cite{parrott2018earning}. This suggests $w=\lambda p_d/N=\$22.6$, and so  $p_d=22.6N/\lambda= \$10.2$ per trip.

Note that our estimates (\ref{NYsolutions}) are  solutions to the profit-maximizing problem (\ref{optimalpricing}). We now utilize these solutions to `reverse-engineer' the model parameters $(N_0,\lambda_0, \alpha, \beta)$. In particular, we select $(N_0, \lambda_0\blue{, } \alpha, \beta)$ so that the solutions to (\ref{optimalpricing}) match the real data (\ref{NYsolutions}). This can be realized by substituting (\ref{NYsolutions}) into (\ref{1srordercondition_proposition1}) and solving the first-order conditions (\ref{1srordercondition_proposition1}). We obtain:
\begin{equation}
\label{parameters_unregulated}
N_0 = 13,512, \,\lambda_0=1512/\text{min}, \, \alpha = 3.2\$/\text{min},\, \beta = 1, \,e_p=0.0262, \,e_d=1, M=226, \mu=1/16.3 \text{min}^{-1}.
\end{equation}
Empirical study suggests that value of travel time (VOT) in the range \$20 to \$100 per hour \cite{schwieterman2018uber} and value of waiting time at 2 to 3 times that of  in-vehicle travel time \cite{quarmby1967choice}. Our estimate of $\alpha = \$3.2 \text{ per min}$ corresponds to a VOT between \$64 and \$96 per hour.  

\begin{figure*}[bt]%
\begin{minipage}[b]{0.32\linewidth}
\centering
\include{figure1_unreg} 
\vspace*{-0.3in}
\caption{Number of drivers under different potential passengers. }
\label{figure1_unreg}
\end{minipage}
\begin{minipage}[b]{0.005\linewidth}
\hfill
\end{minipage}
\begin{minipage}[b]{0.32\linewidth}
\centering
\include{figure2_unreg}
\vspace*{-0.3in}
\caption{Arrival rates of Passengers (per minute).} 
\label{figure2_unreg}
\end{minipage}
\begin{minipage}[b]{0.005\linewidth}
\hfill
\end{minipage}
\begin{minipage}[b]{0.32\linewidth}
\centering
\include{figure3_unreg}
\vspace*{-0.3in}
\caption{Occupancy rate under different potential passengers.}
\label{figure3_unreg}
\end{minipage}
\begin{minipage}[b]{0.32\linewidth}
\centering
\include{figure4_unreg}
\vspace*{-0.3in}
\caption{Per mile ride price and driver payment.} 
\label{figure4_unreg}
\end{minipage}
\begin{minipage}[b]{0.005\linewidth}
\hfill
\end{minipage}
\begin{minipage}[b]{0.32\linewidth}
\centering
\include{figure5_unreg} 
\vspace*{-0.3in}
\caption{Driver wage per hour under different potential passengers.}
\label{figure5_unreg}
\end{minipage}
\begin{minipage}[b]{0.005\linewidth}
\hfill
\end{minipage}
\begin{minipage}[b]{0.32\linewidth}
\centering
\include{figure6_unreg}
\vspace*{-0.3in}
\caption{Passenger travel cost under different potential passengers. }
\label{figure6_unreg}
\end{minipage}
\end{figure*}

{\bf Results:}
We vary $\lambda_0$ between $1000$ and $2000$ to study how the platform decision varies at different times of the day ($\lambda_0$ is large during peak hours). 
The  results are shown in Figures \ref{figure1_unreg}-\ref{figure6_unreg}. As $\lambda_0$ increases, the number of passengers ($\lambda$) and drivers ($N$) both increase. At the same time, occupancy rate  (Figure  \ref{figure3_unreg}), and the ride price increase (Figure \ref{figure4_unreg}). At peak hours, the drivers benefit since they earn more (Figure \ref{figure5_unreg}), but the passengers travel at a higher cost (\ref{costoftravel}) due to the increased trip fare.

Note  that as the number ($\lambda_0$ of potential passengers doubles from 1000 to 2,000 riders per minute, the profit-maximizing fare ($p_f$) per ride set by the platform increases  by 13 percent from \$15.8 to \$17.8 per trip,  driver payment ($p_d$) increases by 16 percent from \$9.4 to \$10.9 per trip, and the platform's share increases by 8 percent from \$6.4 to \$6.9.  The 33 percent increase in driver wages   from \$17.8 to \$25.1 per hour is  due jointly to the  increases in per trip payment to the driver and the vehicle occupancy (from 0.55 to 0.63). By the same token, a driver's hourly wage declines by 33 percent from peak to off-peak hours.  Thus in the absence of a wage floor,  drivers bear most of the risk of shifts in demand.

\section{TNC decisions with wage floor} \label{sec5}
This section is devoted to  platform pricing with a wage floor.  A driver minimum wage $w$ imposes the constraint $\lambda p_d/N\geq w$\footnote{New York City Taxi and Limousine Commission imposes the minimum driver payment \eqref{utilization}.
We assume a constant speed, so per-minute price can be transformed to per-mile price, and we use $p_d$ to represent the sum of the first and second term in \eqref{utilization}. If we neglect the constant shared ride bonus,   (\ref{utilization}) is proportional to $\lambda p_d/N$.}. After a wage floor is imposed the platform may find it prohibitive to hire all drivers who wish to join and, thus, limit the entry of new drivers. We capture this by relaxing \eqref{supply_constraint} to the inequality \eqref{supply_constraint_wage}. The profit maximizing problem subject to a wage floor is
\begin{equation}
\label{optimalpricing_wage}
 \hspace{-4cm} \max_{p_f\geq 0, p_d\geq 0, N} \lambda(p_f-p_d)
\end{equation}
\begin{subnumcases}{\label{constraint_optimapricing_wage}}
\lambda=\lambda_0\left[ 1-F_p\bigg(\alpha t_w(N-\lambda/\mu)+\beta p_f\bigg) \right] \label{demand_constraint_wage}\\
N\leq N_0 F_d\left( \dfrac{\lambda p_d}{N}   \right) \label{supply_constraint_wage} \\
\dfrac{\lambda p_d}{N} \geq w \label{minimum_wage}
\end{subnumcases}
Constraint (\ref{supply_constraint_wage}) indicates that the platform can hire up to the number of all available drivers. This introduces an additional decision variable $N$. It can be solved via numerical computation as in \cite{bai2018coordinating}. Similarly to Proposition \ref{feasibility}, we can show that (\ref{constraint_optimapricing_wage}) has at least one non-trivial solution if $F_p(\alpha t_w(N_0))<1$.

\subsection{A Cheap-Lunch Theorem}
\label{maintheoremsec}
{\bf Example:} 
Consider an example for which we calculate the profit-maximizing prices (\ref{optimalpricing_wage})  for different wage floors $w$. 
We assume that the waiting time function $t_w$ is of form (\ref{pickup}), and that the reservation cost of passengers and the reservation wage of drivers are both  uniformly distributed. We set the model parameters as (\ref{NYsolutions}) and (\ref{parameters_unregulated}).  We emphasize that these assumptions are only needed for numerical simulations. Our analysis does not depend on these assumptions or model parameters. 
\begin{figure*}[bt]%
\begin{minipage}[b]{0.32\linewidth}
\centering
\include{figure1_wage} 
\vspace*{-0.3in}
\caption{Number of drivers under different wage floors. }
\label{figure1_wage}
\end{minipage}
\begin{minipage}[b]{0.005\linewidth}
\hfill
\end{minipage}
\begin{minipage}[b]{0.32\linewidth}
\centering
\include{figure2_wage}
\vspace*{-0.3in}
\caption{Arrival rates of Passengers (per minute).} 
\label{figure2_wage}
\end{minipage}
\begin{minipage}[b]{0.005\linewidth}
\hfill
\end{minipage}
\begin{minipage}[b]{0.32\linewidth}
\centering
\include{figureoccupancy_wage}
\vspace*{-0.3in}
\caption{Occupancy rate under different wage floors.}
\label{figureoccupancy_wage}
\end{minipage}
\begin{minipage}[b]{0.32\linewidth}
\centering
\include{figure5_wage}
\vspace*{-0.3in}
\caption{Per mile ride price and driver payment.} 
\label{figure5_wage}
\end{minipage}
\begin{minipage}[b]{0.005\linewidth}
\hfill
\end{minipage}
\begin{minipage}[b]{0.32\linewidth}
\centering
\include{commission_wage} 
\vspace*{-0.3in}
\caption{Commission rate defined as percentage of service fee in $p_f$. }
\label{commission_wage}
\end{minipage}
\begin{minipage}[b]{0.005\linewidth}
\hfill
\end{minipage}
\begin{minipage}[b]{0.32\linewidth}
\centering
\include{figure3_wage}
\vspace*{-0.3in}
\caption{Driver wage per hour under different wage floors.}
\label{figure3_wage}
\end{minipage}
\begin{minipage}[b]{0.32\linewidth}
\centering
\include{figure4_wage}
\vspace*{-0.3in}
\caption{Passenger pickup time under different wage floor.}
\label{figure4_wage}
\end{minipage}
\begin{minipage}[b]{0.005\linewidth}
\hfill
\end{minipage}
\begin{minipage}[b]{0.32\linewidth}
\centering
\include{GC_wage} 
\vspace*{-0.3in}
\caption{Total cost of passengers under different caps. }
\label{GC_wage}
\end{minipage}
\begin{minipage}[b]{0.005\linewidth}
\hfill
\end{minipage}
\begin{minipage}[b]{0.32\linewidth}
\centering
\include{figure6_wage}
\vspace*{-0.3in}
\caption{Platform profit ($\$$/hour) under different wage floors.}
\label{figure6_wage}
\end{minipage}
\end{figure*}
Figures \ref{figure1_wage}-\ref{figure6_wage} reveal the market response to different levels of the wage floor, including number of drivers, arrival rates of passengers, vehicle occupancy rate, driver wage, passenger pickup time, platform prices, and platform profit. The response has three distinct regimes:
\begin{itemize}
\item $w< \$22.6 $: the wage floor constraint (\ref{minimum_wage}) is inactive and the  solution to  (\ref{optimalpricing_wage}) is the same as that to (\ref{optimalpricing}) because even in the absence of  the minimum wage constraint the platform sets $w = \$ 22.6$ to attract enough drivers.  
\item $\$22.6 \leq w<\$36.3$: both (\ref{supply_constraint_wage}) and (\ref{minimum_wage}) constraints are active. As the minimum wage increases, the platform hires all drivers whose reservation wage is below the minimum wage, the ride cost \eqref{costoftravel} goes down, the quality of service (pickup time) improves, driver wage increases, more passengers are served, and the platform profit reduces. 
\item $w\geq \$36.3 $: only the wage floor constraint (\ref{minimum_wage}) is active. As the minimum wage exceeds $\$36.3$, the platform hires fewer drivers than wish to work,  both ride fare ($p_f$) and pickup time ($t_p$) increase,  fewer passengers take  the ride-hailing option, the drivers who are hired earn more, and the platform profit reduces further. 
\end{itemize}
According to Theorem~\ref{freelunchtheorem} below, the qualitative behavior of most variables, including number of drivers, arrival of passengers, driver wage, travel cost, and platform rent remains consistent with Figures \ref{figure1_wage}-\ref{figure6_wage} for all model parameters. The behavior of $p_f$ may depend on  model parameters.  When $\alpha$ is small,   the trip fare $p_f$ may decrease in the second regime of Figure \ref{figure5_wage}. This is because for small $\alpha$ passengers are more sensitive to trip fare, and  the platform may find it more effective to attract passengers by reducing  the trip fare. However, we emphasize that the total travel cost (\ref{costoftravel}) as the sum of $p_f$ and pickup time always decreases in the second regime.

\begin{remark}
When the wage floor reaches \$39 per hour, the platform profit is 0. In this case, the platform may exit the market. This regulatory risk associated with the ride-sharing business model is explicitly called out by Uber and Lyft in their IPO registration statements \cite{Uber_Sec,Lyft2019IPO}.  In practice, it is unlikely that regulations will drive  platform revenue to zero. New York city's  wage floor of \$27.86  per hour (before expenses) will predictably lead   to 10.5\% decrease in platform profit from \$76K to  \$68K per hour.  
\end{remark}



Assuming the optimum solution to  (\ref{optimalpricing_wage}) is unique,  write it as a function of $w$: $(N^*(w), \lambda^*(w), p_f^*(w), p_d^*(w))$. We have the following theorem.
\begin{theorem}
\label{freelunchtheorem}
Assume that  (\ref{optimalpricing_wage}) has a unique solution.\footnote{For almost all parameter values there will not be multiple solutions with the same optimal value.} For any parameters $(N_0, \lambda_0, \alpha, \beta)$ and any distributions $F_p(\cdot)$ and $F_d(\cdot)$ that satisfy $F_p(\alpha t_w(N_0))<1$ and $N_0>\lambda_0/\mu$, we have $\nabla_+ N^*(\tilde{w})>0$ and $\nabla_+\lambda^*(\tilde{w})>0$,
where $\nabla_+$ denotes the right-hand derivative, and $\tilde{w}$ is the optimal driver wage in absence of regulation, i.e., the  solution to (\ref{optimalpricing}). 
\end{theorem}
The proof of Theorem \ref{freelunchtheorem} can be found in Appendix D.  Theorem \ref{freelunchtheorem} holds for every pickup time $t_w$ that satisfies Assumption \ref{assump_pickuptimefunc} and it does not assume any specific formula for pickup time or a specific matching algorithm utilized by the platform.  This implies that the second regime always exists: when $w\leq \tilde{w}$, the minimum floor constraint (\ref{minimum_wage}) is inactive, so the solution is in regime 1.  When $w=\tilde{w}$, the right-hand derivative of $N$ and $\lambda$ are both strictly positive. This corresponds to the beginning of the second regime, where the platform hires more drivers and serves more passengers. The increase in the number of drivers and passengers implies that the wage of drivers increases and the total cost of passengers decreases.

\textbf{Discussion:} The effect of minimum wage on labor markets has been  the subject of many studies in labor economics since its inception as part of \text{Fair Labor Standard Act} of 1938, and it still remains a contentious topic among economists. We provide a brief overview of the existing literature on the effect of the minimum wage regulation on employment. We then discuss how our result connects to the current literature. 

There is mixed empirical evidence on whether imposing a minimum wage has a positive or negative effect on employment; for instance, the authors in \cite{card1994minimum,card2000minimum} and \cite{neumark2000minimum} use the data from fast-food industry in Pennsylvania and New Jersey and draw drastically different conclusions on the impact  of an increase in the minimum wage. A recent meta-study \cite{belman2014does} found that it is equally likely to find positive or negative employment effect of the minimum wage in the literature; a similar observation is made in \cite{schmitt2013does}.\footnote{We refer the interested reader to \cite{card2015myth,Neumark} for surveys of studies on the effects of minimum wage.} 

A similar division in economic theory literature exists regarding the direction of the employment effect of minimum wage \cite{mcconnell2016contemporary}. 
One strand of work that assumes that the labor market is perfectly competitive concludes that minimum wage has a negative effect on employment. Another strand of work considers a monopsony framework where the employer has bargaining power over the wage, and the labor demand is upward-sloping. This work contends that the minimum wage may actually increase  employment \cite{laing2011labor}. 

Our results above are similar to those of the monopsony framework in the literature. In a ride-sharing market, a TNC has market power over drivers as it explicitly sets price $p_d$ for drivers. Moreover, given the significant size of its drivers	(e.g. Uber is the largest for-profit employer in New York city if we consider drivers as employees \cite{parrott2018earning}), TNC does not face a perfectly competitive labor market.      
	\footnote{We note that in labor economics, by monopsony they do not only refer to the traditional \textit{company town} with a single employer with full market power. The term monopsony applies more broadly to cases where the employer has some market power to set wages and faces upward-sloping labor supply \cite{boal1997monopsony}.}  

Our model is different from the standard economic equilibrium model in which  supply and demand are equal. We consider a model where the supply (drivers $N$) must exceed the demand (riders $\lambda$), and the difference between the supply and demand (idle cars $N_I$) contributes to the total cost the riders faces through waiting times $t_w(N_I)$. Nevertheless, we can use the monopsony framework to provide an intuitive explanation of Theorem \ref{freelunchtheorem} below.

Consider Figure \ref{fig:monopsony-noreg}, where curve $W(L)$ depicts the wage corresponds to every employment level $L$ and curve $MRP$ represents the resulting marginal-revenue product equivalent to employment level $L$. We note that in deriving the MRP curve, we ignore the effect of waiting time $t_w$ and assume that $c=\beta p_f$.
The intersection of $W(L)$ and $MRP$ (point $E$) determines the outcome in a competitive labor market with employment $L^*_c$ and wage $\omega^*_c$. However, a TNC does not face a perfectly competitive labor market, and sets wages to maximize its profit. From $W(L)$ we can determine the marginal cost of labor $MCL$ defined as the marginal cost the TNC has to pay to hire one more driver; we note that to hire an additional driver the TNC has to increase the wage for all of his existing drivers, thus, MLC curve lies above $W(L)$.  Figures \ref{fig:monopsony-noreg}-\ref{fig:monopsony-regime3} depicts the resulting MCL curves for tree different regimes depending on the value of the minimum wage $\omega_m$. The optimal employment level and wage can be determined by the intersection of $MRP$ and $MCL$ curves, i.e. point $A$ in the first regime, point $B$ in the second regime, and point $C$ in the third regime.

As the minimum wage $\omega_m$ increases, it is easy to verify that the number of drivers is constant in the first regime (Figure \ref{fig:monopsony-noreg}), increases in the second regime (Figure \ref{fig:monopsony-regime2}), and decreases in the third regime (Figure \ref{fig:monopsony-regime3}). Ignoring the effect of idle vehicles and waiting time $t_w$ on cost $c$, the number of riders follow a similar pattern as the number of drivers. Consequently, the cost for riders is constant in the first regime, decreases in the second regime, and increases in the third regime. 
	The above monopsony argument does not capture the presence of idle vehicles and their effect on waiting time $t_w$. The results of Theorem \ref{freelunchtheorem} establishes the result formally incorporating the impact of $t_w$ on total cost $c$ to riders.

\begin{figure*}[bt]%
	\begin{minipage}[b]{0.32\linewidth}
		\centering
		\include{Monopsony-Fig1}
		\vspace*{-0.3in}
		\caption{Regulated market with wage floor $w_m< \omega^*_1$ (first regime)}
		\label{fig:monopsony-noreg}
	\end{minipage}
	\begin{minipage}[b]{0.005\linewidth}
		\hfill
	\end{minipage}
	\begin{minipage}[b]{0.32\linewidth}
		\centering
		\include{Monopsony-Fig2}
		\vspace*{-0.3in}
		\caption{Regulated market with wage floor $\omega^*_1\leq w_m\leq \omega^*_c$ (second regime) }
		\label{fig:monopsony-regime2}
	\end{minipage}
	\begin{minipage}[b]{0.005\linewidth}
		\hfill
	\end{minipage}		
	\begin{minipage}[b]{0.32\linewidth}
		\centering
		\include{Monopsony-Fig3}
		\vspace*{-0.3in}
		\caption{Regulated market with wage floor $w_m> \omega^*_c$ (third regime) }
		\label{fig:monopsony-regime3}
	\end{minipage}
	\begin{minipage}[b]{0.005\linewidth}
		\hfill
	\end{minipage}		
\end{figure*}

\section{TNC decisions with cap on number of drivers}
\label{sec6}
Let $N_{cap}$ be the cap on the total number of drivers. With a cap constraint, the platform  pricing problem is
\begin{equation}
\label{optimalpricing_cap}
 \hspace{-4cm} \max_{p_f\geq 0, p_d\geq 0} \lambda(p_f-p_d)
\end{equation}
\begin{subnumcases}{\label{constraint_optimapricing_cap}}
\lambda=\lambda_0\left[ 1-F_p\bigg(\alpha t_w(N-\lambda/\mu)+\beta p_f\bigg) \right] \label{demand_constraint_cap}\\
N= N_0 F_d\left( \dfrac{\lambda p_d}{N}   \right) \label{supply_constraint_cap} \\
N\leq N_{cap} \label{max_cap}
\end{subnumcases}
It is unnecessary to  relax (\ref{supply_constraint_cap}), since the platform can always lower $p_d$ to increase its profit. As with Proposition \ref{feasibility} we can  show that (\ref{constraint_optimapricing_cap}) admits a non-trivial solution if $F_p(\alpha t_w(N_0))<1$.
This is a non-convex program. One approach is via numerical computation as in \cite{bai2018coordinating}. This is suitable for a small problem. Alternatively, we can find the optimal solution based on first order conditions for the following  special case:
\begin{proposition}
\label{1storderconditionsufficient_cap}
Assume that (a) the waiting time function $t_w$ satisfies (\ref{pickup}); (b) the reservation cost and the reservation wage are uniformly distributed as $F_p(c)=\min\{e_p c, 1\}$ and $F_d(w)=\min\{e_d w,1\}$, with $e_p\in \mathbb{R}$ and $e_d \in\mathbb{R}$; (c) the profit maximizing problem (\ref{optimalpricing_cap}) has at least one solution at which the objective value is positive.  Then the first order conditions of (\ref{optimalpricing_cap}) admit a unique solution $(p_f,p_d,\lambda,N)$, which is the globally optimal solution to (\ref{optimalpricing_cap}).
\end{proposition}
The proof is deferred to Appendix E, and the first order conditions of (\ref{optimalpricing_cap}) are defined in the proof.
\begin{figure*}[bt]%
\begin{minipage}[b]{0.32\linewidth}
\centering
\include{figure1_cap} 
\vspace*{-0.3in}
\caption{Number of drivers under different caps. }
\label{figure1_cap}
\end{minipage}
\begin{minipage}[b]{0.005\linewidth}
\hfill
\end{minipage}
\begin{minipage}[b]{0.32\linewidth}
\centering
\include{figure2_cap}
\vspace*{-0.3in}
\caption{Arrival rates of Passengers (per minute).} 
\label{figure2_cap}
\end{minipage}
\begin{minipage}[b]{0.005\linewidth}
\hfill
\end{minipage}
\begin{minipage}[b]{0.32\linewidth}
\centering
\include{figure3_cap}
\vspace*{-0.3in}
\caption{Occupancy rate under different caps.}
\label{figure3_cap}
\end{minipage}
\begin{minipage}[b]{0.32\linewidth}
\centering
\include{figure4_cap} 
\vspace*{-0.3in}
\caption{Per mile ride price and driver payment under different caps.}
\label{figure4_cap}
\end{minipage}
\begin{minipage}[b]{0.005\linewidth}
\hfill
\end{minipage}
\begin{minipage}[b]{0.32\linewidth}
\centering
\include{commission_cap} 
\vspace*{-0.3in}
\caption{Commission rate defined as percentage of service fee in $p_f$.}
\label{commission_cap}
\end{minipage}
\begin{minipage}[b]{0.005\linewidth}
\hfill
\end{minipage}
\begin{minipage}[b]{0.32\linewidth}
\centering
\include{figure6_cap}
\vspace*{-0.3in}
\caption{Driver wage per hour under different caps.}
\label{figure6_cap}
\end{minipage}
\begin{minipage}[b]{0.32\linewidth}
\centering
\include{figure5_cap}
\vspace*{-0.3in}
\caption{Passenger pickup time under different caps.} 
\label{figure5_cap}
\end{minipage}
\begin{minipage}[b]{0.005\linewidth}
\hfill
\end{minipage}
\begin{minipage}[b]{0.32\linewidth}
\centering
\include{GC_cap}
\vspace*{-0.3in}
\caption{total cost of passengers under different caps.} 
\label{GC_cap}
\end{minipage}
\begin{minipage}[b]{0.005\linewidth}
\hfill
\end{minipage}
\begin{minipage}[b]{0.32\linewidth}
\centering
\include{figure7_cap} 
\vspace*{-0.3in}
\caption{Platform profit ($\$$/hour) under different caps.}
\label{figure7_cap}
\end{minipage}
\end{figure*}
{\bf Example:} 
Consider an example where the platform solves the profit-maximizing problem (\ref{optimalpricing_cap}) for different levels of $N_{cap}$. We assume that the waiting time function $t_w$ is of form (\ref{pickup}), and that the reservation cost of passengers and the reservation wage of drivers are both  uniformly distributed. We set the model parameters as (\ref{NYsolutions}) and (\ref{parameters_unregulated}).

Figures \ref{figure1_cap}-\ref{figure7_cap} exhibit the market response to different caps on the total number of vehicles. These responses include the arrival rates of passengers, occupancy rate, platform prices, driver wage, pickup time and platform profit.  It is more instructive to ``read'' the figures from right to left, as the cap decreases.   
As the cap decreases, the supply of vehicles drops (Figure \ref{figure1_cap}), so it is more difficult for passengers to find a ride. In this case,  pickup time increases (Figure \ref{figure5_cap}), and the number of rides decreases (Figure \ref{figure2_cap}).  Here are some interesting observations:
\begin{itemize}
\item The platform loses passengers faster than  it loses drivers. This is evidenced by the drop in occupancy (Figure \ref{figure3_cap}). 
\item The pickup time increases at an increasing rate (Figure \ref{figure5_cap}).  This is just a counterpart of the aforementioned  network externality.
\item Both trip fare and driver wage drop (Figures \ref{figure4_cap}, \ref{figure6_cap}).
\end{itemize}
These observations can be explained.  As the cap reduces the number of drivers, the passenger pickup time increases. Since $t_p$ is a decreasing convex function of $N$, it has an increasing derivative as $N$ decreases. Therefore customers leave the platform at an increasing rate as $N_{cap}$ decreases. This rate is greater than the decreasing rate of $N_{cap}$, so occupancy rate decreases. In this case,  the platform loses passengers quickly, and has to reduce trip prices to keep passengers from leaving. This further squeezes driver pay (Figure \ref{figure6_cap}). 

A surprising fact is that the cap on number of drivers hurts the earning of drivers (Figure \ref{figure6_cap}). This is contrary to the common belief that limiting their number will protect drivers, as expressed in the regulation freezing the number of TNC vehicles in New York  for one year \cite{cnbc_NY}.  This happens because the platform hires  drivers with lowest reservation wage first.
That is, with a smaller cap on the number of drivers, the platform responds by reducing driver pay and hiring drivers  with lower reservation costs. Thus the benefit of limiting supply is intercepted by the platform. This is in contrast with the situation of taxis that need a  medallion to operate.  A limit on the number of medallions will increase their value
and benefit their owners\footnote{The platform revenue (sum of platform rent and  driver payments) divided by the number of drivers increases as the cap decreases.}, who may be taxi drivers.  In the TNC case, the platform accumulates the increased value.
This conclusion  holds in general and is not affected by the model parameters.
\begin{figure*}[bt]%
\begin{minipage}[b]{0.32\linewidth}
\centering
\include{figure1_tax} 
\vspace*{-0.3in}
\caption{Number of drivers under different surcharge. }
\label{figure_tax1}
\end{minipage}
\begin{minipage}[b]{0.005\linewidth}
\hfill
\end{minipage}
\begin{minipage}[b]{0.32\linewidth}
\centering
\include{figure2_tax}
\vspace*{-0.3in}
\caption{Arrival rates of Passengers (per minute).} 
\label{figure_tax2}
\end{minipage}
\begin{minipage}[b]{0.005\linewidth}
\hfill
\end{minipage}
\begin{minipage}[b]{0.32\linewidth}
\centering
\include{figure3_tax} 
\vspace*{-0.3in}
\caption{Per mile ride price and driver payment under different surcharge.}
\label{figure_tax3}
\end{minipage}
\begin{minipage}[b]{0.32\linewidth}
\centering
\include{figure4_tax}
\vspace*{-0.3in}
\caption{Driver wage per hour under different surcharge.}
\label{figure_tax4}
\end{minipage}
\begin{minipage}[b]{0.005\linewidth}
\hfill
\end{minipage}
\begin{minipage}[b]{0.32\linewidth}
\centering
\include{figure5_tax}
\vspace*{-0.3in}
\caption{Total cost of passengers under different surcharge.} 
\label{figure_tax5}
\end{minipage}
\begin{minipage}[b]{0.005\linewidth}
\hfill
\end{minipage}
\begin{minipage}[b]{0.32\linewidth}
\centering
\include{figure6_tax} 
\vspace*{-0.3in}
\caption{Platform profit ($\$$/hour) under different surcharge.}
\label{figure_tax6}
\end{minipage}
\end{figure*}
\section{TNC decisions with congestion surcharge} \label{sec7}
As of Jan 2019, all trips by \textit{for-hire} vehicles that cross 96th street in NYC incur a congestion surcharge of \$2.75 per TNC trip. We model the likely impact of this policy by adding a congestion surcharge $p_c$ to the travel cost (\ref{costoftravel}).  The profit-maximizing problem for the platform now is
\begin{equation}
\label{optimalpricing_congestion}
 \hspace{-5cm} \max_{p_f\geq 0, p_d\geq 0} \lambda(p_f-p_d)
\end{equation}
\begin{subnumcases}{\label{constraint_optimapricing_congestion}}
\lambda=\lambda_0\left[ 1-F_p\bigg(\alpha t_w(N-\lambda/\mu)+\beta p_f+\beta p_c\bigg) \right] \label{demand_constraint_congestion}\\
N= N_0 F_d\left( \dfrac{\lambda p_d}{N}   \right) \label{supply_constraint_congestion} 
\end{subnumcases}
As with Proposition \ref{feasibility}, we can show that the constraint set (\ref{constraint_optimapricing_congestion}) is non-empty if $F_p\bigg(\alpha t_w(N_0)+\beta p_c\bigg)<1$. Since (\ref{optimalpricing_congestion}) is not a convex program, one approach to solve (\ref{optimalpricing_congestion}) is via brute-force computation \cite{bai2018coordinating}. Alternatively, we can show that the first order conditions are sufficient for global optimization for some special cases:
\begin{proposition}
\label{1storderconditionsufficient2}
Assume that (a) the waiting time function $t_w$ satisfies (\ref{pickup}); (b) the reservation cost and the reservation wage are uniformly distributed as $F_p(c)=\min\{e_p c, 1\}$ and $F_d(w)=\min\{e_d w,1\}$, with $e_p\in \mathbb{R}$ and $e_d \in\mathbb{R}$; (c) the profit maximizing problem (\ref{optimalpricing_congestion}) has at least one solution at which the objective value is positive.  Then the following equations have a unique solution $(p_f,p_d,\lambda,N)$, which is the globally optimal solution to (\ref{optimalpricing_congestion}).
\begin{subnumcases}{\label{1srordercondition_proposition2}}
\dfrac{\partial \lambda}{\partial p_f}(p_f-p_d)+\lambda=0  \label{proposition1_constraint12}\\
\dfrac{\partial \lambda}{\partial p_d}(p_f-p_d)-\lambda=0 \label{proposition1_constraint22}\\
\lambda=\lambda_0\left[ 1-F_p\left(\dfrac{\alpha M}{\sqrt{N-\lambda / \mu}}+\beta p_f+\beta p_c\right) \right] \label{proposition1_constraint32}\\
N=N_0 F_d\left( \dfrac{\lambda p_d}{N}   \right) \label{proposition1_constraint42}
\end{subnumcases} 
\end{proposition}
The proof is similar to that for Proposition \ref{1storderconditionsufficient}, and is therefore omitted. 

We estimate the platform's response to various values of congestion surcharge $p_c$ by numerical simulation. In this example, we impose the same assumptions and model parameters as in Section \ref{parametersetup_section}. Simulation results, presented in Figure \ref{figure_tax1}-\ref{figure_tax6}, show that under a congestion surcharge of \$2.75 per trip, 
the number of TNC vehicles drops by $11.9\%$ from 5089 to 4480,  TNC rides reduce by $17.1\%$ from $187$/min to $155$/min, and  platform revenue shrinks by $37.3\%$ from $\$76,006$/hour to \$47,686 per hour. This also suggests that the TNC business model is  vulnerable to regulatory risk.

\section{Extensions}
\label{extension}
We formulate some extensions of the basic model to examine platform subsidy, platform competition and autonomous mobility on demand.

\subsection{Platform Subsidy}
The ride-hailing platform company is not always a short-term profit maximizer. In its early stages, it tries to grow its business  via subsidies to both passengers and drivers. To model this, we consider a ride-hailing platform that sets   prices  to maximize the number of rides or passengers $\lambda$ subject to a reservation revenue $R$, which may be positive or negative (negative $R$ indicates subsidy):
\begin{equation}
\label{optimalpricing_subsidy}
 \hspace{-5cm} \max_{p_f\geq 0, p_d\geq 0} \lambda
\end{equation}
\begin{subnumcases}{\label{constraint_optimapricing_subsidy}}
\lambda=\lambda_0\left[ 1-F_p\bigg(\alpha t_w(N-\lambda/\mu)+\beta p_f\bigg) \right] \label{demand_constraint_subsidy}\\
N= N_0 F_d\left( \dfrac{\lambda p_d}{N}   \right) \label{supply_constraint_subsidy} \\
\lambda(p_f-p_d)\geq R \label{budget_constraint}
\end{subnumcases}
For notational convenience, let $(\lambda^\star, p_f^\star, p_d^\star)$ be the  solution to (\ref{optimalpricing_subsidy}), and denote $(\tilde{\lambda}, \tilde{p}_f, \tilde{p}_d)$ as the  solution to the non-subsidy case (\ref{optimalpricing}). 

We define subsidy as $\epsilon_f=p_f^\star-\tilde{p}_f$ and $\epsilon_d=\tilde{p}_d-p_d^\star$, where $\epsilon_f$ and $\epsilon_d$ represent the subsidy to passengers and drivers, respectively. Note that this definition essentially compares $(p_f^\star, p_d^\star)$ to the profit-maximizing prices.  
For ease of understanding, we define $B=\tilde{\lambda}(\tilde{p}_f-\tilde{p}_d)-R$ as the  subsidy budget. When reservation revenue is the maximal profit, i.e., $R=\tilde{\lambda}(\tilde{p}_f-\tilde{p}_d)$, the subsidy budget  is $0$, and $\epsilon_f=\epsilon_d=0$.

We estimate the platform's ridership under different levels of subsidy. In the numerical example we impose the same assumptions and model parameters as in Section \ref{parametersetup_section}. Simulation results presented in Figure \ref{figure_subsidy1}-\ref{figure_subsidy3} show that the platform should always subsidize both sides of the market, regardless of the subsidy level. Another interesting observation is that the platform should subsidize drivers more than it does passengers. This conclusion, however, depends on the elasticities of demand and supply because the platform has to grow both sides of the market  to maximize  profit. Under a fixed budget, the platform  allocates more subsidy to the less price-sensitive side  as  it costs more to grow this side of the market by one unit. 
\begin{figure*}[bt]%
\begin{minipage}[b]{0.32\linewidth}
\centering
\include{figure_subsidy1}
\vspace*{-0.3in}
\caption{Subsidies to passengers and drivers under different subsidy budgets.} 
\label{figure_subsidy1}
\end{minipage}
\begin{minipage}[b]{0.005\linewidth}
\hfill
\end{minipage}
\begin{minipage}[b]{0.32\linewidth}
\centering
\include{figure_subsidy2}
\vspace*{-0.3in}
\caption{Arrival rate of passengers under different subsidy budgets.}
\label{figure_subsidy2}
\end{minipage}
\begin{minipage}[b]{0.005\linewidth}
\hfill
\end{minipage}
\begin{minipage}[b]{0.32\linewidth}
\centering
\include{figure_subsidy3}
\vspace*{-0.3in}
\caption{Number of drivers under different subsidy budgets.} 
\label{figure_subsidy3}
\end{minipage}
\end{figure*}

\subsection{Platform Competition}
Consider two platforms (e.g., Uber and Lyft) competing with each other to maximize their  profits. The profits are coupled through the market response to the joint decisions of both platforms: passengers choose the platform with lower overall cost, and drivers work for the platform with a higher wage rate. This subsection modifies the model to capture this competition.

Each platform selects its passenger fare and driver wage.  Passengers and drivers respond to the platform prices until the market settles down. Assume that when the market settles down, both platforms survive with a positive profit. In this case neither passengers nor drivers deviate from their choice of platform at the market equilibrium, so the passenger costs and driver wages for the two platforms are equal. This gives rise to the following profit maximization  problem for one platform, given the pricing decisions $(p_f',p_d')$ of its competitor:
\begin{equation}
\label{optimalpricing_competition}
 \hspace{-5cm} \max_{p_f\geq 0, p_d\geq 0} \lambda(p_f-p_d)
\end{equation}
\begin{subnumcases}{\label{constraint_optimapricing_competition}}
 \alpha t_w(N-\lambda/\mu)+\beta p_f=  \alpha t_w(N'-\lambda'/\mu)+\beta p_f' \label{demand_equal_competition}\\
\dfrac{\lambda p_d}{N}  = \dfrac{\lambda' p_d'}{N'}  \label{supply_equal_competition}\\
\lambda+\lambda'=\lambda_0\left[ 1-F_p\bigg( \alpha t_w(N-\lambda/\mu)+\beta p_f\bigg) \right]
 \label{demand_constraint_competition}\\
N+N'= N_0 F_d\left( \dfrac{\lambda p_d}{N}   \right) \label{supply_constraint_competition}
\end{subnumcases}
Constraints (\ref{demand_equal_competition}) and (\ref{supply_equal_competition}) guarantee that if both platforms have positive number of passengers and drivers, then the passenger cost and driver wage in the two platforms are the same, so  no passenger or driver has an incentive to switch platforms. Note that the market outcomes $(N,\lambda, N', \lambda')$ are not given. Instead, they are governed by the market equilibrium conditions (\ref{demand_equal_competition})-(\ref{supply_constraint_competition}) and the platform prices $(p_f,p_d,p_f',p_d')$.  One  difference between (\ref{optimalpricing_competition})-(\ref{supply_constraint_competition})
and the monopoly case (\ref{optimalpricing})-(\ref{supply_constraint}) is that in the former the waiting time for each TNC depends on the number of its own idle vehicles rather than on the sum of the idle vehicles of the two platforms.  The second difference, by contrast, is that the wage rate (\ref{supply_constraint_competition}) is determined by the sum of the demand for drivers by both platforms.

Analogously, the second platform's decisions $(p_f',p_d')$
will maximize its own profits given the decisions $(p_f,p_d)$ of the first.  The solution of the two
decision problems will be a Nash equilibrium.

Due to non-convexity of (\ref{optimalpricing_competition}), the question of existence and uniqueness of Nash equilibrium remains open. It is possible that the two platforms will split the heterogeneous passengers, with one platform offering a higher fare, lower waiting time, luxury rides to passengers with  higher reservation cost; the emergence of such equilibrium outcomes with product quality differentiation was first demonstrated by \cite{shaked1982relaxing,shaked1983natural} in oligopolies.

\subsection{Autonomous Vehicles}
Autonomous vehicles (AV) will revolutionize road transportation \cite{talebpour2016influence}. AV companies claim they will banish 94 percent of all accidents attributed to human error \cite{waymo_safety}.  So commuters can sit back and relax, work, or entertain themselves.   Eventually there will be hardly any need for human drivers. The impact on our lives will be profound. Uber and Lyft have R\&D efforts to build self-driving ride-hailing vehicles.   Billions of venture capital are flowing into the race to develop AVs.

We  model   a company that owns and operates a fleet of autonomous vehicles to provide autonomous mobility on demand (AMoD) service \cite{spieser2014toward, zhang2018routing}, 
and  compare it  to a ride-hailing service with human drivers.
We modify the model (\ref{optimalpricing})-(\ref{supply_constraint}) to relate the decisions of an
AMoD monopoly and those of a TNC monopoly.  The AMoD monopoly will set its ride rates to
maximize its profit (\ref{amodprofit}) subject to demand (\ref{amoddemand}):
\begin{align} 
 & \max_{p_f\geq 0, N\geq 0} \lambda(p_f-Nc_{av}) \label{amodprofit}\\
 & \text{s.t. }\lambda=\lambda_0\left[ 1-F_p\bigg(\alpha t_w(N-\lambda/\mu)+\beta p_f\bigg) \right] \label{amoddemand}
\end{align}
Here $N$ is the number of deployed  AVs and $c_{av}$ is the per ride investment and operating cost of an AV.  Comparing this with the TNC decision making model (\ref{optimalpricing})-(\ref{supply_constraint}) we get a straightforward formal identification: 
\[
\text{number of AVs deployed $\doteq$ number of drivers hired, and wage rate $w \doteq c_{av}$ vehicle cost} .\]  
Following the NYTLC ruling, we take $w = \$27.86$ per hour or an annual  cost of \$55,000 for 2,000 hours per year of driver (plus vehicle) service.  So for the AMoD monopoly to be as profitable as the TNC monopoly (\ref{amodprofit})
implies that an AV's annual investment and operating cost should be smaller than \$55,000.  How likely is this?

Today's AVs do not meet this cost target.  In records submitted to the California Department of Motor Vehicles (DMV) Waymo  reported that in its 2017 AV tests its safety drivers disengaged autonomous driving once every 5,500 miles \cite{waymo_safety}. Waymo reports a disengagement when its evaluation process identifies the event as having `safety significance', so this  rate is almost 100 times worse than the estimated 500,000 miles per accident in 2015 for human drivers.  With this poor safety performance, each AV will require a safety driver, making its total cost more than twice  today's TNC cost.  
Of course AV safety will improve over time with more and more testing and R\&D but it's anyone's guess as to when an AV will perform as safely as human drivers.  

Alternatively, AMoD service can be scaled back to very controlled environments that reduce the risk of accidents by a factor of 100.  That  direction is also being pursued.  For example,  Waymo is providing  rides to 400 people in the calm, sunny  suburb of Chandler, AZ \cite{rueter2018}. These AMoD rides use AVs with a safety driver.  

One additional piece of evidence also suggests that the cost of AVs is very high.  Two proposed contracts  show the leasing cost of AV cars and shuttles  of well over \$100,000 each per year \cite{avcost}.  EasyMile is charging more than \$27,000 a month per small electric shuttle for cities that sign up for one year of service.  Drive.ai charges \$14,000 monthly per vehicle for one year. Considering that a TNC driver (with car) costs \$55,000 per year or \$4,400 per month, it seems unlikely that these are viable business models, except in selective subsidized niche markets.\footnote{According to \cite{avcost}, ``Arlington, Texas, a suburb of Dallas, hired Drive.ai to run three on-demand self-driving shuttles in the entertainment district. For the yearlong program, the city will foot 20\% of the \$435,000 price tag and a federal grant will cover the rest.''}

\section{Conclusion}
\label{section_conclusion}
This paper analyzed the impact of three regulations on the ride-hailing app-based platforms or TNCs  like Uber and Lyft: (a) a floor under driver wage; (b) a cap on total number of drivers; and (c) a per-trip congestion tax that goes to the public transit. We constructed a general equilibrium model to predict  market responses to the platform's decision on fares and wages, with and without these regulations. We showed that imposing a wage floor  increases driver employment, lowers pickup time, decreases ride cost, and attracts more passengers, over a wide range of parameters. Our analysis suggests that a higher minimum wage benefits both drivers and passengers, at the expense of platform profits.  On the other hand,  a cap on the number of drivers or vehicles  hurts  drivers, as the platform benefits by hiring cheaper drivers when  supply is limited. Variants of our model were analyzed from other perspectives as well, including platform subsidy, platform competition and autonomous vehicles.

Our study advocates a wage floor for TNC drivers. Our simulation shows that increasing driver wage by 23.3\% (from \$22.6  to \$27.86 per hour before expense) will increase the number of TNC vehicles by 23.3\% (from 5089 to 6276), increase TNC ridership by 24.6\% (from 187 to 233 per min), improve the pickup time by 10\% (from 5 min to 4.5 min), decrease the travel cost by 3.6\% (from \$33.4 to \$32.2 per trip), and reduce the platform rent by 10.5\% (from \$76K to \$68K per hour). This indicates that enforcing a minimum wage for drivers benefits both drivers and passengers.  Under the wage floor, the platform is motivated hire more drivers to attract more passengers so as to increase the platform sales. As a consequence, more drivers are hired, more passengers are delivered, each driver earns more, and each passenger spends less. The wage floor squeezes the monopoly profit of the platform and improves the efficiency of the  system. It thus boosts the TNC economy without costing taxpayer money.

A congestion surcharge will relieve traffic congestion by reducing the number of TNC vehicles.  Numerical simulation suggests that  a  surcharge of \$2.75 per trip in New York will reduce   
TNC vehicles by 11.9\% (from 5089 to 4480) and   TNC rides by 17.1\% (from 187 to 155 per min). More importantly, the funds raised from this surcharge can be used to subsidize public transit. In New York, it is estimated that the congestion surcharge will yield  \$1M  per day. The money goes to the Metropolitan Transportation Authority to upgrade the subway system. It can be used  to make public transportation  cleaner, faster, and safer, so that more residents will commute by transit. Increased public transit ridership will improve the efficiency of the city's transportation system, and reduce the environmental footprint of transportation, which accounts for 28\% percent of the total carbon emissions in the US.

Interest in TNC regulation has been driven by concerns about working conditions of TNC drivers and by the deleterious impact
on urban transport of TNC growth. This paper deals only with the impact on driver wage and ride fare.  There is a debate whether TNC drivers are more like `independent contractors' or more like employees \cite{Tirole_regulating}. This paper contributes to that debate in showing that TNC driver wages can be significantly increased and passenger fares decreased at the cost of lower TNC profits.

\section*{Acknowledgments} This research was supported by National Science Foundation EAGER
award 1839843.

\bibliographystyle{unsrt}
\bibliography{resourceprocurement}


\section*{Appendix}

\subsection*{\bf{A: Proof of Proposition \ref{pickup_theory}}}

\begin{proof}
	We prove the result in two steps. (i) First, we consider the  case of a passenger located at the origin and $N_I=N-\lambda/\mu$ idle vehicles  uniformly and independently distributed in a disk of radius $R$ centered at the origin. We show that the expected distance of the passenger to the closest idle vehicle is $\sqrt{N_{I_0}}\mathbb{E}\{d(N_{I_0})\}\frac{1}{\sqrt{N_I}}(1+\mathcal{O}(\max\{N_{I_0}^{-1},N_I^{-1}\}))$. (ii) Second, based on the result of part (i), we show that for a city with any two-dimensional area $\mathcal{A}$, the expected shortest distance to an idle vehicle of a passenger is also given by $\sqrt{N_{I_0}}\mathbb{E}\{d(N_{I_0})\}\frac{1}{\sqrt{N_I}}(1+\mathcal{O}(\max\{N_{I_0}^{-1},N_I^{-1}\}))$.

	\textbf{(i)} To prove the first step, let $d(n):=\min (|x_1|,\cdots,|x_{n}|)$ be the shortest distance to the origin among $n$ idle vehicles where $x_i \in R^2$ is the location of the $i$th idle vehicle. Then the cumulative distribution function (cdf) of $d(n)$ is 
	\begin{align}
	\mathbb{P}\{d(n)\leq r\}=1-\mathbb{P}\{d(n)>r\}=1-\mathbb{P}\{|x_i|>r,\;\;\forall i\}=1-(1-\frac{\pi r^2}{\pi R^2})^{n}\nonumber
	\end{align} 
	Therefore, the probability density function (pdf) of $d$ is 
	\begin{align}
	f_{d(n)}(r)=n(1-\frac{r^2}{R^2})^{n-1}\frac{2 r}{R^2}\nonumber
	\end{align}
	Consequently, 
	\begin{align*}
	\mathbb{E}\{d(n)\}&=\int_0^R f_{d(n)}(r)\times r \times dr=\int_0^R 2n(1-\frac{r^2}{R^2})^{n-1} \frac{r^2}{R^2}dr\\
	&= \frac{2}{3}n(1-\frac{r^2}{R^2})^{n-1} \frac{r^3}{R^2}\Big|_{0}^R+\int_0^R \frac{2^2}{3}n(n-1)(1-\frac{r^2}{R^2})^{n-2} \frac{r^4}{R^4}dr\\
	&=\int_0^R \frac{2^2}{3}n(n-1)(1-\frac{r^2}{R^2})^{n-2} \frac{r^4}{R^4}dr\\
	\cdots\\
	&= 2^{n}\frac{\prod_{i=0}^{n-1}(n-i)}{\prod_{i=1}^{n}(2i+1)}\frac{r^{2n+1}}{R^{2n}}\Big|_{0}^R\\
	&=\frac{(2^{n}n!)^2}{(2n+1)!}R
	\end{align*}
	in which the fifth equality follows from an iterative application of integration by parts (similar to the one that leads to the third equality). Therefore, (for any $N_{I_0}$)
	\begin{align*}
	\mathbb{E}\{d(N_I)\}=
	\begin{cases}
	\mathbb{E}\{d(N_{I_0})\}\prod_{i=N_{I_0}}^{N_I-1}\left(\frac{2i}{2i+1}\right) & \text{if  }  N_I> N_{I_0}\\
	\mathbb{E}\{d(N_{I_0})\}\frac{1}{\prod_{i=N_I}^{N_{I_0}-1}\left(\frac{2i}{2i+1}\right)} & \text{if  } N_I< N_{I_0}
	\end{cases}
	\end{align*}

	Next, we prove that $\prod_{i=m}^{n}\left(\frac{2i}{2i+1}\right)=\frac{\sqrt{m}}{\sqrt{n}}(1 + \mathcal{O}(m))$ for $m<n$. We have,
	
	\begin{align*}
	\ln \left(\prod_{i=m}^{n}\left(\frac{2i}{2i+1}\right)\right)=\sum_{i=m}^{n}\ln\left(1-\frac{1}{2i+1}\right).
	\end{align*}
	For  $x\in [0,\sim 0.43]$, we have $-x-x^2\leq \ln(1-x)\leq -x$. Thus,
	\begin{gather*}
	-\frac{1}{2}\sum_{i=m}^{n}\frac{1}{i+\frac{1}{2}}-\sum_{i=m}^{n}\frac{1}{(2i+1)^2}\leq \sum_{i=m}^{n}\ln\left(1-\frac{1}{2i+1}\right) \leq -\sum_{i=m}^{n}\frac{1}{2i+1}\\
	\Longleftrightarrow\\
	-\frac{1}{2}\sum_{i=m}^{n}\frac{1}{i+\frac{1}{2}}-\frac{1}{2}\sum_{i=m}^{n}\left[\frac{1}{i-\frac{1}{2}}-\frac{1}{i+\frac{1}{2}}\right]\leq \sum_{i=m}^{n}\ln\left(1-\frac{1}{2i+1}\right) \leq -\frac{1}{2}\sum_{i=m}^{n}\frac{1}{2i+1}\\
	\Longrightarrow\\
	-\frac{1}{2}\ln\left(\frac{n}{m}\right)+\mathcal{O}(m^{-1})-\frac{1}{2}\left[\frac{1}{m-\frac{1}{2}}-\frac{1}{n+\frac{1}{2}}\right]\leq \sum_{i=m}^{n}\ln\left(1-\frac{1}{2i+1}\right) \leq -\frac{1}{2}\ln\left(\frac{n}{m}\right)+\mathcal{O}(m^{-1}).
	\end{gather*}
	Therefore, $\prod_{i=m}^{n}\left(\frac{2i}{2i+1}\right)=\frac{\sqrt{m}}{\sqrt{n}}(1+\mathcal{O}(m^{-1}))$ for $m<n$. Consequently, 
	\begin{align*}
	\mathbb{E}\{d(N_I)\}&=
	\begin{cases}
	\mathbb{E}\{d(N_{I_0})\}\frac{\sqrt{N_{I_0}}}{\sqrt{N_I}}(1+\mathcal{O}(N_{I_0}^{-1})) & \text{if  }  N_I> N_{I_0}\\
	\mathbb{E}\{d(N_{I_0})\}\frac{\sqrt{N_{I_0}}}{\sqrt{N_I}}+(1+\mathcal{O}(N_I^{-1})) & \text{if  } N_I< N_{I_0}
	\end{cases}\\&=\sqrt{N_{I_0}}\mathbb{E}\{d(N_{I_0})\}\frac{1}{\sqrt{N_I}}(1+\mathcal{O}(\max\{N_{I_0}^{-1},N_I^{-1}\}))
	\end{align*}
		
	\textbf{(ii)}We cannot directly apply the result of part (i) for a passenger with an arbitrary location inside $\mathcal{A}$. Nevertheless, based on the result of part (i) we can iteratively provide an approximation. For every point $x\in\mathcal{A}$ consider a disk $B_r(x)$ of radius $r>0$  around $x$,  $r>0$. We partition $\mathcal{A}$ into two sets $\mathcal{C}$ and $\mathcal{D}$, where (a)  $\mathcal{C} =\{x\in\mathcal{A}:B_r(x)\subset \mathcal{A}\}$ is the set of points such that $B_r(x)$ is contained in $\mathcal{A}$ , and (b) let $\mathcal{D}:=\{x\in\mathcal{A}:B_r(x)\not\subset \mathcal{A}\}$ is the set of points such that $B_r(x)$ is not completely contained in $\mathcal{A}$. 
	
	In the following, we first determine an approximate of $\mathbb{E}\{d(x)|N_I \text{\hspace*{-2pt} idle\hspace*{-2pt} vehicles}\}$ for $x\in\mathcal{C}$. We then provide an upper bound approximation for $\mathbb{E}\{d(x)|N_I \text{\hspace*{-2pt} idle\hspace*{-2pt} vehicles}\}$ for $x\in\mathcal{D}$. Putting the results of (ii-a) and (ii-b) together we provide an approximate equation for a general shape $\mathcal{A}$; see (Iteration 1). We then use the approximation provided in (Iteration 1), to provide a better approximation for $\mathbb{E}\{d(x)|N_I \text{\hspace*{-2pt} idle\hspace*{-2pt} vehicles}\}$ for $x\in\mathcal{D}$, and improve our approximation for a general shape $\mathcal{A}$; see (Iteration 2). We repeat the above process iteratively, and determine the best approximation by analyzing the limit point of the above iterative process which leads to the final approximation $\sqrt{N_{I_0}}\mathbb{E}\{d(N_{I_0})\}\frac{1}{\sqrt{N_I}}(1+\mathcal{O}(\max\{N_{I_0}^{-1},N_I^{-1}\}))$.

	Let $|\mathcal{A}|$ and $L$ denote the area and length of the  (assumed smooth) boundary  of $\mathcal{A}$, respectively. Then the areas of $\mathcal{C}$ and $D$ are  $|\mathcal{C}|:=A-\mathcal{O}(Lr)$ and $|\mathcal{D}|:=\mathcal{O}(Lr)$, respectively.
	
	\textbf{(ii-a)} First consider an arbitrary point $x\in \mathcal{C}$. Suppose $N_I$ idle vehicles  are uniformly and independently distributed in $\mathcal{A}$. Then the pdf of the number $m$ of these vehicles that lie in $B_r(x)$ is a binomial distribution $B(N_I,\frac{\pi r^2}{\mathcal{A}})$ with mean $N_I\frac{\pi r^2}{\mathcal{A}}$ and variance $N_I\frac{\pi r^2}{\mathcal{A}}\left(1-\frac{\pi r^2}{\mathcal{A}}\right)$. 
	
	Let $l:=\{\max |y-z|:x,y\in \mathcal{A}\}$ denote the largest distance between two points in $\mathcal{A}$. 
	Then, for every realization of $m>1$, the conditional expected shortest distance of an idle vehicle to $x$ is given by $\sqrt{N_{I_0}\pi r^2}\mathbb{E}\{d(N_{I_0} \pi r^2)\}\frac{1}{\sqrt{m}}\left(1+\mathcal{O}\left(\max\{(N_{I_0}\pi r^2)^{-1},m^{-1}\}\right)\right)\mathbf{1}\{m>0\}+\mathcal{O}(l)\mathbf{1}\{m=0\}$ from part (i) where $\mathbf{1}\{\cdot\}$ denotes the indicator function. Note that we modify the expression from part (i) since the shortest distance of an idle vehicle cannot exceed $l$ when $m=0$. Moreover, we substitute $N_{I_0}$ with $N_{I_0}\frac{\pi r^2}{\mathcal{A}}$ to reflect the fact among the $N_{I_0}$ idle vehicles uniformly distributed in $\mathcal{A}$, on average $N_{I_0}\frac{\pi r^2}{\mathcal{A}}$ vehicles are inside $B_r(x)$.
	
	Taking the expectation with respect to $m$, the expected shortest distance of an idle vehicle to $x$ is given by 	
	{\selectfont \fontsize{9}{10}
		\begin{align}
		\mathbb{E}\{\hspace*{-1pt}d(x)|N_I \text{\hspace*{-2pt} idle\hspace*{-2pt} vehicles}\}&\hspace*{-2pt}=\hspace*{-2pt}\mathbb{E}_m\hspace*{-2pt}\left\{\hspace*{-3pt}\sqrt{\hspace*{-1pt}N_{I_0} \frac{\pi r^2}{\mathcal{A}}}\mathbb{E}\{\hspace*{-1pt}d(N_{I_0} \frac{\pi r^2}{\mathcal{A}})\hspace*{-1pt}\}\frac{1}{\sqrt{m}}\hspace*{-2pt}\left(1\hspace*{-1pt}+\hspace*{-1pt}\mathcal{O}\hspace*{-2pt}\left(\max\{(N_{I_0}\frac{\pi r^2}{\mathcal{A}})^{-1}\hspace*{-1pt},\hspace*{-1pt}m^{-1}\}\hspace*{-1pt}\right)\hspace*{-1pt}\right)\hspace*{-2pt}\mathbf{1}\{m\hspace*{-1pt}>\hspace*{-1pt}0\}\hspace*{-1pt}+\hspace*{-1pt}\mathcal{O}(l)\mathbf{1}\{\hspace*{-1pt}m\hspace*{-1pt}=\hspace*{-1pt}0\hspace*{-1pt}\}\hspace*{-3pt}\right\}\hspace*{-2pt}\hspace*{-6pt}\label{proof-prop1-eq1}
		\end{align}
	}
	We can write the first term in this expectation as
	{\selectfont \fontsize{9}{10}
		\begin{gather}
		\mathbb{E}_m\left\{\sqrt{N_{I_0} \frac{\pi r^2}{\mathcal{A}}}\mathbb{E}\{d(N_{I_0} \frac{\pi r^2}{\mathcal{A}})\}\frac{1}{\sqrt{m}}\left(1+\mathcal{O}\left(\max\{(N_{I_0}\frac{\pi r^2}{\mathcal{A}})^{-1},m^{-1}\}\right)\right)\mathbf{1}\{m>0\}\right\}=\nonumber \\
		\sqrt{N_{I_0} \frac{\pi r^2}{\mathcal{A}}}\mathbb{E}\{\hspace*{-1pt}d(N_{I_0}\frac{\pi r^2}{\mathcal{A}})\hspace*{-1pt}\}\frac{1}{\sqrt{N_I\frac{\pi r^2}{\mathcal{A}}}}
		\mathbb{E}_m\hspace*{-2pt}\left\{\frac{1}{\sqrt{1\hspace*{-1pt}+\hspace*{-1pt}\frac{m-N_I\frac{\pi r^2}{\mathcal{A}}}{N_I\frac{\pi r^2}{\mathcal{A}}}}}
		\hspace*{-2pt}\left(\hspace*{-2pt}1\hspace*{-1pt}+\hspace*{-1pt}\mathcal{O}\hspace*{-2pt}\left(\hspace*{-2pt}\max\{(N_{I_0}\frac{\pi r^2}{\mathcal{A}})^{-1},(N_I\frac{\pi r^2}{\mathcal{A}})^{-1}(1\hspace*{-1pt}+\hspace*{-1pt}\frac{m-N_I\frac{\pi r^2}{\mathcal{A}}}{N_I\frac{\pi r^2}{\mathcal{A}}})^{-1}\}\right)\hspace*{-3pt}\right)\hspace*{-3pt}\mathbf{1}\{m\hspace*{-1pt}>\hspace*{-1pt}0\}\hspace*{-2pt}\right\}\nonumber\\
		=\sqrt{N_{I_0}}\mathbb{E}\{\hspace*{-1pt}d(N_{I_0} \frac{\pi r^2}{\mathcal{A}})\hspace*{-1pt}\}\frac{1}{\sqrt{N_I}}\left(1+\mathcal{O}\left(\max\{N_{I_0}^{-1},N_I^{-1}\}\right)\hspace*{-1pt}\right)\mathbb{P}\{m>0\}
		\end{gather} 
	}
	where the last equality is by Taylor expansion for terms $\frac{1}{\sqrt{1\hspace*{-1pt}+\hspace*{-1pt}\frac{m-N_I\frac{\pi r^2}{\mathcal{A}}}{N_I\frac{\pi r^2}{\mathcal{A}}}}}$ and $(1\hspace*{-1pt}+\hspace*{-1pt}\frac{m-N_I\frac{\pi r^2}{\mathcal{A}}}{N_I\frac{\pi r^2}{\mathcal{A}}})^{-1}$, along with the fact that $\mathbb{E}_m\{m-N_I\frac{\pi r^2}{\mathcal{A}}\}=0$ and $\mathbb{E}_m\left\{\left(\frac{m-N_I\frac{\pi r^2}{\mathcal{A}}}{N_I\frac{\pi r^2}{\mathcal{A}}}\right)^2\right\}=\mathcal{O}(N_I^{-1})$  .  
	
	Moreover, we have $\mathbb{P}\{m=0\}=(1-\frac{\pi r^2}{\mathcal{A}})^{N_I}$. Using inequality $\ln((1-x)^n)\leq -nx$, we have
	\begin{align*}
	\mathbb{P}\{m=0\}=\exp (\ln ((1-\frac{\pi r^2}{\mathcal{A}})^{N_I}))\leq \exp (-N_I\frac{\pi r^2}{\mathcal{A}}).
	\end{align*}
	
	Set $r=\mathcal{O}(N_I^{\frac{-1}{2}+\delta})$ for some $\delta\in(0,1/2]$. Then we have 
	\begin{align}\mathbb{P}\{m=0\}=\mathcal{O}(\exp (-N_I^{2\delta})).\label{proof-prop1-eq2}\end{align} 
	
	Consequently, $\mathbb{P}\{m>0\}=1-\mathbb{P}\{m=0\}=1+\mathcal{O}(\exp (-N_I^{2\delta}))$. Substituting $\mathbb{P}\{m>0\}$ and $\mathbb{P}\{m=0\}$ in (\ref{proof-prop1-eq1}), we get
	\begin{align}
	\mathbb{E}\{\hspace*{-1pt}d(x)|N_I \text{\hspace*{-2pt} idle\hspace*{-2pt} vehicles}\}=\mathbb{E}\{\hspace*{-1pt}d(x)|N_{I_0} \text{\hspace*{-2pt} idle\hspace*{-2pt} vehicles}\} \frac{\sqrt{N_{I_0}}}{\sqrt{N_I}}\left(1+\mathcal{O}\left(\hspace*{-2pt}\max\{N_{I_0}^{-1}N_I^{1-2\delta},N_I^{-2\delta}, \exp (-N_I^{2\delta})\}\right)\hspace*{-3pt}\right)\label{proof-prop1-e1}
	\end{align}
	for $x\in \mathcal{C}(x)$ and $\delta\in (0,1/2]$.
	
	\textbf{(ii-b)} For points $x\in\mathcal{D}$, consider the intersection of $B_r(x)$ and $\mathcal{A}$. We assume that $\frac{|B_r(x)\cap \mathcal{A}|}{r^2}=\mathcal{O}(1)$ since  $\mathcal{A}$ has a smooth boundary. (Here $ |\mathcal{R}|$ is the area of $\mathcal{R}$.) Apply an argument similar to the one given in (ii-a) for $B_r(x)\cap \mathcal{A}$, and we can show that with large enough probability the closest idle vehicles to $x$ has a distance smaller than $r$. More precisely,
	\begin{align}
	\mathbb{E}\{\hspace*{-1pt}d(x)|N_I \text{\hspace*{-2pt} idle\hspace*{-2pt} vehicles}\}\leq\mathcal{O}(r)+\mathcal{O}(exp(-N_Ir^{2}))=
	\mathcal{O}(N_I^{\frac{-1}{2}+\delta})+\mathcal{O}(exp(-N_I^{2\delta}))\label{proof-prop1-e2-0}
	\end{align}
	for $x\in\mathcal{D}$, where the last equality follows from $r=\mathcal{O}(N_I^{\frac{-1}{2}+\delta})$. Therefore, from (\ref{proof-prop1-e2-0}),
	{\selectfont \fontsize{9}{10}
	\begin{align*}
	\mathbb{E}\{\hspace*{-1pt}d(x)|N_I \text{\hspace*{-2pt} idle\hspace*{-2pt} vehicles}\}-\frac{\sqrt{N_{I_0}}}{\sqrt{N_I}}\mathbb{E}\{\hspace*{-1pt}d(x)|N_{I_0} \text{\hspace*{-2pt} idle\hspace*{-2pt} vehicles}\}
	&\leq \mathcal{O}(N_I^{\frac{-1}{2}+\delta})+\mathcal{O}(exp(-N_I^{2\delta}))-\frac{\sqrt{N_{I_0}}}{\sqrt{N_I}}\mathbb{E}\{\hspace*{-1pt}d(x)|N_{I_0} \text{\hspace*{-2pt} idle\hspace*{-2pt} vehicles}\}\\
	&\hspace*{-40pt}=\mathcal{O}(N_I^{\frac{-1}{2}+\delta})+\mathcal{O}(exp(-N_I^{2\delta}))+\frac{\sqrt{N_{I_0}}}{\sqrt{N_I}}\left[\mathcal{O}(N_{I_0}^{\frac{-1}{2}+\delta})+\mathcal{O}(exp(-N_{I_0}^{2\delta}))\right]\\
	&\hspace*{-40pt}=\mathcal{O}(N_I^{\frac{-1}{2}+\delta})+\mathcal{O}(exp(-N_I^{2\delta}))+\left(\mathcal{O}(N_I^{\frac{-1}{2}}N_{I_0}^\delta)+\mathcal{O}(N_{I_0}^{\frac{1}{2}}N_I^{-\frac{1}{2}}exp(-N_{I_0}^{2\delta})\right)
	\end{align*}
	}
	for $x\in\mathcal{D}$, where the first equality follows from 
	$$0\hspace*{-2pt}\leq \left|\hspace*{-2pt} \frac{\sqrt{N_{I_0}}}{\sqrt{N_I}}\mathbb{E}\{\hspace*{-1pt}d(x)|N_{I_0} \text{\hspace*{-2pt} idle\hspace*{-2pt} vehicles}\}\right|\hspace*{-2pt}\leq\hspace*{-2pt} \frac{\sqrt{N_{I_0}}}{\sqrt{N_I}} \hspace*{-2pt}	\left[\hspace*{-1pt}\mathcal{O}(N_{I_0}^{\frac{-1}{2}+\delta})\hspace*{-2pt}+\hspace*{-2pt}\mathcal{O}(exp(-N_{I_0}^{2\delta})\hspace*{-1pt})\hspace*{-1pt}\right].$$
	
	Therefore, for $x\in\mathcal{D}$, we can write
	{\selectfont \fontsize{9}{10}
	\begin{align}
	\mathbb{E}\{\hspace*{-1pt}d(x)|N_I \text{\hspace*{-2pt} idle\hspace*{-2pt} vehicles}\}\hspace*{-2pt}=\hspace*{-2pt}\frac{\sqrt{N_{I_0}}}{\sqrt{N_I}}\mathbb{E}\{\hspace*{-1pt}d(x)|N_{I_0} \text{\hspace*{-2pt} idle\hspace*{-2pt} vehicles}\}+\hspace*{-2pt}\mathcal{O}\hspace*{-2pt}\left(\max\{N_I^{\frac{-1}{2}+\delta}\hspace*{-1pt},exp(-N_I^{2\delta}),N_I^{\frac{-1}{2}}N_{I_0}^\delta\hspace*{-1pt},N_{I_0}^{\frac{1}{2}}N_I^{-\frac{1}{2}}exp(-N_{I_0}^{2\delta})\}\hspace*{-2pt}\right)\label{proof-prop1-e20}
	\end{align}
	}

	\textbf{(Iteration 1)} Using (\ref{proof-prop1-e1}) for $x\in\mathcal{C}$ and (\ref{proof-prop1-e20}) for $x\in\mathcal{D}$, along with $\mathbb{P}\{x\in\mathcal{D}\}=\frac{|\mathcal{D}|}{|\mathcal{A}|}=\mathcal{O}(\frac{Lr}{|\mathcal{A}|})=\mathcal{O}(N_I^{-\frac{1}{2}+\delta})$, we have,
	{\selectfont \fontsize{9}{10}
	\begin{align*}
	\mathbb{E}_{x\in\mathcal{A}}\{d(x)|N_I \text{\hspace*{-2pt} idle\hspace*{-2pt} vehicles}\}&=\Big[\mathbb{E}_x\{d(x)\mathbf{1}\{x\in\mathcal{C}\}|N_I \text{\hspace*{-2pt} idle\hspace*{-2pt} vehicles}\}\Big]+\Big[\mathbb{E}_x\{d(x)\mathbf{1}\{x\in\mathcal{D}\}|N_I \text{\hspace*{-2pt} idle\hspace*{-2pt} vehicles}\}\Big]\\
	&\hspace*{-110pt}= \hspace*{-2pt}\left[\frac{\sqrt{N_{I_0}}}{\sqrt{N_I}}\mathbb{E}_{x\in\mathcal{C}}\{d(x)|N_{I_0} \text{\hspace*{-2pt} idle\hspace*{-2pt} vehicles}\}\left(1\hspace*{-2pt}+\hspace*{-2pt}\mathcal{O}\hspace*{-2pt}\left(\hspace*{-2pt}\max\{N_{I_0}^{-1}N_I^{1-2\delta},N_I^{-2\delta},exp(-N_I^{2\delta})\}\right)\hspace*{-3pt}\right)\right](1-\mathbb{P}\{x\hspace*{-1pt}\in\hspace*{-1pt} \mathcal{D}\})\hspace*{-1pt}\\
	&\hspace*{-110pt}\hspace*{10pt}+\hspace*{-2pt}\left[\frac{\sqrt{N_{I_0}}}{\sqrt{N_I}}\mathbb{E}_{x\in\mathcal{D}}\{\hspace*{-1pt}d(x)|N_{I_0} \text{\hspace*{-2pt} idle\hspace*{-2pt} vehicles}\}+\hspace*{-2pt}\mathcal{O}\hspace*{-2pt}\left(\max\{N_I^{\frac{-1}{2}+\delta}\hspace*{-1pt},exp(-N_I^{2\delta}),N_I^{\frac{-1}{2}}N_{I_0}^\delta\hspace*{-1pt},N_{I_0}^{\frac{1}{2}}N_I^{-\frac{1}{2}}exp(-N_{I_0}^{2\delta})\}\right)\right]\hspace*{-2pt}\mathbb{P}\{x\hspace*{-1pt}\in\hspace*{-1pt}\mathcal{D}\}\\
	&\hspace*{-110pt}= \frac{\sqrt{N_{I_0}}}{\sqrt{N_I}}\mathbb{E}_{x\in\mathcal{A}}\{d(x)|N_{I_0} \text{\hspace*{-2pt} idle\hspace*{-2pt} vehicles}\}\hspace*{-2pt}\\
	&\hspace*{-110pt}\hspace*{10pt}+\left[\frac{\sqrt{N_{I_0}}}{\sqrt{N_I}}\mathbb{E}_{x\in\mathcal{C}}\{d(x)|N_{I_0} \text{\hspace*{-2pt} idle\hspace*{-2pt} vehicles}\}\mathcal{O}\hspace*{-2pt}\left(\hspace*{-2pt}\max\{N_{I_0}^{-1}N_I^{1-2\delta},N_I^{-2\delta},exp(-N_I^{2\delta})\}\right)(1-\mathbb{P}\{x\in\mathcal{D}\})\right]\\
	&\hspace*{-110pt}\hspace*{10pt}+\left[\mathcal{O}\hspace*{-2pt}\left(\max\{N_I^{\frac{-1}{2}+\delta}\hspace*{-1pt},exp(-N_I^{2\delta}),N_I^{\frac{-1}{2}}N_{I_0}^\delta\hspace*{-1pt},N_{I_0}^{\frac{1}{2}}N_I^{-\frac{1}{2}}exp(-N_{I_0}^{2\delta})\}\right)\mathbb{P}\{x\in\mathcal{D}\}\right]\\
	&\hspace*{-110pt}= \frac{\sqrt{N_{I_0}}}{\sqrt{N_I}}\mathbb{E}_{x\in\mathcal{A}}\{d(x)|N_{I_0} \text{\hspace*{-2pt} idle\hspace*{-2pt} vehicles}\}\hspace*{-2pt}\\
	&\hspace*{-110pt}\hspace*{10pt}+\frac{\sqrt{N_{I_0}}}{\sqrt{N_I}}\mathbb{E}_{x\in\mathcal{A}}\{d(x)|N_{I_0} \text{\hspace*{-2pt} idle\hspace*{-2pt} vehicles}\}\hspace*{-2pt}\Bigg[\frac{\mathbb{E}_{x\in\mathcal{C}}\{d(x)|N_{I_0} \text{\hspace*{-2pt} idle\hspace*{-2pt} vehicles}\}}{\mathbb{E}_{x\in\mathcal{A}}\{d(x)|N_{I_0} \text{\hspace*{-2pt} idle\hspace*{-2pt} vehicles}\}}\mathcal{O}\hspace*{-2pt}\left(\hspace*{-2pt}\max\{N_{I_0}^{-1}N_I^{1-2\delta},N_I^{-2\delta},exp(-N_I^{2\delta})\}\right)\hspace*{-2pt}(1\hspace*{-2pt}-\hspace*{-2pt}\mathbb{P}\{x\hspace*{-1pt}\in\hspace*{-1pt}\mathcal{D}\hspace*{-1pt}\})\\
	&\hspace*{-110pt}\hspace*{60pt}+\frac{\sqrt{N_I}}{\sqrt{N_{I_0}}}\frac{1}{\mathbb{E}_{x\in\mathcal{A}}\{d(x)|N_{I_0} \text{\hspace*{-2pt} idle\hspace*{-2pt} vehicles}\}}\mathcal{O}\hspace*{-2pt}\left(\max\{N_I^{\frac{-1}{2}+\delta}\hspace*{-1pt},exp(-N_I^{2\delta}),N_I^{\frac{-1}{2}}N_{I_0}^\delta\hspace*{-1pt},N_{I_0}^{\frac{1}{2}}N_I^{-\frac{1}{2}}exp(-N_{I_0}^{2\delta})\}\right)\mathbb{P}\{x\hspace*{-1pt}\in\hspace*{-1pt}\mathcal{D}\}\hspace*{-1pt}\Bigg]\\
	&\hspace*{-110pt}= \frac{\sqrt{N_{I_0}}}{\sqrt{N_I}}\mathbb{E}_{x\in\mathcal{A}}\{d(x)|N_{I_0} \text{\hspace*{-2pt} idle\hspace*{-2pt} vehicles}\}\hspace*{-2pt}\\
	&\hspace*{-110pt}\hspace*{10pt}+\frac{\sqrt{N_{I_0}}}{\sqrt{N_I}}\mathbb{E}_{x\in\mathcal{A}}\{d(x)|N_{I_0} \text{\hspace*{-2pt} idle\hspace*{-2pt} vehicles}\}\hspace*{-2pt}\Bigg[\mathcal{O}(1)\mathcal{O}\hspace*{-2pt}\left(\hspace*{-2pt}\max\{N_{I_0}^{-1}N_I^{1-2\delta},N_I^{-2\delta},exp(-N_I^{2\delta})\}\right)(1-\mathbb{P}\{x\in\mathcal{D}\})\\
	&\hspace*{-110pt}\hspace*{60pt}+\frac{\sqrt{N_I}}{\sqrt{N_{I_0}}}\frac{1}{\mathcal{O}(N_{I_0}^{\frac{-1}{2}})}\mathcal{O}\hspace*{-2pt}\left(\max\{N_I^{\frac{-1}{2}+\delta}\hspace*{-1pt},exp(-N_I^{2\delta}),N_I^{\frac{-1}{2}}N_{I_0}^\delta\hspace*{-1pt},N_{I_0}^{\frac{1}{2}}N_I^{-\frac{1}{2}}exp(-N_{I_0}^{2\delta})\}\right)\mathbb{P}\{x\in\mathcal{D}\}\Bigg]\\
	&\hspace*{-110pt}= \frac{\sqrt{N_{I_0}}}{\sqrt{N_I}}\mathbb{E}_{x\in\mathcal{A}}\{d(x)|N_{I_0} \text{\hspace*{-2pt} idle\hspace*{-2pt} vehicles}\}\hspace*{-2pt}\\
	&\hspace*{-110pt}\hspace*{10pt}+\frac{\sqrt{N_{I_0}}}{\sqrt{N_I}}\mathbb{E}_{x\in\mathcal{A}}\{d(x)|N_{I_0} \text{\hspace*{-2pt} idle\hspace*{-2pt} vehicles}\}\hspace*{-2pt}\Bigg[\mathcal{O}\hspace*{-2pt}\left(\hspace*{-2pt}\max\{N_{I_0}^{-1}N_I^{1-2\delta},N_I^{-2\delta},exp(-N_I^{2\delta})\}\right)(1-\mathcal{O}(N_I^{\frac{-1}{2}+\delta}))\\
	&\hspace*{-110pt}\hspace*{60pt}+\sqrt{N_I}\mathcal{O}\hspace*{-2pt}\left(\max\{N_I^{\frac{-1}{2}+\delta}\hspace*{-1pt},exp(-N_I^{2\delta}),N_I^{\frac{-1}{2}}N_{I_0}^\delta\hspace*{-1pt},N_{I_0}^{\frac{1}{2}}N_I^{-\frac{1}{2}}exp(-N_{I_0}^{2\delta})\}\right)\mathcal{O}(N_I^{\frac{-1}{2}+\delta})\Bigg]\\
	&\hspace*{-110pt}= \frac{\sqrt{N_{I_0}}}{\sqrt{N_I}}\mathbb{E}_{x\in\mathcal{A}}\{d(x)|N_{I_0} \text{\hspace*{-2pt} idle\hspace*{-2pt} vehicles}\}\Bigg(1+\mathcal{O}
	(\max{N_{I_0}^{-1}N_I^{1-2\delta},N_I^{-2\delta},exp(-N_I^{2\delta}),N_I^{\frac{-1}{2}+2\delta},N_I^\delta exp(-N_I^{2\delta}),N_I^{-\frac{1}{2}+\delta}N_{I_0}^\delta,}\\
	&\hspace*{270pt}N_{I_0}^{\frac{1}{2}}N_I^{\frac{-1}{2}+\delta}exp(-N_{I_0}^{2\delta}))\Bigg),	
	\end{align*}
	}
	for $\delta\in(0,\frac{1}{2}]$. 
	

	Setting $\delta=\frac{1}{8}$ we have,
	\begin{align}
	\hspace*{-8pt}\mathbb{E}_{x\in\mathcal{A}}\{d(x)|N_I \text{\hspace*{-2pt} idle\hspace*{-2pt} vehicles}\}\hspace*{-2pt}=\hspace*{-2pt} \frac{\sqrt{N_{I_0}}}{\sqrt{N_I}}\mathbb{E}_x\{d(x)|N_{I_0} \text{\hspace*{-2pt} idle\hspace*{-2pt} vehicles}\}\hspace*{-2pt}\left(\hspace*{-2pt}1\hspace*{-2pt}+\hspace*{-2pt}\mathcal{O}\hspace*{-2pt}\left(\hspace*{-2pt}\max\{N_{I_0}^{-1}N_I^{\frac{3}{4}},N_I^{-\frac{1}{4}},N_I^{-\frac{3}{8}}N_{I_0}^\frac{1}{8}\}\hspace*{-2pt}\right)\hspace*{-5pt}\right)\hspace*{-2pt}.\hspace*{-3pt}\label{proof-prop1-R1}
	\end{align}
	where we neglect $\mathcal{O}(exp(-N_I^{\frac{1}{4	}}))$, $\mathcal{O}(N_I^{0.125}exp(-N_I^{\frac{1}{4	}}))$ and $\mathcal{O}(N_{I_0}^{0.5}N_I^{-0.375}exp(-N_I^{\frac{1}{4}}))$ in comparison to $\mathcal{O}\left(\hspace*{-2pt}N_I^{-\frac{1}{4}}\right)$.
	
	\textbf{(Iteration 2)} We can now use (\ref{proof-prop1-R1}) to provide a better approximation for $x\in\mathcal{D}$ in part (ii-b) and iterate the process. Note that the number of idle vehicles in $\mathcal{D}$ (on average) is equal to $\frac{|\mathcal{D}|}{|\mathcal{A}|}N_I=\mathcal{O}(r)N_I=N_I^{\frac{1}{2}+\delta}$; similarly, we need to substitute $N_{I_0}$ with $N_{I_0}\mathcal{O}(r)=N_{I_0}N_I^{\frac{-1}{2}+\delta}$ in (\ref{proof-prop1-R1}).\footnote{We note that the number of vehicles in $\mathcal{D}$ is a random variable. Nevertheless, considering the variation in the number of vehicles around the average of $N_I^{\frac{1}{2}+\delta}$ one can follow a similar approach we used in (ii-a) using Taylor expansion, and show that the error term due to such a randomness is smaller than other error terms and does not affect the result.} Therefore, 
	{\selectfont \fontsize{9}{10}
	\begin{align}
	\mathbb{E}\{\hspace*{-1pt}d(x)|N_I \text{\hspace*{-2pt} idle\hspace*{-2pt} vehicles}\}=\frac{\sqrt{N_{I_0}}}{\sqrt{N_I}}\mathbb{E}_x\{d(x)|N_{I_0} \text{\hspace*{-2pt} idle\hspace*{-2pt} vehicles}\}\left(1+\mathcal{O}\left(\hspace*{-2pt}\max\{N_{I_0}^{-1}N_I^{\frac{7}{8}-\frac{1}{4}\delta},N_I^{-\frac{1}{8}-\frac{1}{4}\delta},N_{I_0}^{\frac{1}{8}}N_I^{\frac{-1}{4}-\frac{1}{4}\delta}\}\right)\hspace*{-3pt}\right)\label{proof-prop1-e21}
	\end{align}}
	for $x\in\mathcal{D}$.\footnote{We note that the closest idle vehicles for $x\in\mathcal{D}$ does not necessarily belong to $\mathcal{D}$ and can be inside $\mathcal{C}$. Nevertheless, we can still use (\ref{proof-prop1-R1}) for the expected distance of the closest idle vehicle for $x$. This is because for area $\mathcal{D}$ we can follow an argument similar to the one that leads to (\ref{proof-prop1-R1}), and divide it to an interior region $\tilde{\mathcal{C}}$ and an exterior region $\tilde{\mathcal{D}}$, ; however, we consider the border of $\mathcal{D}$ that separates it from $\mathcal{C}$ as a part of the interior region $\tilde{\mathcal{C}}$ (and not $\tilde{\mathcal{D}}$). Consequently, equation (\ref{proof-prop1-R1}) is also applicable for $x\in\mathcal{D}$. }
	
	Using (\ref{proof-prop1-e1}) for $x\in\mathcal{C}$ and (\ref{proof-prop1-e21}) for $x\in\mathcal{D}$, along with $|\mathcal{D}|=\mathcal{O}(r)=\mathcal{O}(N^{\frac{-1}{2}+\delta})$, we have
	{\selectfont \fontsize{9}{10}
	\begin{align*}
	\mathbb{E}_{x\in\mathcal{A}}\{d(x)|N_I \text{\hspace*{-2pt} idle\hspace*{-2pt} vehicles}\}&=\mathbb{E}_x\{d(x)\mathbf{1}\{x\in\mathcal{C}\}|N_I \text{\hspace*{-2pt} idle\hspace*{-2pt} vehicles}\}+\mathbb{E}_x\{d(x)\mathbf{1}\{x\in\mathcal{D}\}|N_I \text{\hspace*{-2pt} idle\hspace*{-2pt} vehicles}\}\\
	&\hspace*{-120pt}= \frac{\sqrt{N_{I_0}}}{\sqrt{N_I}}\mathbb{E}_x\{d(x)|N_{I_0} \text{\hspace*{-2pt} idle\hspace*{-2pt} vehicles}\}\left(1\hspace{-2pt}+\hspace{-2pt}\mathcal{O}\hspace{-2pt}\left(\hspace*{-2pt}\max\{N_{I_0}^{-1}N_I^{1-2\delta},N_I^{-2\delta},N_{I_0}^{-1}N_I^{\frac{3}{8}+\frac{3}{4}\delta},N_I^{-\frac{5}{8}+\frac{3}{4}\delta},N_{I_0}^{\frac{1}{8}}N_I^{-\frac{3}{4}+\frac{3}{4}\delta},exp(-N_I^{2\delta})\}\right)\hspace*{-3pt}\right)\hspace*{-3pt}.
	\end{align*}
	}
	Setting $\delta =\frac{5}{22}$, we have an improved approximation
	{\selectfont \fontsize{9}{10}
	\begin{align}
	\mathbb{E}_{x\in\mathcal{A}}\{d(x)|N_I \text{\hspace*{-2pt} idle\hspace*{-2pt} vehicles}\}= \frac{\sqrt{N_{I_0}}}{\sqrt{N_I}}\mathbb{E}_x\{d(x)|N_{I_0} \text{\hspace*{-2pt} idle\hspace*{-2pt} vehicles}\}\left(1+\mathcal{O}\left(\hspace*{-2pt}\max\{N_{I_0}^{-1}N_I^{\frac{6}{11}},N_I^{-\frac{5}{11}},N_{I_0}^{\frac{1}{8}}N_I^{-\frac{5}{11}-\frac{1}{8}}\}\right)\hspace*{-3pt}\right)\hspace*{-3pt}\label{proof-prop1-R2},
	\end{align}
	}
	where we neglect $\mathcal{O}(exp(-N_I^{\frac{5}{11}}))$ in comparison to $\mathcal{O}\left(\hspace*{-2pt}N_I^{-\frac{5}{11}}\right)$.

	\textbf{(Iteration $K$)} We can iterate the same process similar to the one described in iteration 2. Assume that at the end of iteration $K-1$ we show that
	{\selectfont \fontsize{9}{10}
	\begin{align}
	\hspace*{-8pt}\mathbb{E}_{x\in\mathcal{A}}\{d(x)|N_I \text{\hspace*{-2pt} idle\hspace*{-2pt} vehicles}\}\hspace*{-2pt}=\hspace*{-2pt} \frac{\sqrt{N_{I_0}}}{\sqrt{N_I}}\mathbb{E}_x\{d(x)|N_{I_0} \text{\hspace*{-2pt} idle\hspace*{-2pt} vehicles}\}\hspace*{-3pt}\left(\hspace*{-2pt}1\hspace*{-2pt}+\hspace*{-2pt}\mathcal{O}\hspace*{-2pt}\left(\hspace*{-2pt}\max\{N_{I_0}^{-1}N_I^{\alpha({K-1})}\hspace*{-2pt},\hspace*{-1pt}N_I^{-(1-\alpha({K-1}))},\hspace*{-1pt}N_{I_0}^{\frac{1}{8}}N_I^{-(1-\alpha({K-1}))-\frac{1}{8}}\}\hspace*{-2pt}\right)\hspace*{-3pt}\right)\hspace*{-3pt}.\hspace*{-3pt}\label{proof-prop1-RK-1}
	\end{align}}
	where $\alpha(K-1)\in(0,\frac{1}{2}]$. Note that in  iteration 2 we have $\alpha({2})=\frac{6}{11}$. We now use (\ref{proof-prop1-RK-1}) to provide a better approximation for $x\in\mathcal{D}$ and iterate the process another time. Once again note that the number of idle vehicles in $\mathcal{D}$ (on average) is equal to $\frac{|\mathcal{D}|}{|\mathcal{A}|}N_I=\mathcal{O}(r)N_I=N_I^{\frac{1}{2}+\delta}$; similarly, we substitute $N_{I_0}$ with $N_{I_0}\mathcal{O}(r)=N_{I_0}N_I^{\frac{-1}{2}+\delta}$ in (\ref{proof-prop1-R2}). Consequently, 
	{\selectfont \fontsize{9}{10}\begin{align}
	&\hspace*{20pt}\mathbb{E}\{\hspace*{-1pt}d(x)|N_I \text{\hspace*{-2pt} idle\hspace*{-2pt} vehicles}\}=\nonumber\\
	&\hspace*{-7pt}\frac{\sqrt{N_{I_0}}}{\sqrt{N_I}}\mathbb{E}_x\{d(x)|N_{I_0} \text{\hspace*{-2pt} idle\hspace*{-2pt} vehicles}\}\Bigg(1+\mathcal{O}\Bigg(\hspace*{-2pt}\max\{N_{I_0}^{-1}N_I^{\frac{1+\alpha({K-1})}{2}-(1-\alpha({K-1}))\delta},N_I^{-\frac{1-\alpha({K-1})}{2}-(1-\alpha({K-1}))\delta}\nonumber\\
	&\hspace*{265pt},N_{I_0}^{\frac{1}{8}}N_I^{-\frac{1-\alpha({K-1})}{2}-(1-\alpha({K-1}))\delta-\frac{1}{8}}\}\Bigg)\hspace*{-3pt}\Bigg)\hspace*{-1pt},\hspace*{-3pt}\label{proof-prop1-e2K}
	\end{align}}
	for $x\in\mathcal{D}$.
	
	Using (\ref{proof-prop1-e1}) for $x\in\mathcal{C}$ and (\ref{proof-prop1-e2K}) for $x\in\mathcal{D}$, along with $|\mathcal{D}|=\mathcal{O}(r)=\mathcal{O}(N^{\frac{-1}{2}+\delta})$, we have
	{\fontsize{9}{10}\selectfont
		\begin{align*}
		\mathbb{E}_{x\in\mathcal{A}}\{d(x)|N_I \text{\hspace*{-2pt} idle\hspace*{-2pt} vehicles}\}&=\mathbb{E}_x\{d(x)\mathbf{1}\{x\in\mathcal{C}\}|N_I \text{\hspace*{-2pt} idle\hspace*{-2pt} vehicles}\}+\mathbb{E}_x\{d(x)\mathbf{1}\{x\in\mathcal{D}\}|N_I \text{\hspace*{-2pt} idle\hspace*{-2pt} vehicles}\}\\
		&\hspace*{-110pt}= \hspace*{-2pt}\frac{\sqrt{N_{I_0}}}{\sqrt{N_I}}\mathbb{E}_x\{d(x)|N_{I_0} \text{\hspace*{-2pt} idle\hspace*{-2pt} vehicles}\}\hspace*{-3pt}\Bigg(\hspace*{-2pt}1\hspace*{-2pt}+\hspace*{-2pt}\mathcal{O}\hspace*{-2pt}\Bigg(\hspace*{-2pt}\max\{N_{I_0}^{-1}N_I^{1-2\delta},N_I^{-2\delta},N_{I_0}^{-1}N_I^{\frac{\alpha({K-1})}{2}+\alpha({K-1})\delta},N_I^{-\frac{2-\alpha({K-1})}{2}+\alpha({K-1})\delta},\\
		&\hspace*{215pt}N_{I_0}^{\frac{1}{8}}N_I^{-\frac{2-\alpha({K-1})}{2}+\alpha({K-1})\delta-\frac{1}{8}},exp(-N_I^{2\delta})\}\hspace*{-2pt}\Bigg)\hspace*{-3pt}\Bigg)\hspace*{-1pt}.
		\end{align*}
	}
	Setting $\delta =\frac{2-\alpha({K-1})}{2(2+\alpha({K-1}))}$, we have an improved approximation for iteration $K$ as
	{\fontsize{9}{10}\selectfont
	\begin{align}
	\hspace*{-3pt}\mathbb{E}_{x\in\mathcal{A}}\{d(x)|N_I \text{\hspace*{-2pt} idle\hspace*{-2pt} vehicles}\}&\hspace*{-2pt}=\hspace*{-2pt} \frac{\sqrt{N_{I_0}}}{\sqrt{N_I}}\mathbb{E}_x\{d(x)|N_{I_0} \text{\hspace*{-2pt} idle\hspace*{-2pt} vehicles}\}\hspace*{-2pt}\left(\hspace*{-2pt}1\hspace*{-2pt}+\hspace*{-2pt}\mathcal{O}\hspace*{-2pt}\left(\hspace*{-2pt}\max\{N_{I_0}^{-1}N_I^{\frac{2\alpha({K-1})}{2+\alpha({K-1})}},N_I^{-\frac{2-\alpha({K-1})}{2+\alpha({K-1})}},N_{I_0}^{\frac{1}{8}}N_I^{-\frac{2-\alpha({K-1})}{2+\alpha({K-1})}-\frac{1}{8}}\}\hspace*{-2pt}\right)\hspace*{-3pt}\right)\nonumber\\
	&=\frac{\sqrt{N_{I_0}}}{\sqrt{N_I}}\mathbb{E}_x\{d(x)|N_{I_0} \text{\hspace*{-2pt} idle\hspace*{-2pt} vehicles}\}\left(1+\mathcal{O}\left(\hspace*{-2pt}\max\{N_{I_0}^{-1}N_I^{\alpha(K)},N_I^{1-\alpha(K)},N_{I_0}^{\frac{1}{8}}N_I^{1-\alpha(K)-\frac{1}{8}}\}\right)\hspace*{-3pt}\right)\hspace*{-3pt}\label{proof-prop1-R3},
	\end{align}}
	where $\alpha(K)=\frac{2\alpha({K-1})}{2+\alpha({K-1})}$ and we neglect $\mathcal{O}(\exp (-N_I^{\frac{2-\alpha({K-1})}{2+\alpha({K-1})}}))$ in comparison to $\mathcal{O}\left(\hspace*{-2pt}N_I^{-\frac{2-\alpha({K-1})}{2+\alpha({K-1})}}\right)$.

	Iterating the above process, it is easy to show that the sequence of $\alpha(K)=\frac{2\alpha({K-1})}{2+\alpha({K-1})}$ converges to $\alpha^*=0$. 
	
	Therefore, we have 
	\begin{align}
	\mathbb{E}_{x\in\mathcal{A}}\{d(x)|N_I \text{\hspace*{-2pt} idle\hspace*{-2pt} vehicles}\}&= \frac{\sqrt{N_{I_0}}}{\sqrt{N_I}}\mathbb{E}_x\{d(x)|N_{I_0} \text{\hspace*{-2pt} idle\hspace*{-2pt} vehicles}\}\left(1+\mathcal{O}\left(\hspace*{-2pt}\max\{N_{I_0}^{-1},N_I^{-1},N_{I_0}^{\frac{1}{8}}N_I^{-\frac{9}{8}}\}\right)\hspace*{-3pt}\right).\nonumber
	\end{align}

\end{proof}

\subsection*{\bf{B: Proof of Proposition \ref{feasibility}}}

\begin{proof}
We can represent $(N,\lambda, p_f, p_d)$ as a function of total travel cost $c$ and driver wage $w$:
\begin{subnumcases}{\label{variablefuncwc}}
\lambda=\lambda_0\left[ 1-F_p(c)\right]\label{demand_cw}\\
N=N_0 F_d(w)\label{supply_cw} \\
p_f=\dfrac{1}{\beta}\left[ c-\alpha t_w\bigg( N_0F_d(w)-\lambda_0[1-F_p(c)]/\mu  \bigg)\right] \label{pf_cw}\\
p_d= \dfrac{wN}{\lambda}=w \dfrac{N_0 F_d(w)}{\lambda_0\left[ 1-F_p(c)\right]} \label{pd_cw}
\end{subnumcases}
where (\ref{pf_cw}) and (\ref{pd_cw}) are obtained based on the definition of $c$ and $w$, i.e., (\ref{costoftravel}) and (\ref{driverwage}). 
To prove Proposition \ref{feasibility}, it suffices to show that there exists $c>0$ and $w>0$ such that
\begin{subnumcases}{\label{variablefuncwc2}}
 c> \alpha t_w\bigg( N_0F_d(w)-\lambda_0[1-F_p(c)]/\mu  \bigg) \label{1_cw}\\
F_p(c)<1\label{2_cw}\\
F_d(w)>0\label{3_cw}
\end{subnumcases}
Note that (\ref{1_cw}) corresponds to $p_f>0$, and (\ref{2_cw}) and (\ref{3_cw}) guarantee that $p_d>0$, $N>0$ and $\lambda>0$. As $t_w$ is decreasing (see Assumption \ref{assump_pickuptimefunc}), it suffices to prove that there exists $c>0$ and $w>0$ such that $F_p(c)<1$,  $F_d(w)>0$ and that:
\begin{equation}
\label{proveCW1}
c> \alpha t_w\bigg( N_0-\lambda_0[1-F_p(c)]/\mu  \bigg).
\end{equation}
The left hand side of (\ref{proveCW1}) is an increasing function of $c$, while the right hand side is decreasing function. Define $c^*=\inf_c\{c|F_p(c)=1\}$. It suffices to show that $c^*>\alpha t_w(N_0)$. This is equivalent to $F_p(\alpha t_w(N_0))<1$, which completes the proof. 
\end{proof}

\subsection*{\bf{C: Proof of Proposition \ref{1storderconditionsufficient}}}

\begin{proof}
To prove Proposition \ref{1storderconditionsufficient}, we first show that there is at most one solution to  (\ref{1srordercondition_proposition1}). Then we show that this solution exists, and it coincides with the globally optimal solution to (\ref{optimalpricing}).

{\bf Uniqueness}: let $f_p(c)$ and $f_d(w)$ be the probability density function of $c$ and $w$, respectively. Since $c$ and $w$ are subject to uniform distribution, we have: $f_p(c)=e_p* \mathbbm{1}_{e_pc\leq 1}$ and $f_d(w)=e_d*\mathbbm{1}_{e_dw\leq 1}$, where $\mathbbm{1}_A$ is the indicator function of $A$. Assume for the moment that (\ref{1srordercondition_proposition1}) admits a non-trivial solution. Denote it as $(\tilde{p}_f,\tilde{p}_d,\tilde{\lambda},\tilde{N})$, and denote $\tilde{c}$ and $\tilde{w}$ as the corresponding passenger cost and driver wage, respectively. We first show that  $e_p\tilde{c}\leq 1$ and $e_d\tilde{w}\leq 1$.  This is because if $e_p\tilde{c}>1$, then $\lambda=N=0$, this is a trivial solution. If $e_d\tilde{w}>1$, we can decrease $\tilde{p}_d$ without affecting $\tilde{\lambda}$, indicating that $\dfrac{\partial \Pi}{\partial p_d}<0$, which contradicts with (\ref{proposition1_constraint2}). Therefore, we can rewrite (\ref{proposition1_constraint3}) and (\ref{proposition1_constraint4}) as:
\begin{align*}
\begin{cases}
\lambda=\lambda_0 \left[ 1-e_p \left(\dfrac{\alpha M}{\sqrt{N-\lambda / \mu}}+\beta p_f\right) \right] \\
N=N_0 e_d \left( \dfrac{\lambda p_d}{N} \right).
\end{cases}
\end{align*}
This can be further simplified to:
\begin{equation}
\label{simplifiedDS}
\lambda=\lambda_0 \left[ 1-e_p \left(\dfrac{\alpha M}{\sqrt{\sqrt{N_0e_d\lambda p_d}-\lambda / \mu}}+\beta p_f\right) \right]
\end{equation}
It suffices to show that there exists at most one set of $(\tilde{p}_f,\tilde{p}_d,\tilde{\lambda})$ that satisfies (\ref{proposition1_constraint1}), (\ref{proposition1_constraint2}) and (\ref{simplifiedDS}). Using the implicit function theorem on (\ref{simplifiedDS}), we can derive $\dfrac{\partial \lambda}{p_f}$ and $\dfrac{\partial \lambda}{p_d}$, thus (\ref{proposition1_constraint1})- (\ref{proposition1_constraint2}) becomes:
\begin{subnumcases}{\label{complexpartial}}
\dfrac{-\lambda_0\beta e_p\left(\sqrt{N_0e_dp_d\lambda}-\lambda/\mu\right)^{3/2}}{\left(\sqrt{N_0e_dp_d\lambda}-\lambda/\mu\right)^{3/2}-\dfrac{1}{2}\alpha M e_p\lambda_0 \left[ \dfrac{1}{2}\sqrt{N_0e_dp_d/\lambda} -1/\mu  \right]}(p_f-p_d)+\lambda=0  \label{complexpartial1}\\
\dfrac{\dfrac{1}{4}\lambda_0e_p\alpha M \sqrt{N_0e_d\lambda/p_d}}{\left(\sqrt{N_0e_dp_d\lambda}-\lambda/\mu\right)^{3/2}-\dfrac{1}{2}\alpha M e_p\lambda_0 \left[ \dfrac{1}{2}\sqrt{N_0e_dp_d/\lambda} -1/\mu  \right]}(p_f-p_d)-\lambda=0\label{complexpartial2}
\end{subnumcases} 
This reduces to:
\begin{subnumcases}{\label{simplypartial}}
\dfrac{1}{4}\lambda_0e_p \alpha M p_f \sqrt{N_0 e_d\lambda/p_d}=\lambda \left(\sqrt{N_0e_dp_d\lambda}-\lambda/\mu\right)^{3/2}+ \alpha M e_p \lambda_0 \lambda/ 2\mu \label{simplypartial1}\\
\beta\sqrt{p_d/\lambda} \left(\sqrt{N_0e_dp_d\lambda}-\lambda/\mu\right)^{3/2}=\dfrac{1}{4} \alpha M \sqrt{N_0e_d}
\label{simplypartial2} \\
p_f\lambda_0 e_p \beta \sqrt{N_0 e_d}= \sqrt{N_0e_d} \lambda+ 2\sqrt{\lambda p_d} e_p\lambda_0 \beta/\mu
\label{simplypartial3},
\end{subnumcases} 
where (\ref{simplypartial1}) directly follows from (\ref{complexpartial2}), (\ref{simplypartial2}) follows from (\ref{complexpartial1}) and (\ref{complexpartial2}), and (\ref{simplypartial3}) are derived by plugging (\ref{simplypartial2}) into (\ref{simplypartial1}). 

Assume there is another solution, denoted  $(p_f',p_f',\lambda',N')$. Without loss of generality, suppose $p_f'\leq \tilde{p}_f$. If $p_f'< \tilde{p}_f$, then there are three cases:
\begin{itemize}
\item $p_d'\geq \tilde{p}_d$. Based on (\ref{simplifiedDS}), we have $\lambda'>\tilde{\lambda}$. However, (\ref{simplypartial3}) dictates that $\lambda'<\tilde{\lambda}$. A contradiction.
\item $p_d'< \tilde{p}_d$ and $\sqrt{N_0e_d p_d'\lambda'}-\lambda'/\mu \geq \sqrt{N_0e_d\tilde{p}_d\tilde{\lambda}}-\tilde{\lambda}/\mu$. Note that (\ref{simplypartial1}) is equivalent to:
\begin{equation}
\label{temp1}
\dfrac{1}{4}\lambda_0e_p \alpha M \sqrt{N_0 e_d\lambda} p_f=\sqrt{p_d\lambda} \left[ \left(\sqrt{N_0e_dp_d\lambda}-\lambda/\mu\right)^{3/2}+ \alpha M e_p \lambda_0 / 2\mu \right]
\end{equation}
Therefore, based on (\ref{temp1}), $p_d'\lambda'<\tilde{p}_d\tilde{\lambda}$.  On the other hand, since $\sqrt{N_0e_d p_d'\lambda'}-\lambda'/\mu \geq \sqrt{N_0e_d\tilde{p}_d\tilde{\lambda}}-\tilde{\lambda}/\mu$ and $p_f'<\tilde{p}_f$, based on (\ref{simplifiedDS}), we have $\lambda'>\tilde{\lambda}$. This indicates that $\sqrt{N_0e_d p_d'\lambda'}-\lambda'/\mu < \sqrt{N_0e_d\tilde{p}_d\tilde{\lambda}}-\tilde{\lambda}/\mu$. A contradiction.

\item$p_d'< \tilde{p}_d$ and $\sqrt{N_0e_d p_d'\lambda'}-\lambda'/\mu < \sqrt{N_0e_d\tilde{p}_d\tilde{\lambda}}-\tilde{\lambda}/\mu$. Note that (\ref{simplypartial2}) is equivalent to:
\begin{equation}
\label{temp2}
\beta\sqrt{p_d} (\sqrt{N_0e_dp_d}-\sqrt{\lambda}/\mu) \left(\sqrt{N_0e_dp_d\lambda}-\lambda/\mu\right)^{1/2}=\dfrac{1}{4} \alpha M \sqrt{N_0e_d}
\end{equation}
As $p_d'< \tilde{p}_d$ and $\sqrt{N_0e_d p_d'\lambda'}-\lambda'/\mu < \sqrt{N_0e_d\tilde{p}_d\tilde{\lambda}}-\tilde{\lambda}/\mu$, (\ref{temp2}) indicates that $\lambda'<\tilde{\lambda}$. On the other hand, (\ref{simplypartial2}) also indicates that  $p_d'/\lambda'>\tilde{p}_d/\tilde{\lambda}$. Note that  (\ref{simplifiedDS}) can be written as
\begin{equation}
\label{temp3}
\lambda=\lambda_0 \left[ 1-e_p \left(\dfrac{\alpha M}{\sqrt{ \sqrt{\lambda N_0e_d p_d}-\lambda / \mu}}+\beta p_f\right) \right]
\end{equation}
As $p_f'<\tilde{p}_f$ and $p_d'/\lambda'>\tilde{p}_d/\tilde{\lambda}$, (\ref{temp3}) indicates that $\lambda'>\tilde{\lambda}$. A contradiction. 
\end{itemize}

If $p_f'=\tilde{p}_f$ and $p_d'< \tilde{p}_d$, then we can find a contradiction by exactly the same argument. In addition, if $p_f'=\tilde{p}_f$ and $p_d'= \tilde{p}_d$, then  $\lambda$ and $N$ are uniquely determined by (\ref{simplypartial}). 



{\bf Existence and Optimality}: Based on (\ref{variablefuncwc}), we can represent $(N,\lambda, p_f, p_d)$ as a function of total travel cost $c$ and driver wage $w$. Therefore, the platform profit $\lambda(p_f-p_d)$ is a  function of $c$ and $w$: 
\begin{equation}
\label{profitCW}
R=\lambda_0[1-F_p(c)]*\left[ \dfrac{c}{\beta}-\dfrac{\alpha}{\beta}   t_w\bigg( N_0F_d(w)-\lambda_0[1-F_p(c)]/\mu  \bigg)- w \dfrac{N_0 F_d(w)}{\lambda_0\left[ 1-F_p(c)\right]}\right]
\end{equation}
Note that $R$ is well-defined if $\lambda>0$, $N>0$, and the number of idle vehicle is positive, i.e., $N_0F_d(w)
>\lambda_0 [1-F_d(c)]/\mu$. Therefore, $R$ is a continuous function of $c$ and $w$ defined in 
\[\mathcal{D}=\{(c,w)|0\leq c<1/e_r, 0<w\leq 1/e_r, N_0e_d w>\lambda_0(1-e_pc)/\mu\}.\]
It can be verified that $R\leq 0$ on the boundary of $\mathcal{D}$. Based on our assumption, there exists $(c^*, w^*)\in \mathcal{D}$ at which $R>0$. This indicates that the optimal solution to (\ref{optimalpricing}) is in the interior of $\mathcal{D}$, which solves (\ref{1srordercondition_proposition1}). Since (\ref{1srordercondition_proposition1}) has at most one solution, this completes the proof.

\end{proof}

\subsection*{\bf{D: Proof of Theorem \ref{freelunchtheorem}}}
\begin{proof}
We will first prescribe a procedure to compute the optimal solution to (\ref{optimalpricing_wage}), then we use the procedure to prove that $\nabla_+N^*(\tilde{w})>0$ and $\nabla_+\lambda^*(\tilde{w})>0$.

{\bf Computation procedure}: Let $(p_f^*,p_d^*,\lambda^*,N^*)$ denote the optimal solution to (\ref{optimalpricing_wage}). Note that  (\ref{supply_constraint_wage}) at equality uniquely defines an increasing mapping from $N$ to $\lambda p_d$. Denote this mapping as $\lambda p_d=g(N)$, then (\ref{supply_constraint_wage})  reduces to $\lambda p_d \geq g(N)$, (\ref{minimum_wage}) reduces to $\lambda p_d\geq Nw$. Therefore, (\ref{supply_constraint_wage}) and (\ref{minimum_wage}) can be combined as $\lambda p_d\geq \max \{g(N),wN\}$,  and  (\ref{optimalpricing_wage}) is equivalent to:
\begin{align}
\label{simplieidoptimalpricing_wage}
&\max_{p_f,N}\quad \lambda p_f-\max \{g(N),wN\} \\
& \text{s.t. } \lambda=\lambda_0 \left[ 1-F_p \left(\dfrac{\alpha M}{\sqrt{N-\lambda / \mu}}+\beta p_f\right) \right]. \nonumber
\end{align}
Consider the following problem:
\begin{align}
\label{simplieidoptimalpricing_wage1}
&\max_{p_f,N}\quad \lambda p_f-g(N) \\
& \text{s.t. } \lambda=\lambda_0 \left[ 1-F_p \left(\dfrac{\alpha M}{\sqrt{N-\lambda / \mu}}+\beta p_f\right) \right]. \label{temp7}
\end{align}
Let $(\bar{p_f},\bar{N},\bar{\lambda})$ be the optimal solution to (\ref{simplieidoptimalpricing_wage1}), and denote the optimal value as $\bar{R}$. To facilitate our discussion, define $\Pi_1(N)$ as a function representing the optimal value of (\ref{simplieidoptimalpricing_wage1}) for any {\em given} $N$, then $\bar{R}=\max_N \Pi_1(N)$.  Similarly, define:
\begin{align}
\label{simplieidoptimalpricing_wage2}
&\max_{p_f,N}\quad \lambda p_f-wN \\
& \text{s.t. } \lambda=\lambda_0 \left[ 1-F_p \left(\dfrac{\alpha M}{\sqrt{N-\lambda / \mu}}+\beta p_f\right) \right]. \label{constraint4case2}
\end{align}
Let $(\hat{p_f},\hat{N},\hat{\lambda})$ be the optimal solution to (\ref{simplieidoptimalpricing_wage2}), and denote the optimal value as $\hat{R}$.  Define $\Pi_2(N)$ as the optimal value of (\ref{simplieidoptimalpricing_wage2}) for any given $N$, then $\hat{R}=\max_N \Pi_2(N)$.  
\begin{lemma}
\label{wage_lemma}
There are three cases for the solution to (\ref{optimalpricing_wage}):\\
 (a) if $\Pi_1(\bar{N})< \Pi_2(\bar{N})$, then the solution to  (\ref{optimalpricing_wage}) is given by that to (\ref{simplieidoptimalpricing_wage1}); \\
 (b) if $\Pi_2(\hat{N})< \Pi_1(\hat{N})$, then the solution to  (\ref{optimalpricing_wage}) is given by that to (\ref{simplieidoptimalpricing_wage2}); \\
 (c) if $\Pi_1(\bar{N})\geq \Pi_2(\bar{N})$ and $\Pi_2(\hat{N})\geq \Pi_1(\hat{N})$, then the solution to  (\ref{optimalpricing_wage}) is given by: 
\begin{equation}
\label{optimalpricing_wage_bothactive}
 \hspace{-4cm} \max_{p_f, N} \quad \lambda p_f-wN
\end{equation}
\begin{subnumcases}{\label{temp5}}
\lambda=\lambda_0\left[ 1-F_p\left(\dfrac{\alpha M}{\sqrt{N-\lambda / \mu}}+\beta p_f\right) \right] \label{temp6}\\
N= N_0 F_d(w) \label{supply_constraint_wage1}
\end{subnumcases}
\end{lemma}
The proof of Lemma \ref{wage_lemma}
 can be found in Appendix F. It provides a procedure to compute the solution to (\ref{optimalpricing_wage}): first compute $\bar{N}$ and $\hat{N}$ by solving (\ref{simplieidoptimalpricing_wage1}) and (\ref{simplieidoptimalpricing_wage2}) respectively, then identify which case of (a), (b), (c) applies, and obtain $N^*$ correspondingly.

{\bf Driver}: We first prove that $\nabla_+N^*(\tilde{w})>0$. Based on Lemma \ref{wage_lemma}, there are three cases: (a) $\Pi_1(\bar{N})< \Pi_2(\bar{N})$, (b) $\Pi_2(\hat{N})< \Pi_1(\hat{N})$,  (c) $\Pi_1(\bar{N})\geq \Pi_2(\bar{N})$ and $\Pi_2(\hat{N})\geq \Pi_1(\hat{N})$. In case (a), the constraint (\ref{minimum_wage}) is inactive at the optimal solution to (\ref{optimalpricing_wage}) (see proof of Lemma \ref{wage_lemma}). Therefore, case (a) corresponds to  $w\leq \tilde{w}$. When $w> \tilde{w}$, the minimum wage constraint is active, and either case (b) and (c) applies.

Note that in case (c), both (\ref{supply_constraint_wage}) and (\ref{minimum_wage}) are active at the optimal solution (see proof of Lemma \ref{wage_lemma}). Therefore, when $w$ increases, $N$ also increases. Then it suffices to show that there exists $\bar{w}>\tilde{w}$, such that the condition of case (c) holds when $\tilde{w}<w\leq \bar{w}$. 

Assume not, then there exists $\bar{w}'>\tilde{w}$ such that the conditions of case (b) holds when $\tilde{w}\leq w\leq \bar{w}'$. In this case, we have $\Pi_1(\bar{N})< \Pi_2(\bar{N})$ for $w\leq \tilde{w}$ and $\Pi_2(\hat{N})< \Pi_1(\hat{N})$ for $\tilde{w}\leq w\leq \bar{w}'$. Since $\Pi_2(N)$ is continuous with respect to $w$, let $w$ approaches $\tilde{w}$ from the let, then we have $\Pi_1(\tilde{N})\leq  \Pi_2(\tilde{N})$.  Let $w$ approaches $\tilde{w}$ from the right, then we have $\Pi_1(\tilde{N})\leq  \Pi_2(\tilde{N})$. This implies that $\Pi_1(\tilde{N})=  \Pi_2(\tilde{N})$.

In other words, the optimal solution to (\ref{simplieidoptimalpricing_wage1}) and (\ref{simplieidoptimalpricing_wage2}) are the same when $w=\tilde{w}$. This indicates that partial derivative of their objective functions with respect to $N$ at $N=N^*$ are the same, i.e., $\nabla g(N^*)=\tilde{w}$. Note that $w=\tilde{w}$ indicates $N^*=\tilde{N}$. Therefore, we have $\nabla g(\tilde{N})=\tilde{w}$. Since $x=g(N)$ is defined as $N=N_0F_d(x/N)$, applying implicit function theorem, we can obtain $\nabla g(N)$ and derive that $\nabla g(\tilde{N})=\tilde{w}$ is equivalent to
\begin{equation}
\dfrac{1+f_d(\tilde{\lambda}\tilde{p}_2/\tilde{N})N_0\tilde{\lambda}\tilde{p}_2/\tilde{N}^2}{f_d(\tilde{\lambda}\tilde{p}_2/\tilde{N})N_0/\tilde{N}}=\dfrac{\tilde{\lambda}\tilde{p}_2}{\tilde{N}},
\end{equation}
which is clearly impossible. A contradiction.

{\bf Passenger}: We now prove that $\nabla_+\lambda^*(\tilde{w})>0$. Consider $w$ such that $\tilde{w}<w<\bar{w}$. In this regime, both (\ref{supply_constraint_wage}) and (\ref{minimum_wage}) are active at the optimal solution. Therefore (\ref{optimalpricing_wage}) is equivalent to:
\begin{equation}
\label{optimalpricing_wage_proof}
 \hspace{-4cm} \max_{p_f\geq 0} \lambda p_f-Nw
\end{equation}
\begin{subnumcases}{\label{constraint_wage_proof}}
\lambda=\lambda_0\left[ 1-F_p\bigg(\alpha t_w(N-\lambda/\mu)+\beta p_f\bigg) \right] \label{demand_wage_proof}\\
N= N_0 F_d(w) \label{supply_wage_proof} 
\end{subnumcases}
We can plug (\ref{demand_cw}), (\ref{supply_cw}) and (\ref{pf_cw}) into the objective function of (\ref{optimalpricing_wage_proof}), transforming (\ref{optimalpricing_wage_proof}) to:
\begin{equation}
\label{optimalpricing_wage_final}
 \max_{c} \quad \lambda_0[1-F_p(c)]\cdot \dfrac{1}{\beta}\left[ c-\alpha t_w\bigg( N_0F_d(w)-\lambda_0[1-F_p(c)]/\mu  \bigg)\right] -N_0F_d(w)\cdot w.
\end{equation}
Note that $w$ is exogenous. The first order condition dictates that the derivative of (\ref{optimalpricing_wage_final}) with respect to $c$ is 0 at the optimal solution:
\begin{equation}
\label{firstordercond_proof}
\Phi(c,w)=-\lambda_0 f_p(c)\cdot \dfrac{1}{\beta}\left[ c-\alpha t_w(N_I)\right]+ \lambda_0[1-F_p(c)]\cdot \dfrac{1}{\beta} \left[ 1-\alpha  t_w'(N_I)\lambda_0f_p(c)/\mu \right]=0,
\end{equation}
where $N_I=N_0F_d(w)-\lambda_0[1-F_p(c)]/\mu$. Based on (\ref{demand_cw}), $\lambda$ is a decreasing function of $c$, thus it suffices to show that the positive partial derivative of $c$ with respect to $w$ is strictly negative, i.e.,  $\dfrac{\partial c}{\partial_+ w}<0$.  According to implicit function theorem, we have:
\[\dfrac{\partial c}{\partial_+ w}=-\dfrac{\partial 
\Phi}{\partial_+ w}/ \dfrac{\partial 
\Phi}{\partial c}.\]
If $c^*$ is local maximum for (\ref{optimalpricing_wage_final}), then $\Phi(c,w)>0$ in the neighborhood $c^*-\epsilon<c<c^*$ and $\Phi(c,w)<0$ in $c^*<c<c^*+\epsilon$. This indicates that $\dfrac{\partial 
\Phi}{\partial c^*}<0$, thus it suffices to show that $\dfrac{\partial 
\Phi}{\partial_+ w}<0$. We have:
\[\dfrac{\partial 
\Phi}{\partial_+ w}=\lambda_0 F_p(c)\dfrac{\alpha}{\beta}t_w'(N_I)N_0f_d(w)-\lambda_0[1-F_p(c)]\dfrac{\alpha \lambda_0}{\beta \mu} t_w''(N_I)f_p(c)N_0 f_d(w)<0.\]
This completes the proof.
\end{proof}
\subsection*{\bf{E: Proof of Proposition \ref{1storderconditionsufficient_cap}}}
\begin{proof}
To solve (\ref{optimalpricing_cap}), we first obtain the solution to  (\ref{optimalpricing}) and denote it as $(\tilde{p}_f,\tilde{p}_d,\tilde{\lambda},\tilde{N})$. If $\tilde{N}\leq N_{cap}$, then the cap constraint (\ref{max_cap}) is inactive, and (\ref{optimalpricing_cap}) reduces to  (\ref{optimalpricing}).

If $\tilde{N}> N_{cap}$, then (\ref{max_cap}) is active. In this case, (\ref{demand_constraint_cap}) and (\ref{supply_constraint_cap}) is equivalent to
\begin{subnumcases}{\label{simplypartial_cap}}
\lambda=\lambda_0 \left[ 1-e_p \left(\dfrac{\alpha M}{\sqrt{N_{cap}-\lambda / \mu}}+\beta p_f\right) \right]
\label{simplypartial_cap2} \\
N_{cap}=N_0e_d\dfrac{\lambda p_d}{N_{cap}}
\label{simplypartial_cap3}
\end{subnumcases} 
According to (\ref{simplypartial_cap3}), $\lambda p_d=\dfrac{N_{cap}^2}{N_0 e_d}$. Therefore, the objective function is $R=\lambda p_f-\dfrac{N_{cap}^2}{N_0 e_d}$. Apply implicit function theorem on (\ref{simplypartial_cap2}), we can obtain $\dfrac{\partial \lambda}{\partial p_f}$ and further $\dfrac{\partial R}{\partial p_f}$. The first order condition becomes
\begin{equation}
\label{temp4}
\lambda_0e_p\beta\mu p_f=\lambda_0 e_p \alpha M \dfrac{\lambda}{2(N_{cap}-\lambda/\mu)^{3/2}}+\lambda \mu.
\end{equation}

{\bf Uniqueness}: we first show that   there is at most one set of $\lambda>0$, $p_f> 0$ and  $p_d> 0$  that satisfy (\ref{simplypartial_cap2}), (\ref{simplypartial_cap3}), and (\ref{temp4}). Assume not, i.e., $(\lambda,p_f,p_d)$ and $(\lambda',p_f',p_d')$ are both solutions to (\ref{simplypartial_cap2}), (\ref{simplypartial_cap3}), and (\ref{temp4}). If $p_f=p_f'$, then it is easy to verify that $p_d=p_d'$ and $\lambda'=\lambda$. Therefore, we consider $p_f'>p_f$ without loss of generality. Based on (\ref{simplypartial_cap2}), we have $\lambda'<\lambda$. However, (\ref{temp4}) dictates that $\lambda'>\lambda$, A contradiction. 

{\bf Existence and Optimality}:  Note that for any $p_f\geq 0$ such that
\begin{equation}
\label{boundaryofpf}
e_p\left(\dfrac{\alpha M}{\sqrt{N_{cap}}}+\beta p_f\right)<1.
\end{equation}
There is a unique $\lambda>0$ that satisfies (\ref{simplypartial_cap2}). This is because the righthand side of (\ref{simplypartial_cap2}) is a concave function $\lambda$ which has a unique intersection with the left-hand side of (\ref{simplypartial_cap2}). We can therefore view $\lambda$ as a function of $p_f$ determined by (\ref{simplypartial_cap2}), denoted by $\lambda=h_1(p_f)$. Apply implicit function theorem on (\ref{simplypartial_cap2}), we have $\dfrac{\partial \lambda}{\partial p_f}<0$. Therefore, it is a decreasing function such that $h_1(0)>0$, and it can be verified that   $h_1(p_f)\rightarrow 0$ for  sufficiently large $p_f$ that satisfies (\ref{boundaryofpf}).
On the other hand, (\ref{temp4}) prescribes $\lambda$ as an increasing function of $p_f$. We denote it as as $\lambda=h_2(p_f)$, and we have $h_2(0)=0$. To prove existence, it suffices to show that $h_1(p_f)$ intersects with $h_2(p_f)$ at some $p_f^*>0$, which is clearly true. 

$p_f^*$ is either local minimum or local maximum. Note that $p_f$ is bounded between $0$ and an upper bound determined by (\ref{boundaryofpf}). On the  boundary the objective value $R$ is smaller than that in the interior where $\lambda>0$ and $p_f>0$ (since the revenue $\lambda p_f$ on the boundary is $0$). This indicates that $p_f^*$ has to be local maximum. Since there is a unique solution to the first order conditions, $p_f^*$ is the globally optimal solution. This completes the proof. 

\end{proof}

\subsection*{\bf{F: Proof of Lemma \ref{wage_lemma}}}
\begin{proof}
Let $(p_f^*,\lambda^*,N^*)$, $(\bar{p}_f,\bar{\lambda},\bar{N})$ and $(\hat{p}_f,\hat{\lambda},\hat{N})$ be the optimal solution (\ref{simplieidoptimalpricing_wage}), (\ref{simplieidoptimalpricing_wage1}) and (\ref{simplieidoptimalpricing_wage2}), respectively. Define $R^*$, $\bar{R}$ and $\hat{R}$ as the corresponding optimal values. Note that $R^*\leq \bar{R}$ and $R^*\leq \hat{R}$. This is because the objective value of (\ref{simplieidoptimalpricing_wage}) is smaller than that of (\ref{simplieidoptimalpricing_wage1}) and (\ref{simplieidoptimalpricing_wage2}). We consider the following three cases:
\begin{itemize}
\item $\Pi_1(\bar{N})< \Pi_2(\bar{N})$. Given $N$, the optimization problems (\ref{simplieidoptimalpricing_wage1}) and (\ref{simplieidoptimalpricing_wage2}) are equivalent (since the second term of the objective functions are constants). Therefore, the corresponding optimal solution (i.e., $\lambda$ and $p_f$) are the same when $N$ are the same, and $\Pi_1(N)- \Pi_2(N)=wN-g(N)$ for any $N$. This indicates that $g(\bar{N})>w\bar{N}$. It further implies that the objective value of  (\ref{simplieidoptimalpricing_wage}) can attain $\bar{R}$ when $N=\bar{N}$. Since $R^*$ is upper bounded by $\bar{R}$, we conclude that $N^*=\bar{N}$. 
\item $\Pi_2(\hat{N})< \Pi_1(\hat{N})$. In this case we have $N^*=\hat{N}$. Proof is the same as case (a).
\item $\Pi_1(\bar{N})\geq \Pi_2(\bar{N})$ and $\Pi_2(\hat{N})\geq \Pi_1(\hat{N})$. In this case, we can show that $g(N^*)=wN^*$. Assume not, then without loss of generality,  consider the case where $g(N^*)>wN^*$. Consider the following problem:
\begin{equation}
\label{temp15}
 \hspace{-4cm} \max_{p_f, N} \quad \lambda p_f-g(N)
\end{equation}
\begin{subnumcases}{\label{temp16}}
\lambda=\lambda_0\left[ 1-F_p\left(\dfrac{\alpha M}{\sqrt{N-\lambda / \mu}}+\beta p_f\right) \right] \label{temp16_1}\\
g(N)\geq wN \label{temp16_2}
\end{subnumcases}
We conclude that $N^*$ is the optimal solution to (\ref{temp15}). This is because if not, then there exists another solution that satisfies (\ref{temp16}) and obtains a higher value than $R^*$, which contradicts with the fact that $N^*$ is optimal solution to (\ref{simplieidoptimalpricing_wage}). Since $g(N^*)>wN^*$, the constraint (\ref{temp16_2}) is inactive at the optimal solution to (\ref{temp15}), therefore, it is equivalent to (\ref{simplieidoptimalpricing_wage1}), i.e., $N^*=\bar{N}$. In this case, $g(N^*)>wN^*$ implies $g(\bar{N})>w\bar{N}$. Since $\Pi_1(N)- \Pi_2(N)=wN-g(N)$ for any $N$, we have  $\Pi_1(\bar{N})< \Pi_2(\bar{N})$. This contradicts with the assumption that $\Pi_1(\bar{N})\geq  \Pi_2(\bar{N})$. Therefore, $g(N^*)=wN^*$. This indicates that both (\ref{supply_constraint_wage}) and (\ref{minimum_wage}) are active at the optimal solution to (\ref{optimalpricing_wage}), which leads to (\ref{optimalpricing_wage_bothactive}).

This completes the proof.
\end{itemize}
\end{proof}

\end{document}

%% file: format.tex
\usepackage{epsfig}
\usepackage{subfigure}
\usepackage{amssymb,amsmath,amsfonts,layout,graphicx}
\usepackage{makeidx}
\usepackage{babel}
\usepackage{tikz}
\usepackage{sublabel}
\usepackage{cases}
\usepackage{algorithm}
\usepackage{algorithmic}

\usepackage
[
        letterpaper,
        left=1.91cm,
        right=1.91cm,
        top=1.91cm,
        bottom=2cm,
]
{geometry}

\newtheorem{lemma}{Lemma}
\newtheorem{assumption}{Assumption}
\newtheorem{proposition}{Proposition}
\newtheorem{theorem}{Theorem}

\newtheorem{corollary}{Corollary}
\newtheorem{remark}{Remark}

\newcommand{\x}{{\bf x}}

%% file: figure1_unreg.tex
%
%
\begin{tikzpicture}

\pgfplotsset{every axis/.append style={
font=\small,
thin,
tick style={ultra thin}}}
\pgfplotsset{every axis y label/.style={
at={(-0.48,0.5)},
xshift=32pt,
rotate=90}}

\begin{axis}[%
width=1.794in,
height=1.03in,
at={(1.358in,0.0in)},
scale only axis,
xmin=1000,
xmax=2000,
xtick={1000,1250,1500,1750,2000},
xticklabels={{1K},{1.25K},{1.5K},{1.75K},{2K}},
xlabel style={font=\color{white!15!black}},
xlabel={Potential Passenger /min},
ymin=4000,
ymax=6000,
ytick={4000,4500,5000,5500,6000},
yticklabels={{4K},{4.5K},{5K},{5.5K},{6K}},
ylabel style={font=\color{white!15!black}},
ylabel={Number of Drivers},
axis background/.style={fill=white},
legend style={legend cell align=left, align=left, draw=white!15!black}
]
\addplot [color=black, line width=1.0pt]
  table[row sep=crcr]{%
1000	4262.62595441961\\
1010.10101010101	4282.14313549826\\
1020.20202020202	4301.50393564711\\
1030.30303030303	4320.71046500294\\
1040.40404040404	4339.76478948349\\
1050.50505050505	4358.6689321383\\
1060.60606060606	4377.42487444381\\
1070.70707070707	4396.03455754549\\
1080.80808080808	4414.49988344986\\
1090.90909090909	4432.82271616885\\
1101.0101010101	4451.00488281899\\
1111.11111111111	4469.04817467756\\
1121.21212121212	4486.95434819799\\
1131.31313131313	4504.72512598628\\
1141.41414141414	4522.36219774059\\
1151.51515151515	4539.86722115549\\
1161.61616161616	4557.24182279277\\
1171.71717171717	4574.48759892019\\
1181.81818181818	4591.60611631987\\
1191.91919191919	4608.59891306747\\
1202.0202020202	4625.46749928369\\
1212.12121212121	4642.21335785927\\
1222.22222222222	4658.83794515464\\
1232.32323232323	4675.34269167533\\
1242.42424242424	4691.72900272436\\
1252.52525252525	4707.99825903239\\
1262.62626262626	4724.15181736671\\
1272.72727272727	4740.19101111997\\
1282.82828282828	4756.11715087948\\
1292.92929292929	4771.93152497792\\
1303.0303030303	4787.63540002612\\
1313.13131313131	4803.2300214289\\
1323.23232323232	4818.7166138843\\
1333.33333333333	4834.09638186727\\
1343.43434343434	4849.37051009806\\
1353.53535353535	4864.54016399622\\
1363.63636363636	4879.60649012059\\
1373.73737373737	4894.57061659595\\
1383.83838383838	4909.43365352668\\
1393.93939393939	4924.19669339814\\
1404.0404040404	4938.86081146611\\
1414.14141414141	4953.42706613468\\
1424.24242424242	4967.8964993232\\
1434.34343434343	4982.27013682254\\
1444.44444444444	4996.54898864116\\
1454.54545454545	5010.73404934123\\
1464.64646464646	5024.82629836532\\
1474.74747474747	5038.82670035391\\
1484.84848484848	5052.73620545401\\
1494.94949494949	5066.5557496193\\
1505.05050505051	5080.286254902\\
1515.15151515152	5093.92862973685\\
1525.25252525253	5107.48376921738\\
1535.35353535354	5120.95255536474\\
1545.45454545455	5134.33585738947\\
1555.55555555556	5147.63453194626\\
1565.65656565657	5160.84942338208\\
1575.75757575758	5173.98136397783\\
1585.85858585859	5187.03117418377\\
1595.9595959596	5199.99966284884\\
1606.06060606061	5212.88762744422\\
1616.16161616162	5225.69585428115\\
1626.26262626263	5238.42511872336\\
1636.36363636364	5251.07618539409\\
1646.46464646465	5263.64980837806\\
1656.56565656566	5276.14673141846\\
1666.66666666667	5288.56768810905\\
1676.76767676768	5300.91340208169\\
1686.86868686869	5313.18458718922\\
1696.9696969697	5325.38194768411\\
1707.07070707071	5337.50617839265\\
1717.17171717172	5349.55796488525\\
1727.27272727273	5361.5379836425\\
1737.37373737374	5373.44690221753\\
1747.47474747475	5385.28537939445\\
1757.57575757576	5397.05406534328\\
1767.67676767677	5408.75360177118\\
1777.77777777778	5420.38462207032\\
1787.87878787879	5431.94775146235\\
1797.9797979798	5443.44360713968\\
1808.08080808081	5454.87279840346\\
1818.18181818182	5466.23592679868\\
1828.28282828283	5477.53358624606\\
1838.38383838384	5488.76636317122\\
1848.48484848485	5499.93483663094\\
1858.58585858586	5511.03957843663\\
1868.68686868687	5522.08115327523\\
1878.78787878788	5533.06011882742\\
1888.88888888889	5543.97702588337\\
1898.9898989899	5554.83241845593\\
1909.09090909091	5565.62683389154\\
1919.19191919192	5576.36080297873\\
1929.29292929293	5587.03485005437\\
1939.39393939394	5597.64949310767\\
1949.49494949495	5608.20524388212\\
1959.59595959596	5618.70260797522\\
1969.69696969697	5629.14208493623\\
1979.79797979798	5639.52416836192\\
1989.89898989899	5649.84934599035\\
2000	5660.11809979279\\
};

\end{axis}
\end{tikzpicture}%

%% file: figure2_unreg.tex
%
%
\begin{tikzpicture}
\pgfplotsset{every axis y label/.style={
at={(-0.44,0.5)},
xshift=32pt,
rotate=90}}

\begin{axis}[%
width=1.794in,
height=1.03in,
at={(1.358in,0.0in)},
scale only axis,
xmin=1000,
xmax=2000,
xtick={1000,1250,1500,1750,2000},
xticklabels={{1K},{1.25K},{1.5K},{1.75K},{2K}},
xlabel style={font=\color{white!15!black}},
xlabel={Potential Passenger/min},
ymin=140,
ymax=220,
ylabel style={font=\color{white!15!black}},
ylabel={Ride Arrival/min},
axis background/.style={fill=white},
legend style={legend cell align=left, align=left, draw=white!15!black}
]
\addplot [color=black, line width=1.0pt]
  table[row sep=crcr]{%
000	143.66800749413\\
1010.10101010101	144.675248582715\\
1020.20202020202	145.675622163584\\
1030.30303030303	146.669197484996\\
1040.40404040404	147.656042897234\\
1050.50505050505	148.636225865706\\
1060.60606060606	149.609812983908\\
1070.70707070707	150.576869986222\\
1080.80808080808	151.537461760567\\
1090.90909090909	152.491652360879\\
1101.0101010101	153.439505019439\\
1111.11111111111	154.381082159035\\
1121.21212121212	155.316445404955\\
1131.31313131313	156.245655596826\\
1141.41414141414	157.168772800277\\
1151.51515151515	158.085856318441\\
1161.61616161616	158.996964703293\\
1171.71717171717	159.90215576682\\
1181.81818181818	160.801486592032\\
1191.91919191919	161.6950135438\\
1202.0202020202	162.58279227954\\
1212.12121212121	163.464877759731\\
1222.22222222222	164.341324258273\\
1232.32323232323	165.212185372689\\
1242.42424242424	166.07751403416\\
1252.52525252525	166.937362517419\\
1262.62626262626	167.79178245047\\
1272.72727272727	168.640824824177\\
1282.82828282828	169.484540001681\\
1292.92929292929	170.32297772768\\
1303.0303030303	171.156187137556\\
1313.13131313131	171.984216766359\\
1323.23232323232	172.807114557645\\
1333.33333333333	173.62492787217\\
1343.43434343434	174.437703496452\\
1353.53535353535	175.245487651183\\
1363.63636363636	176.048325999514\\
1373.73737373737	176.846263655202\\
1383.83838383838	177.639345190626\\
1393.93939393939	178.427614644671\\
1404.0404040404	179.211115530486\\
1414.14141414141	179.989890843113\\
1424.24242424242	180.763983066993\\
1434.34343434343	181.533434183353\\
1444.44444444444	182.298285677465\\
1454.54545454545	183.058578545791\\
1464.64646464646	183.814353303016\\
1474.74747474747	184.565649988955\\
1484.84848484848	185.312508175357\\
1494.94949494949	186.054966972594\\
1505.05050505051	186.793065036243\\
1515.15151515152	187.526840573554\\
1525.25252525253	188.256331349825\\
1535.35353535354	188.981574694659\\
1545.45454545455	189.702607508132\\
1555.55555555556	190.41946626685\\
1565.65656565657	191.132187029914\\
1575.75757575758	191.840805444789\\
1585.85858585859	192.545356753071\\
1595.9595959596	193.245875796169\\
1606.06060606061	193.942397020885\\
1616.16161616162	194.634954484918\\
1626.26262626263	195.323581862265\\
1636.36363636364	196.00831244854\\
1646.46464646465	196.689179166212\\
1656.56565656566	197.366214569754\\
1666.66666666667	198.039450850707\\
1676.76767676768	198.70891984267\\
1686.86868686869	199.374653026202\\
1696.9696969697	200.036681533655\\
1707.07070707071	200.695036153919\\
1717.17171717172	201.349747337098\\
1727.27272727273	202.000845199116\\
1737.37373737374	202.648359526235\\
1747.47474747475	203.29231977952\\
1757.57575757576	203.932755099217\\
1767.67676767677	204.569694309072\\
1777.77777777778	205.203165920579\\
1787.87878787879	205.833198137157\\
1797.9797979798	206.459818858272\\
1808.08080808081	207.08305568348\\
1818.18181818182	207.702935916418\\
1828.28282828283	208.31948656873\\
1838.38383838384	208.932734363925\\
1848.48484848485	209.542705741188\\
1858.58585858586	210.149426859116\\
1868.68686868687	210.75292359941\\
1878.78787878788	211.353221570504\\
1888.88888888889	211.950346111131\\
1898.9898989899	212.544322293852\\
1909.09090909091	213.135174928509\\
1919.19191919192	213.722928565643\\
1929.29292929293	214.307607499847\\
1939.39393939394	214.889235773079\\
1949.49494949495	215.46783717791\\
1959.59595959596	216.043435260739\\
1969.69696969697	216.616053324946\\
1979.79797979798	217.185714434004\\
1989.89898989899	217.752441414542\\
2000	218.31625685936\\
};

\end{axis}
\end{tikzpicture}%

%% file: figure3_unreg.tex
%
%
\begin{tikzpicture}
\pgfplotsset{every axis y label/.style={
at={(-0.46,0.5)},
xshift=32pt,
rotate=90}}
\begin{axis}[%
width=1.794in,
height=1.03in,
at={(1.358in,0.0in)},
scale only axis,
xmin=1000,
xmax=2000,
xtick={1000,1250,1500,1750,2000},
xticklabels={{1K},{1.25K},{1.5K},{1.75K},{2K}},
xlabel style={font=\color{white!15!black}},
xlabel={Potential Passegner/minute},
ymin=0.54,
ymax=0.63,
ylabel style={font=\color{white!15!black}},
ylabel={Occupancy Rate},
axis background/.style={fill=white},
legend style={legend cell align=left, align=left, draw=white!15!black}
]
\addplot [color=black, line width=1.0pt]
  table[row sep=crcr]{%
1000	0.549376967905495\\
1010.10101010101	0.550707082243261\\
1020.20202020202	0.552019172082705\\
1030.30303030303	0.553313613205465\\
1040.40404040404	0.554590770692752\\
1050.50505050505	0.5558509993147\\
1060.60606060606	0.55709464390238\\
1070.70707070707	0.558322039703398\\
1080.80808080808	0.559533512721929\\
1090.90909090909	0.560729380044003\\
1101.0101010101	0.561909950148796\\
1111.11111111111	0.563075523206644\\
1121.21212121212	0.564226391364448\\
1131.31313131313	0.565362839019099\\
1141.41414141414	0.566485143079525\\
1151.51515151515	0.567593573217926\\
1161.61616161616	0.568688392110703\\
1171.71717171717	0.569769855669608\\
1181.81818181818	0.570838213263573\\
1191.91919191919	0.571893707931654\\
1202.0202020202	0.572936576587535\\
1212.12121212121	0.573967050215959\\
1222.22222222222	0.574985354061492\\
1232.32323232323	0.575991707809947\\
1242.42424242424	0.576986325762825\\
1252.52525252525	0.577969417005088\\
1262.62626262626	0.578941185566552\\
1272.72727272727	0.579901830577207\\
1282.82828282828	0.580851546416715\\
1292.92929292929	0.581790522858359\\
1303.0303030303	0.582718945207679\\
1313.13131313131	0.583636994436027\\
1323.23232323232	0.584544847309263\\
1333.33333333333	0.585442676511798\\
1343.43434343434	0.586330650766191\\
1353.53535353535	0.587208934948471\\
1363.63636363636	0.588077690199391\\
1373.73737373737	0.588937074031749\\
1383.83838383838	0.589787240433979\\
1393.93939393939	0.590628339970129\\
1404.0404040404	0.5914605198764\\
1414.14141414141	0.592283924154375\\
1424.24242424242	0.59309869366107\\
1434.34343434343	0.593904966195945\\
1444.44444444444	0.59470287658498\\
1454.54545454545	0.595492556761956\\
1464.64646464646	0.596274135847019\\
1474.74747474747	0.59704774022267\\
1484.84848484848	0.597813493607252\\
1494.94949494949	0.59857151712604\\
1505.05050505051	0.599321929380038\\
1515.15151515152	0.600064846512551\\
1525.25252525253	0.600800382273627\\
1535.35353535354	0.601528648082455\\
1545.45454545455	0.602249753087782\\
1555.55555555556	0.602963804226429\\
1565.65656565657	0.603670906279986\\
1575.75757575758	0.604371161929733\\
1585.85858585859	0.605064671809869\\
1595.9595959596	0.605751534559112\\
1606.06060606061	0.606431846870702\\
1616.16161616162	0.607105703540909\\
1626.26262626263	0.607773197516055\\
1636.36363636364	0.608434419938134\\
1646.46464646465	0.609089460189064\\
1656.56565656566	0.609738405933624\\
1666.66666666667	0.610381343161124\\
1676.76767676768	0.611018356225848\\
1686.86868686869	0.611649527886309\\
1696.9696969697	0.612274939343373\\
1707.07070707071	0.612894670277273\\
1717.17171717172	0.613508798883554\\
1727.27272727273	0.614117401907999\\
1737.37373737374	0.614720554680549\\
1747.47474747475	0.615318331148272\\
1757.57575757576	0.615910803907392\\
1767.67676767677	0.616498044234433\\
1777.77777777778	0.617080122116478\\
1787.87878787879	0.617657106280603\\
1797.9797979798	0.618229064222485\\
1808.08080808081	0.618796062234238\\
1818.18181818182	0.61935816543147\\
1828.28282828283	0.61991543777962\\
1838.38383838384	0.620467942119573\\
1848.48484848485	0.621015740192588\\
1858.58585858586	0.621558892664552\\
1868.68686868687	0.622097459149596\\
1878.78787878788	0.622631498233077\\
1888.88888888889	0.623161067493955\\
1898.9898989899	0.623686223526579\\
1909.09090909091	0.624207021961904\\
1919.19191919192	0.62472351748816\\
1929.29292929293	0.625235763870977\\
1939.39393939394	0.625743813972992\\
1949.49494949495	0.626247719772959\\
1959.59595959596	0.626747532384361\\
1969.69696969697	0.627243302073555\\
1979.79797979798	0.627735078277455\\
1989.89898989899	0.628222909620764\\
2000	0.628706843932786\\
};

\end{axis}
\end{tikzpicture}%

%% file: figure4_unreg.tex
%
%
\definecolor{mycolor1}{rgb}{0.63529,0.07843,0.18431}%
\begin{tikzpicture}

\pgfplotsset{every axis/.append style={
font=\small,
thin,
tick style={ultra thin}}}
\pgfplotsset{every axis y label/.style={
at={(-0.42,0.5)},
xshift=32pt,
rotate=90}}

\begin{axis}[%
width=1.794in,
height=1.03in,
at={(1.358in,0.0in)},
scale only axis,
xmin=1000,
xmax=2000,
xtick={1000,1250,1500,1750,2000},
xticklabels={{1K},{1.25K},{1.5K},{1.75K},{2K}},
xlabel style={font=\color{white!15!black}},
xlabel={Potential Passenger/min},
ymin=9.02145411203814,
ymax=18.0214541120381,
ylabel style={font=\color{white!15!black}},
ylabel={Price (\$/trip)},
axis background/.style={fill=white},
legend style={at={(0.65,0.4)}, anchor=south west, legend cell align=left, align=left, draw=white!15!black}
]
\addplot [color=black, line width=1.0pt]
  table[row sep=crcr]{%
1000	15.7592639701024\\
1010.10101010101	15.7896063107345\\
1020.20202020202	15.8196464741895\\
1030.30303030303	15.8493899030915\\
1040.40404040404	15.8788418868397\\
1050.50505050505	15.9080075675807\\
1060.60606060606	15.9368919458852\\
1070.70707070707	15.9654998861478\\
1080.80808080808	15.9938361217244\\
1090.90909090909	16.0219052598248\\
1101.0101010101	16.0497117861724\\
1111.11111111111	16.0772600694454\\
1121.21212121212	16.1045543655112\\
1131.31313131313	16.1315988214669\\
1141.41414141414	16.1583974794945\\
1151.51515151515	16.1849542805423\\
1161.61616161616	16.2112730678428\\
1171.71717171717	16.2373575902729\\
1181.81818181818	16.2632115055686\\
1191.91919191919	16.288838383399\\
1202.0202020202	16.3142417083077\\
1212.12121212121	16.3394248825301\\
1222.22222222222	16.3643912286899\\
1232.32323232323	16.3891439923842\\
1242.42424242424	16.4136863446603\\
1252.52525252525	16.4380213843911\\
1262.62626262626	16.4621521405538\\
1272.72727272727	16.4860815744161\\
1282.82828282828	16.5098125816353\\
1292.92929292929	16.533347994274\\
1303.0303030303	16.556690582736\\
1313.13131313131	16.5798430576277\\
1323.23232323232	16.602808071546\\
1333.33333333333	16.6255882207995\\
1343.43434343434	16.6481860470624\\
1353.53535353535	16.6706040389676\\
1363.63636363636	16.692844633639\\
1373.73737373737	16.7149102181678\\
1383.83838383838	16.7368031310339\\
1393.93939393939	16.7585256634762\\
1404.0404040404	16.7800800608124\\
1414.14141414141	16.8014685237122\\
1424.24242424242	16.8226932094247\\
1434.34343434343	16.8437562329625\\
1444.44444444444	16.8646596682445\\
1454.54545454545	16.8854055491989\\
1464.64646464646	16.9059958708281\\
1474.74747474747	16.9264325902368\\
1484.84848484848	16.946717627626\\
1494.94949494949	16.9668528672524\\
1505.05050505051	16.9868401583571\\
1515.15151515152	17.0066813160617\\
1525.25252525253	17.0263781222367\\
1535.35353535354	17.04593232634\\
1545.45454545455	17.065345646229\\
1555.55555555556	17.0846197689464\\
1565.65656565657	17.1037563514807\\
1575.75757575758	17.1227570215026\\
1585.85858585859	17.1416233780782\\
1595.9595959596	17.1603569923597\\
1606.06060606061	17.1789594082553\\
1616.16161616162	17.1974321430774\\
1626.26262626263	17.2157766881715\\
1636.36363636364	17.2339945095264\\
1646.46464646465	17.2520870483649\\
1656.56565656566	17.2700557217176\\
1666.66666666667	17.2879019229793\\
1676.76767676768	17.3056270224491\\
1686.86868686869	17.3232323678543\\
1696.9696969697	17.3407192848594\\
1707.07070707071	17.3580890775598\\
1717.17171717172	17.3753430289621\\
1727.27272727273	17.3924824014496\\
1737.37373737374	17.4095084372355\\
1747.47474747475	17.426422358803\\
1757.57575757576	17.4432253693324\\
1767.67676767677	17.4599186531175\\
1777.77777777778	17.4765033759695\\
1787.87878787879	17.4929806856102\\
1797.9797979798	17.5093517120542\\
1808.08080808081	17.5256175679812\\
1818.18181818182	17.5417793490973\\
1828.28282828283	17.557838134488\\
1838.38383838384	17.5737949869605\\
1848.48484848485	17.5896509533772\\
1858.58585858586	17.6054070649814\\
1868.68686868687	17.621064337713\\
1878.78787878788	17.636623772517\\
1888.88888888889	17.6520863556437\\
1898.9898989899	17.6674530589411\\
1909.09090909091	17.68272484014\\
1919.19191919192	17.6979026431318\\
1929.29292929293	17.7129873982392\\
1939.39393939394	17.7279800224802\\
1949.49494949495	17.7428814198253\\
1959.59595959596	17.7576924814488\\
1969.69696969697	17.7724140859733\\
1979.79797979798	17.7870470997088\\
1989.89898989899	17.8015923768852\\
2000	17.8160507598801\\
};
\addlegendentry{$p_f$}

\addplot [color=mycolor1, dashed, line width=1.0pt]
  table[row sep=crcr]{%
1000	9.36031951120671\\
1010.10101010101	9.38046602738736\\
1020.20202020202	9.40048065006858\\
1030.30303030303	9.42036447637654\\
1040.40404040404	9.44011860003235\\
1050.50505050505	9.45974411081057\\
1060.60606060606	9.47924209404418\\
1070.70707070707	9.49861363017272\\
1080.80808080808	9.51785979433023\\
1090.90909090909	9.53698165597012\\
1101.0101010101	9.5559802785242\\
1111.11111111111	9.57485671909331\\
1121.21212121212	9.59361202816734\\
1131.31313131313	9.61224724937234\\
1141.41414141414	9.63076341924277\\
1151.51515151515	9.64916156701714\\
1161.61616161616	9.66744271445518\\
1171.71717171717	9.68560787567506\\
1181.81818181818	9.70365805700917\\
1191.91919191919	9.72159425687715\\
1202.0202020202	9.73941746567473\\
1212.12121212121	9.75712866567747\\
1222.22222222222	9.77472883095816\\
1232.32323232323	9.79221892731686\\
1242.42424242424	9.80959991222273\\
1252.52525252525	9.82687273476671\\
1262.62626262626	9.84403833562431\\
1272.72727272727	9.86109764702772\\
1282.82828282828	9.87805159274648\\
1292.92929292929	9.89490108807625\\
1303.0303030303	9.91164703983483\\
1313.13131313131	9.92829034636509\\
1323.23232323232	9.94483189754413\\
1333.33333333333	9.96127257479828\\
1343.43434343434	9.97761325112346\\
1353.53535353535	9.99385479111042\\
1363.63636363636	10.0099980509746\\
1373.73737373737	10.02604387859\\
1383.83838383838	10.0419931135271\\
1393.93939393939	10.0578465870944\\
1404.0404040404	10.0736051223823\\
1414.14141414141	10.0892695343114\\
1424.24242424242	10.1048406296821\\
1434.34343434343	10.1203192072275\\
1444.44444444444	10.1357060576681\\
1454.54545454545	10.1510019637688\\
1464.64646464646	10.1662077003974\\
1474.74747474747	10.181324034585\\
1484.84848484848	10.1963517255877\\
1494.94949494949	10.2112915249499\\
1505.05050505051	10.2261441765683\\
1515.15151515152	10.2409104167575\\
1525.25252525253	10.2555909743156\\
1535.35353535354	10.270186570592\\
1545.45454545455	10.2846979195544\\
1555.55555555556	10.2991257278572\\
1565.65656565657	10.3134706949104\\
1575.75757575758	10.3277335129483\\
1585.85858585859	10.3419148670992\\
1595.9595959596	10.3560154354551\\
1606.06060606061	10.370035889141\\
1616.16161616162	10.3839768923852\\
1626.26262626263	10.3978391025891\\
1636.36363636364	10.4116231703969\\
1646.46464646465	10.4253297397655\\
1656.56565656566	10.4389594480342\\
1666.66666666667	10.4525129259943\\
1676.76767676768	10.4659907979579\\
1686.86868686869	10.4793936818274\\
1696.9696969697	10.4927221891641\\
1707.07070707071	10.5059769252563\\
1717.17171717172	10.5191584891878\\
1727.27272727273	10.5322674739052\\
1737.37373737374	10.5453044662853\\
1747.47474747475	10.5582700472021\\
1757.57575757576	10.5711647915927\\
1767.67676767677	10.5839892685237\\
1777.77777777778	10.5967440412564\\
1787.87878787879	10.6094296673116\\
1797.9797979798	10.6220466985342\\
1808.08080808081	10.634595681157\\
1818.18181818182	10.647077155864\\
1828.28282828283	10.6594916578531\\
1838.38383838384	10.6718397168983\\
1848.48484848485	10.6841218574116\\
1858.58585858586	10.6963385985035\\
1868.68686868687	10.7084904540442\\
1878.78787878788	10.7205779327228\\
1888.88888888889	10.732601538107\\
1898.9898989899	10.7445617687017\\
1909.09090909091	10.7564591180068\\
1919.19191919192	10.7682940745751\\
1929.29292929293	10.7800671220692\\
1939.39393939394	10.7917787393172\\
1949.49494949495	10.803429400369\\
1959.59595959596	10.8150195745511\\
1969.69696969697	10.8265497265211\\
1979.79797979798	10.8380203163213\\
1989.89898989899	10.8494317994326\\
2000	10.8607846268263\\
};
\addlegendentry{$p_d$}

\end{axis}
\end{tikzpicture}%

%% file: figure5_unreg.tex
%
%
\begin{tikzpicture}

\begin{axis}[%
width=1.794in,
height=1.03in,
at={(1.358in,0.0in)},
scale only axis,
xmin=1000,
xmax=2000,
xtick={1000,1250,1500,1750,2000},
xticklabels={{1K},{1.25K},{1.5K},{1.75K},{2K}},
xlabel style={font=\color{white!15!black}},
xlabel={Potential Passengers/min},
ymin=18,
ymax=26,
ylabel style={font=\color{white!15!black}},
ylabel={Driver Wage/hour},
axis background/.style={fill=white},
legend style={legend cell align=left, align=left, draw=white!15!black}
]
\addplot [color=black, line width=1.0pt]
  table[row sep=crcr]{%
1000	18.9288734418162\\
1010.10101010101	19.0155426111332\\
1020.20202020202	19.1015173458786\\
1030.30303030303	19.1868070164527\\
1040.40404040404	19.2714207968942\\
1050.50505050505	19.3553676708788\\
1060.60606060606	19.4386564374699\\
1070.70707070707	19.5212957166351\\
1080.80808080808	19.6032939545389\\
1090.90909090909	19.6846594286251\\
1101.0101010101	19.7654002524972\\
1111.11111111111	19.845524380609\\
1121.21212121212	19.9250396127734\\
1131.31313131313	20.0039535984982\\
1141.41414141414	20.0822738411588\\
1151.51515151515	20.1600077020135\\
1161.61616161616	20.2371624040706\\
1171.71717171717	20.3137450358128\\
1181.81818181818	20.389762554787\\
1191.91919191919	20.4652217910646\\
1202.0202020202	20.5401294505786\\
1212.12121212121	20.6144921183429\\
1222.22222222222	20.6883162615591\\
1232.32323232323	20.7616082326164\\
1242.42424242424	20.8343742719879\\
1252.52525252525	20.9066205110294\\
1262.62626262626	20.9783529746832\\
1272.72727272727	21.0495775840931\\
1282.82828282828	21.1203001591321\\
1292.92929292929	21.1905264208477\\
1303.0303030303	21.2602619938283\\
1313.13131313131	21.3295124084935\\
1323.23232323232	21.3982831033112\\
1333.33333333333	21.4665794269453\\
1343.43434343434	21.5344066403359\\
1353.53535353535	21.6017699187152\\
1363.63636363636	21.6686743535615\\
1373.73737373737	21.7351249544934\\
1383.83838383838	21.8011266511075\\
1393.93939393939	21.8666842947599\\
1404.0404040404	21.9318026602962\\
1414.14141414141	21.9964864477293\\
1424.24242424242	22.0607402838691\\
1434.34343434343	22.1245687239041\\
1444.44444444444	22.1879762529389\\
1454.54545454545	22.2509672874859\\
1464.64646464646	22.313546176917\\
1474.74747474747	22.3757172048726\\
1484.84848484848	22.4374845906328\\
1494.94949494949	22.4988524904498\\
1505.05050505051	22.5598249988444\\
1515.15151515152	22.6204061498663\\
1525.25252525253	22.6805999183219\\
1535.35353535354	22.7404102209673\\
1545.45454545455	22.7998409176714\\
1555.55555555556	22.8588958125463\\
1565.65656565657	22.9175786550501\\
1575.75757575758	22.9758931410595\\
1585.85858585859	23.0338429139154\\
1595.9595959596	23.0914315654409\\
1606.06060606061	23.1486626369346\\
1616.16161616162	23.2055396201377\\
1626.26262626263	23.2620659581768\\
1636.36363636364	23.3182450464841\\
1646.46464646465	23.3740802336935\\
1656.56565656566	23.4295748225158\\
1666.66666666667	23.4847320705917\\
1676.76767676768	23.5395551913242\\
1686.86868686869	23.5940473546915\\
1696.9696969697	23.6482116880392\\
1707.07070707071	23.702051276855\\
1717.17171717172	23.755569165524\\
1727.27272727273	23.8087683580667\\
1737.37373737374	23.8616518188596\\
1747.47474747475	23.9142224733393\\
1757.57575757576	23.9664832086897\\
1767.67676767677	24.0184368745144\\
1777.77777777778	24.0700862834928\\
1787.87878787879	24.1214342120217\\
1797.9797979798	24.1724834008429\\
1808.08080808081	24.2232365556559\\
1818.18181818182	24.2736963477175\\
1828.28282828283	24.3238654144277\\
1838.38383838384	24.3737463599037\\
1848.48484848485	24.4233417555395\\
1858.58585858586	24.4726541405556\\
1868.68686868687	24.5216860225348\\
1878.78787878788	24.5704398779473\\
1888.88888888889	24.6189181526652\\
1898.9898989899	24.6671232624648\\
1909.09090909091	24.7150575935189\\
1919.19191919192	24.7627235028791\\
1929.29292929293	24.8101233189472\\
1939.39393939394	24.8572593419374\\
1949.49494949495	24.9041338443285\\
1959.59595959596	24.9507490713073\\
1969.69696969697	24.9971072412023\\
1979.79797979798	25.043210545909\\
1989.89898989899	25.0890611513064\\
2000	25.1346611976651\\
};

\end{axis}
\end{tikzpicture}%

%% file: figure6_unreg.tex
%
%
\begin{tikzpicture}
\pgfplotsset{every axis y label/.style={
at={(-0.46,0.5)},
xshift=32pt,
rotate=90}}

\begin{axis}[%
width=1.794in,
height=1.03in,
at={(1.358in,0.0in)},
scale only axis,
xmin=1000,
xmax=2000,
xtick={1000,1250,1500,1750,2000},
xticklabels={{1K},{1.25K},{1.5K},{1.75K},{2K}},
xlabel style={font=\color{white!15!black}},
xlabel={Potential Passenger/min},
ymin=32.6,
ymax=34,
ylabel style={font=\color{white!15!black}},
ylabel={Travel Cost (\$/trip)},
axis background/.style={fill=white},
legend style={legend cell align=left, align=left, draw=white!15!black}
]
\addplot [color=black, line width=1.0pt]
  table[row sep=crcr]{%
1000	32.6311662081556\\
1010.10101010101	32.6479141161858\\
1020.20202020202	32.6645868919424\\
1030.30303030303	32.6811841840524\\
1040.40404040404	32.6977056876774\\
1050.50505050505	32.7141511418014\\
1060.60606060606	32.730520326678\\
1070.70707070707	32.7468130614285\\
1080.80808080808	32.7630292017797\\
1090.90909090909	32.7791686379327\\
1101.0101010101	32.7952312925556\\
1111.11111111111	32.8112171188901\\
1121.21212121212	32.8271260989669\\
1131.31313131313	32.8429582419216\\
1141.41414141414	32.8587135824064\\
1151.51515151515	32.8743921790895\\
1161.61616161616	32.8899941132401\\
1171.71717171717	32.9055194873899\\
1181.81818181818	32.9209684240702\\
1191.91919191919	32.9363410646174\\
1202.0202020202	32.9516375680445\\
1212.12121212121	32.9668581099741\\
1222.22222222222	32.9820028816284\\
1232.32323232323	32.9970720888752\\
1242.42424242424	33.0120659513243\\
1252.52525252525	33.0269847014727\\
1262.62626262626	33.0418285838953\\
1272.72727272727	33.0565978544791\\
1282.82828282828	33.0712927796977\\
1292.92929292929	33.0859136359245\\
1303.0303030303	33.1004607087819\\
1313.13131313131	33.114934292525\\
1323.23232323232	33.1293346894576\\
1333.33333333333	33.1436622093787\\
1343.43434343434	33.1579171690579\\
1353.53535353535	33.1720998917388\\
1363.63636363636	33.1862107066672\\
1373.73737373737	33.2002499486448\\
1383.83838383838	33.2142179576057\\
1393.93939393939	33.2281150782146\\
1404.0404040404	33.2419416594864\\
1414.14141414141	33.2556980544252\\
1424.24242424242	33.269384619682\\
1434.34343434343	33.28300171523\\
1444.44444444444	33.2965497040568\\
1454.54545454545	33.310028951873\\
1464.64646464646	33.3234398268347\\
1474.74747474747	33.3367826992819\\
1484.84848484848	33.3500579414891\\
1494.94949494949	33.3632659274299\\
1505.05050505051	33.3764070325529\\
1515.15151515152	33.3894816335702\\
1525.25252525253	33.402490108256\\
1535.35353535354	33.4154328352568\\
1545.45454545455	33.4283101939103\\
1555.55555555556	33.4411225640749\\
1565.65656565657	33.4538703259677\\
1575.75757575758	33.4665538600109\\
1585.85858585859	33.4791735466871\\
1595.9595959596	33.491729766401\\
1606.06060606061	33.5042228993504\\
1616.16161616162	33.516653325402\\
1626.26262626263	33.5290214239758\\
1636.36363636364	33.5413275739345\\
1646.46464646465	33.5535721534798\\
1656.56565656566	33.5657555400539\\
1666.66666666667	33.5778781102468\\
1676.76767676768	33.5899402397083\\
1686.86868686869	33.601942303066\\
1696.9696969697	33.6138846738465\\
1707.07070707071	33.6257677244026\\
1717.17171717172	33.6375918258438\\
1727.27272727273	33.6493573479712\\
1737.37373737374	33.6610646592161\\
1747.47474747475	33.6727141265826\\
1757.57575757576	33.6843061155936\\
1767.67676767677	33.6958409902399\\
1777.77777777778	33.7073191129328\\
1787.87878787879	33.7187408444595\\
1797.9797979798	33.7301065439417\\
1808.08080808081	33.7414165687965\\
1818.18181818182	33.7526712747004\\
1828.28282828283	33.7638710155556\\
1838.38383838384	33.7750161434584\\
1848.48484848485	33.7861070086705\\
1858.58585858586	33.7971439595917\\
1868.68686868687	33.8081273427354\\
1878.78787878788	33.8190575027056\\
1888.88888888889	33.8299347821756\\
1898.9898989899	33.8407595218693\\
1909.09090909091	33.8515320605432\\
1919.19191919192	33.8622527349704\\
1929.29292929293	33.8729218799267\\
1939.39393939394	33.883539828177\\
1949.49494949495	33.8941069104644\\
1959.59595959596	33.9046234554994\\
1969.69696969697	33.9150897899511\\
1979.79797979798	33.9255062384396\\
1989.89898989899	33.9358731235292\\
2000	33.9461907657225\\
};

\end{axis}
\end{tikzpicture}%

%% file: figure1_wage.tex
%
%
\begin{tikzpicture}

\pgfplotsset{every axis y label/.style={
at={(-0.44,0.5)},
xshift=32pt,
rotate=90}}

\begin{axis}[%
width=1.794in,
height=1.03in,
at={(1.358in,0.0in)},
scale only axis,
xmin=20,
xmax=40,
xtick={20,25,30,35,40},
xticklabels={{20},{25},{30},{35},{40}},
xlabel style={font=\color{white!15!black}},
xlabel={Wage  Floor (\$/hour)},
ymin=5000,
ymax=8500,
ytick={5000,6000,7000, 8000},
yticklabels={{5K},{6K},{7K}, {8K}},
ylabel style={font=\color{white!15!black}},
ylabel={Number of Drivers},
axis background/.style={fill=white},
legend style={legend cell align=left, align=left, draw=white!15!black}
]
\addplot [color=black, line width=1.0pt]
  table[row sep=crcr]{%
20	5089.64440901078\\
20.1919191919192	5089.64440901078\\
20.3838383838384	5089.64440901078\\
20.5757575757576	5089.64440901078\\
20.7676767676768	5089.64440901078\\
20.959595959596	5089.64440901078\\
21.1515151515152	5089.64440901078\\
21.3434343434343	5089.64440901078\\
21.5353535353535	5089.64440901078\\
21.7272727272727	5089.64440901078\\
21.9191919191919	5089.64440901078\\
22.1111111111111	5089.64440901078\\
22.3030303030303	5089.64440901078\\
22.4949494949495	5089.64440901078\\
22.6868686868687	5108.89544411668\\
22.8787878787879	5152.11406096362\\
23.0707070707071	5195.33267781055\\
23.2626262626263	5238.55129465749\\
23.4545454545455	5281.76991150442\\
23.6464646464646	5324.98852835136\\
23.8383838383838	5368.20714519829\\
24.030303030303	5411.42576204523\\
24.2222222222222	5454.64437889217\\
24.4141414141414	5497.8629957391\\
24.6060606060606	5541.08161258604\\
24.7979797979798	5584.30022943297\\
24.989898989899	5627.51884627991\\
25.1818181818182	5670.73746312684\\
25.3737373737374	5713.95607997378\\
25.5656565656566	5757.17469682071\\
25.7575757575758	5800.39331366765\\
25.9494949494949	5843.61193051458\\
26.1414141414141	5886.83054736152\\
26.3333333333333	5930.04916420846\\
26.5252525252525	5973.26778105539\\
26.7171717171717	6016.48639790233\\
26.9090909090909	6059.70501474926\\
27.1010101010101	6102.9236315962\\
27.2929292929293	6146.14224844313\\
27.4848484848485	6189.36086529007\\
27.6767676767677	6232.579482137\\
27.8686868686869	6275.79809898394\\
28.0606060606061	6319.01671583088\\
28.2525252525253	6362.23533267781\\
28.4444444444444	6405.45394952474\\
28.6363636363636	6448.67256637168\\
28.8282828282828	6491.89118321862\\
29.020202020202	6535.10980006555\\
29.2121212121212	6578.32841691249\\
29.4040404040404	6621.54703375942\\
29.5959595959596	6664.76565060636\\
29.7878787878788	6707.98426745329\\
29.979797979798	6751.20288430023\\
30.1717171717172	6794.42150114716\\
30.3636363636364	6837.6401179941\\
30.5555555555556	6880.85873484104\\
30.7474747474747	6924.07735168797\\
30.9393939393939	6967.29596853491\\
31.1313131313131	7010.51458538184\\
31.3232323232323	7053.73320222878\\
31.5151515151515	7096.95181907571\\
31.7070707070707	7140.17043592265\\
31.8989898989899	7183.38905276958\\
32.0909090909091	7226.60766961652\\
32.2828282828283	7269.82628646345\\
32.4747474747475	7313.04490331039\\
32.6666666666667	7356.26352015732\\
32.8585858585859	7399.48213700426\\
33.0505050505051	7442.7007538512\\
33.2424242424242	7485.91937069813\\
33.4343434343434	7529.13798754507\\
33.6262626262626	7572.356604392\\
33.8181818181818	7615.57522123894\\
34.010101010101	7658.79383808587\\
34.2020202020202	7702.01245493281\\
34.3939393939394	7745.23107177974\\
34.5858585858586	7788.44968862668\\
34.7777777777778	7831.66830547361\\
34.969696969697	7874.88692232055\\
35.1616161616162	7918.10553916748\\
35.3535353535354	7961.32415601442\\
35.5454545454545	8004.54277286136\\
35.7373737373737	8047.76138970829\\
35.9292929292929	8090.98000655523\\
36.1212121212121	8134.19862340216\\
36.3131313131313	8177.4172402491\\
36.5050505050505	8120.56757385936\\
36.6969696969697	8060.40688590634\\
36.8888888888889	8000.20206926436\\
37.0808080808081	7939.95076867615\\
37.2727272727273	7879.63967498404\\
37.4646464646465	7819.26867548742\\
37.6565656565657	7758.81054099207\\
37.8484848484849	7698.25990437203\\
38.040404040404	7637.61727909\\
38.2323232323232	7576.8580920251\\
38.4242424242424	7515.98350076592\\
38.6161616161616	7454.96497570913\\
38.8080808080808	7393.7876937524\\
39	7332.45617776892\\
};

\end{axis}
\end{tikzpicture}%

%% file: figure2_wage.tex
%
%

\pgfplotsset{every axis/.append style={
font=\small,
thin,
tick style={ultra thin}}}

\begin{tikzpicture}

\pgfplotsset{every axis y label/.style={
at={(-0.43,0.5)},
xshift=32pt,
rotate=90}}

\begin{axis}[%
width=1.794in,
height=1.03in,
at={(1.358in,0.0in)},
scale only axis,
xmin=20,
xmax=40,
xtick={20,25,30,35,40},
xticklabels={{20},{25},{30},{35},{40}},
xlabel style={font=\color{white!15!black}},
xlabel={Wage Floor (\$/hour)},
ymin=180,
ymax=320,
ylabel style={font=\color{white!15!black}},
ylabel={Ride Arrival/min},
axis background/.style={fill=white},
legend style={legend cell align=left, align=left, draw=white!15!black}
]
\addplot [color=black, line width=1.0pt]
  table[row sep=crcr]{%
20	187.289529766413\\
20.1919191919192	187.289529766413\\
20.3838383838384	187.289529766413\\
20.5757575757576	187.289529766413\\
20.7676767676768	187.289529766413\\
20.959595959596	187.289529766413\\
21.1515151515152	187.289529766413\\
21.3434343434343	187.289529766413\\
21.5353535353535	187.289529766413\\
21.7272727272727	187.289529766413\\
21.9191919191919	187.289529766413\\
22.1111111111111	187.289529766413\\
22.3030303030303	187.289529766413\\
22.4949494949495	187.289529766413\\
22.6868686868687	188.056154793706\\
22.8787878787879	189.770968193523\\
23.0707070707071	191.483053492817\\
23.2626262626263	193.192745455729\\
23.4545454545455	194.900589399409\\
23.6464646464646	196.605434098762\\
23.8383838383838	198.30846107428\\
24.030303030303	200.008428214676\\
24.2222222222222	201.706607926635\\
24.4141414141414	203.402062567612\\
24.6060606060606	205.094789108066\\
24.7979797979798	206.785092016741\\
24.989898989899	208.472365385692\\
25.1818181818182	210.157851326205\\
25.3737373737374	211.840277431596\\
25.5656565656566	213.520279905208\\
25.7575757575758	215.197283134492\\
25.9494949494949	216.871862731997\\
26.1414141414141	218.543686963122\\
26.3333333333333	220.212207481177\\
26.5252525252525	221.878304367453\\
26.7171717171717	223.54197762195\\
26.9090909090909	225.202621336722\\
27.1010101010101	226.860235511769\\
27.2929292929293	228.514850442489\\
27.4848484848485	230.167011446032\\
27.6767676767677	231.815537001906\\
27.8686868686869	233.461669221397\\
28.0606060606061	235.104470461961\\
28.2525252525253	236.744513306606\\
28.4444444444444	238.381556906924\\
28.6363636363636	240.015570967517\\
28.8282828282828	241.646555488386\\
29.020202020202	243.27451046953\\
29.2121212121212	244.898860298402\\
29.4040404040404	246.520150292152\\
29.5959595959596	248.139016654123\\
29.7878787878788	249.753671955876\\
29.979797979798	251.36590362585\\
30.1717171717172	252.974499848154\\
30.3636363636364	254.580036235337\\
30.5555555555556	256.181967470247\\
30.7474747474747	257.781203929573\\
30.9393939393939	259.376470177089\\
31.1313131313131	260.968706884881\\
31.3232323232323	262.557944348346\\
31.5151515151515	264.143516068744\\
31.7070707070707	265.72548263687\\
31.8989898989899	267.304419665271\\
32.0909090909091	268.879690950606\\
32.2828282828283	270.451357083668\\
32.4747474747475	272.019418064458\\
32.6666666666667	273.584419210126\\
32.8585858585859	275.145784908125\\
33.0505050505051	276.703545453852\\
33.2424242424242	278.257640256512\\
33.4343434343434	279.8081299069\\
33.6262626262626	281.354984109618\\
33.8181818181818	282.898202864667\\
34.010101010101	284.438392079992\\
34.2020202020202	285.974339939702\\
34.3939393939394	287.506652351743\\
34.5858585858586	289.035359611511\\
34.7777777777778	290.560401128213\\
34.969696969697	292.081837492643\\
35.1616161616162	293.599638409403\\
35.3535353535354	295.113197970549\\
35.5454545454545	296.623122084025\\
35.7373737373737	298.129410749832\\
35.9292929292929	299.632063967969\\
36.1212121212121	301.13050612589\\
36.3131313131313	302.625282540744\\
36.5050505050505	300.658352359822\\
36.6969696969697	298.569481690217\\
36.8888888888889	296.471796574785\\
37.0808080808081	294.365053135576\\
37.2727272727273	292.24891660845\\
37.4646464646465	290.123390022949\\
37.6565656565657	287.987563002385\\
37.8484848484849	285.841161373411\\
38.040404040404	283.68418816557\\
38.2323232323232	281.516004145977\\
38.4242424242424	279.336612344174\\
38.6161616161616	277.14479791473\\
38.8080808080808	274.939988274638\\
39	272.722756006905\\
};

\end{axis}
\end{tikzpicture}%

%% file: figureoccupancy_wage.tex
%
%
\begin{tikzpicture}

\pgfplotsset{every axis/.append style={
font=\small,
thin,
tick style={ultra thin}}}
\pgfplotsset{every axis y label/.style={
at={(-0.45,0.5)},
xshift=32pt,
rotate=90}}

\begin{axis}[%
width=1.794in,
height=1.03in,
at={(1.358in,0in)},
scale only axis,
xmin=20,
xmax=40,
xtick={20,25,30,35,40},
xticklabels={{20},{25},{30},{35},{40}},
xlabel style={font=\color{white!15!black}},
xlabel={Wage Floor (\$/hour)},
ymin=0.5,
ymax=0.7,
ylabel style={font=\color{white!15!black}},
ylabel={Occupancy Rate},
axis background/.style={fill=white},
legend style={legend cell align=left, align=left, draw=white!15!black}
]
\addplot [color=black, line width=1.0pt]
  table[row sep=crcr]{%
20	0.599809945423257\\
20.1919191919192	0.599809945423257\\
20.3838383838384	0.599809945423257\\
20.5757575757576	0.599809945423257\\
20.7676767676768	0.599809945423257\\
20.959595959596	0.599809945423257\\
21.1515151515152	0.599809945423257\\
21.3434343434343	0.599809945423257\\
21.5353535353535	0.599809945423257\\
21.7272727272727	0.599809945423257\\
21.9191919191919	0.599809945423257\\
22.1111111111111	0.599809945423257\\
22.3030303030303	0.599809945423257\\
22.4949494949495	0.599809945423257\\
22.6868686868687	0.599995704877338\\
22.8787878787879	0.600387868931589\\
23.0707070707071	0.600764949136667\\
23.2626262626263	0.601128360457198\\
23.4545454545455	0.601480121330293\\
23.6464646464646	0.601816991482233\\
23.8383838383838	0.602142918125296\\
24.030303030303	0.602454422042568\\
24.2222222222222	0.602755648365826\\
24.4141414141414	0.603044059559429\\
24.6060606060606	0.603319946572899\\
24.7979797979798	0.603584488904733\\
24.989898989899	0.603836192931296\\
25.1818181818182	0.60407892251959\\
25.3737373737374	0.604309251559886\\
25.5656565656566	0.604529260572352\\
25.7575757575758	0.60473756267991\\
25.9494949494949	0.604936023227719\\
26.1414141414141	0.605123940436079\\
26.3333333333333	0.605300037579421\\
26.5252525252525	0.605466972811737\\
26.7171717171717	0.605624943573077\\
26.9090909090909	0.605772511839086\\
27.1010101010101	0.605909898609477\\
27.2929292929293	0.606037399013355\\
27.4848484848485	0.60615665627286\\
27.6767676767677	0.606264751851263\\
27.8686868686869	0.60636514245493\\
28.0606060606061	0.606455567514048\\
28.2525252525253	0.60653769706967\\
28.4444444444444	0.606611086146542\\
28.6363636363636	0.606675833902937\\
28.8282828282828	0.606732112923041\\
29.020202020202	0.606780091225639\\
29.2121212121212	0.606818506142251\\
29.4040404040404	0.606848887318208\\
29.5959595959596	0.606872947002153\\
29.7878787878788	0.606886464035543\\
29.979797979798	0.606893956428039\\
30.1717171717172	0.606892631967079\\
30.3636363636364	0.606884030020191\\
30.5555555555556	0.606866995920314\\
30.7474747474747	0.6068438306842\\
30.9393939393939	0.606811664522355\\
31.1313131313131	0.606772851039133\\
31.3232323232323	0.606727582427668\\
31.5151515151515	0.606674445830079\\
31.7070707070707	0.606613722438586\\
31.8989898989899	0.606546855326461\\
32.0909090909091	0.606472519730332\\
32.2828282828283	0.606390984702375\\
32.4747474747475	0.606302377884152\\
32.6666666666667	0.606208032230701\\
32.8585858585859	0.60610678030749\\
33.0505050505051	0.605998808774351\\
33.2424242424242	0.605884102082997\\
33.4343434343434	0.605762907390887\\
33.6262626262626	0.605635270574398\\
33.8181818181818	0.605501301311275\\
34.010101010101	0.605362396340807\\
34.2020202020202	0.605216074148482\\
34.3939393939394	0.605063734045129\\
34.5858585858586	0.604905539615595\\
34.7777777777778	0.604741461672955\\
34.969696969697	0.604571722501272\\
35.1616161616162	0.60439635243564\\
35.3535353535354	0.604214202644411\\
35.5454545454545	0.60402661678092\\
35.7373737373737	0.603833682424623\\
35.9292929292929	0.603635485283728\\
36.1212121212121	0.603430955783452\\
36.3131313131313	0.603221281303222\\
36.5050505050505	0.603496134831864\\
36.6969696969697	0.603776288273978\\
36.8888888888889	0.604046028129056\\
37.0808080808081	0.60430480060897\\
37.2727272727273	0.604552687331784\\
37.4646464646465	0.604789457126473\\
37.6565656565657	0.605015066695862\\
37.8484848484849	0.605229102194969\\
38.040404040404	0.605431261888227\\
38.2323232323232	0.605621856955352\\
38.4242424242424	0.605800529065297\\
38.6161616161616	0.605966657218318\\
38.8080808080808	0.606119893415845\\
39	0.606260823813771\\
};

\end{axis}
\end{tikzpicture}%

%% file: figure5_wage.tex
%
%
\definecolor{mycolor1}{rgb}{0.60000,0.20000,0.00000}%
\begin{tikzpicture}
\pgfplotsset{every axis/.append style={
font=\small,
thin,
tick style={ultra thin}}}
\pgfplotsset{every axis y label/.style={
at={(-0.41,0.5)},
xshift=32pt,
rotate=90}}

\begin{axis}[%
width=1.794in,
height=1.03in,
at={(1.358in,0.0in)},
scale only axis,
xmin=20,
xmax=40,
xtick={20,25,30,35,40},
xticklabels={{20},{25},{30},{35},{40}},
xlabel style={font=\color{white!15!black}},
xlabel={Wage Floor (\$/hour)},
ymin=10,
ymax=18,
ylabel style={font=\color{white!15!black}},
ylabel={Prices (\$/trip)},
axis background/.style={fill=white},
legend style={at={(0.61,0.131)}, anchor=south west, legend cell align=left, align=left, draw=white!15!black}
]
\addplot [color=black, line width=1.0pt]
  table[row sep=crcr]{%
20	17.0003526018022\\
20.1919191919192	17.0003526018022\\
20.3838383838384	17.0003526018022\\
20.5757575757576	17.0003526018022\\
20.7676767676768	17.0003526018022\\
20.959595959596	17.0003526018022\\
21.1515151515152	17.0003526018022\\
21.3434343434343	17.0003526018022\\
21.5353535353535	17.0003526018022\\
21.7272727272727	17.0003526018022\\
21.9191919191919	17.0003526018022\\
22.1111111111111	17.0003526018022\\
22.3030303030303	17.0003526018022\\
22.4949494949495	17.0003526018022\\
22.6868686868687	17.0081304825901\\
22.8787878787879	17.0256670218759\\
23.0707070707071	17.042774260133\\
23.2626262626263	17.0594292084721\\
23.4545454545455	17.0755912350543\\
23.6464646464646	17.0913685350741\\
23.8383838383838	17.1066648219409\\
24.030303030303	17.1215943010177\\
24.2222222222222	17.1360541681558\\
24.4141414141414	17.1501311281633\\
24.6060606060606	17.1638311523981\\
24.7979797979798	17.1771345605481\\
24.989898989899	17.1900973844501\\
25.1818181818182	17.2026228144922\\
25.3737373737374	17.2148212449759\\
25.5656565656566	17.2266459721\\
25.7575757575758	17.2381487551587\\
25.9494949494949	17.2492883800317\\
26.1414141414141	17.2600965077553\\
26.3333333333333	17.2706209582182\\
26.5252525252525	17.2807970278186\\
26.7171717171717	17.2906299174318\\
26.9090909090909	17.3001712817171\\
27.1010101010101	17.3094251927344\\
27.2929292929293	17.3183933646342\\
27.4848484848485	17.3270387411383\\
27.6767676767677	17.3354544292931\\
27.8686868686869	17.3435510782224\\
28.0606060606061	17.3514025901422\\
28.2525252525253	17.3589702534795\\
28.4444444444444	17.3662756554242\\
28.6363636363636	17.3733246288232\\
28.8282828282828	17.3801207178653\\
29.020202020202	17.3866674149522\\
29.2121212121212	17.3930090423897\\
29.4040404040404	17.3991097400392\\
29.5959595959596	17.4049281964412\\
29.7878787878788	17.4105931034081\\
29.979797979798	17.4159823575546\\
30.1717171717172	17.421183395141\\
30.3636363636364	17.4261592263632\\
30.5555555555556	17.43095004782\\
30.7474747474747	17.4354966725496\\
30.9393939393939	17.4398883102534\\
31.1313131313131	17.444063974509\\
31.3232323232323	17.4480246278019\\
31.5151515151515	17.4518172746247\\
31.7070707070707	17.4554401031914\\
31.8989898989899	17.4588577319979\\
32.0909090909091	17.4621143094566\\
32.2828282828283	17.4652080160098\\
32.4747474747475	17.4681409894139\\
32.6666666666667	17.4708805213018\\
32.8585858585859	17.4734658369636\\
33.0505050505051	17.4758970552216\\
33.2424242424242	17.4781800245941\\
33.4343434343434	17.4803129191501\\
33.6262626262626	17.4822996221419\\
33.8181818181818	17.4841420832644\\
34.010101010101	17.4858049841751\\
34.2020202020202	17.4873649279486\\
34.3939393939394	17.4887863579512\\
34.5858585858586	17.4900693123205\\
34.7777777777778	17.4912193193537\\
34.969696969697	17.4922345554168\\
35.1616161616162	17.4931186683669\\
35.3535353535354	17.4939093718947\\
35.5454545454545	17.4945721637494\\
35.7373737373737	17.4951088103474\\
35.9292929292929	17.4955210602597\\
36.1212121212121	17.4958440878425\\
36.3131313131313	17.4960475851828\\
36.5050505050505	17.4957519552724\\
36.6969696969697	17.4952405997308\\
36.8888888888889	17.4945018256979\\
37.0808080808081	17.4935392849968\\
37.2727272727273	17.4923388486505\\
37.4646464646465	17.4908934849124\\
37.6565656565657	17.489192391776\\
37.8484848484849	17.4872325946529\\
38.040404040404	17.4850069036465\\
38.2323232323232	17.4824911051881\\
38.4242424242424	17.4796779551774\\
38.6161616161616	17.4765679639621\\
38.8080808080808	17.4731517758212\\
39	17.4693918363891\\
};
\addlegendentry{$p_f$}

\addplot [color=mycolor1, dotted, line width=1.0pt]
  table[row sep=crcr]{%
20	10.2366457722024\\
20.1919191919192	10.2366457722024\\
20.3838383838384	10.2366457722024\\
20.5757575757576	10.2366457722024\\
20.7676767676768	10.2366457722024\\
20.959595959596	10.2366457722024\\
21.1515151515152	10.2366457722024\\
21.3434343434343	10.2366457722024\\
21.5353535353535	10.2366457722024\\
21.7272727272727	10.2366457722024\\
21.9191919191919	10.2366457722024\\
22.1111111111111	10.2366457722024\\
22.3030303030303	10.2366457722024\\
22.4949494949495	10.2366457722024\\
22.6868686868687	10.2721835225904\\
22.8787878787879	10.3523144987332\\
23.0707070707071	10.432602795068\\
23.2626262626263	10.5130294133413\\
23.4545454545455	10.5935640362073\\
23.6464646464646	10.674268622982\\
23.8383838383838	10.7550783728502\\
24.030303030303	10.8360600974575\\
24.2222222222222	10.9171442660237\\
24.4141414141414	10.9983811503823\\
24.6060606060606	11.0797703649912\\
24.7979797979798	11.161294956418\\
24.989898989899	11.2429871518069\\
25.1818181818182	11.3247795131252\\
25.3737373737374	11.4067402334243\\
25.5656565656566	11.4888346244168\\
25.7575757575758	11.5710932796456\\
25.9494949494949	11.6534848710094\\
26.1414141414141	11.7360268982797\\
26.3333333333333	11.8187484631541\\
26.5252525252525	11.9016020024391\\
26.7171717171717	11.9845872601335\\
26.9090909090909	12.0677364644847\\
27.1010101010101	12.1510493463424\\
27.2929292929293	12.2345240353982\\
27.4848484848485	12.3181311208039\\
27.6767676767677	12.4019336369069\\
27.8686868686869	12.4858649284181\\
28.0606060606061	12.5699749849435\\
28.2525252525253	12.6542330301158\\
28.4444444444444	12.7386518048908\\
28.6363636363636	12.8232328037482\\
28.8282828282828	12.9079759170039\\
29.020202020202	12.9928810500308\\
29.2121212121212	13.0779788612763\\
29.4040404040404	13.1632401575271\\
29.5959595959596	13.2486309201177\\
29.7878787878788	13.3342465468134\\
29.979797979798	13.4199915788377\\
30.1717171717172	13.5059307032326\\
30.3636363636364	13.592033190449\\
30.5555555555556	13.6783281703062\\
30.7474747474747	13.7647670631275\\
30.9393939393939	13.8514180125692\\
31.1313131313131	13.9382308434802\\
31.3232323232323	14.0252039975331\\
31.5151515151515	14.1123731524603\\
31.7070707070707	14.1997351694039\\
31.8989898989899	14.2872593927274\\
32.0909090909091	14.3749799362827\\
32.2828282828283	14.4628936963383\\
32.4747474747475	14.5510008192515\\
32.6666666666667	14.6392722837878\\
32.8585858585859	14.7277390413845\\
33.0505050505051	14.8163996508182\\
33.2424242424242	14.9052575481201\\
33.4343434343434	14.9943096914245\\
33.6262626262626	15.0835579167487\\
33.8181818181818	15.1730024483461\\
34.010101010101	15.2626110083144\\
34.2020202020202	15.3524488499776\\
34.3939393939394	15.4424837301354\\
34.5858585858586	15.5327142842761\\
34.7777777777778	15.6231440404733\\
34.969696969697	15.7137700235079\\
35.1616161616162	15.8045941520738\\
35.3535353535354	15.8956493613399\\
35.5454545454545	15.9869033635939\\
35.7373737373737	16.0783564766862\\
35.9292929292929	16.1700090265624\\
36.1212121212121	16.2618924317346\\
36.3131313131313	16.3539776294936\\
36.5050505050505	16.4329227890885\\
36.6969696969697	16.5116511329766\\
36.8888888888889	16.5905924628319\\
37.0808080808081	16.6697658507232\\
37.2727272727273	16.7491730463527\\
37.4646464646465	16.8288244825789\\
37.6565656565657	16.9087254734042\\
37.8484848484849	16.9888917764887\\
38.040404040404	17.0693362151122\\
38.2323232323232	17.1500544310995\\
38.4242424242424	17.2310609809004\\
38.6161616161616	17.3123781329511\\
38.8080808080808	17.3940206671764\\
39	17.4759766487147\\
};
\addlegendentry{$p_d$}

\end{axis}
\end{tikzpicture}%

%% file: commission_wage.tex
%
%
\begin{tikzpicture}

\pgfplotsset{every axis/.append style={
font=\small,
thin,
tick style={ultra thin}}}
\pgfplotsset{every axis y label/.style={
at={(-0.45,0.5)},
xshift=32pt,
rotate=90}}

\begin{axis}[%
width=1.794in,
height=1.03in,
at={(1.358in,0in)},
scale only axis,
xmin=20,
xmax=40,
xtick={20,25,30,35,40},
xticklabels={{20},{25},{30},{35},{40}},
xlabel style={font=\color{white!15!black}},
xlabel={Wage Floor (\$/hour)},
ymin=-0.000376934262356323,
ymax=0.4,
ylabel style={font=\color{white!15!black}},
ylabel={Commission Rate},
axis background/.style={fill=white},
legend style={legend cell align=left, align=left, draw=white!15!black}
]
\addplot [color=black, line width=1.0pt]
  table[row sep=crcr]{%
20	0.39785685556209\\
20.1919191919192	0.39785685556209\\
20.3838383838384	0.39785685556209\\
20.5757575757576	0.39785685556209\\
20.7676767676768	0.39785685556209\\
20.959595959596	0.39785685556209\\
21.1515151515152	0.39785685556209\\
21.3434343434343	0.39785685556209\\
21.5353535353535	0.39785685556209\\
21.7272727272727	0.39785685556209\\
21.9191919191919	0.39785685556209\\
22.1111111111111	0.39785685556209\\
22.3030303030303	0.39785685556209\\
22.4949494949495	0.39785685556209\\
22.6868686868687	0.396042761248498\\
22.8787878787879	0.391958359961357\\
23.0707070707071	0.387857714018293\\
23.2626262626263	0.383740845905896\\
23.4545454545455	0.379607775193175\\
23.6464646464646	0.375458518662403\\
23.8383838383838	0.371293090453505\\
24.030303030303	0.367111502179824\\
24.2222222222222	0.362913763058056\\
24.4141414141414	0.358699880007267\\
24.6060606060606	0.354469857771637\\
24.7979797979798	0.350223699006651\\
24.989898989899	0.345961404385228\\
25.1818181818182	0.3416829726927\\
25.3737373737374	0.337388400895926\\
25.5656565656566	0.333077684244287\\
25.7575757575758	0.328750816343731\\
25.9494949494949	0.324407789222205\\
26.1414141414141	0.320048593412762\\
26.3333333333333	0.315673218018826\\
26.5252525252525	0.311281650766456\\
26.7171717171717	0.3068738780852\\
26.9090909090909	0.3024498851501\\
27.1010101010101	0.29800965595075\\
27.2929292929293	0.293553173334063\\
27.4848484848485	0.289080419058686\\
27.6767676767677	0.284591373852284\\
27.8686868686869	0.280086017442175\\
28.0606060606061	0.27556432861025\\
28.2525252525253	0.271026285238362\\
28.4444444444444	0.266471864339435\\
28.6363636363636	0.261901042102565\\
28.8282828282828	0.2573137939292\\
29.020202020202	0.252710094468297\\
29.2121212121212	0.248089917655821\\
29.4040404040404	0.24345323673455\\
29.5959595959596	0.238800024304225\\
29.7878787878788	0.234130252334526\\
29.979797979798	0.22944389220649\\
30.1717171717172	0.224740914730302\\
30.3636363636364	0.220021290182737\\
30.5555555555556	0.215284988323579\\
30.7474747474747	0.210531978432326\\
30.9393939393939	0.205762229312817\\
31.1313131313131	0.200975709338823\\
31.3232323232323	0.196172386461149\\
31.5151515151515	0.191352228230121\\
31.7070707070707	0.186515201824802\\
31.8989898989899	0.181661274062491\\
32.0909090909091	0.176790411428133\\
32.2828282828283	0.171902580084896\\
32.4747474747475	0.166997745892383\\
32.6666666666667	0.16207587442782\\
32.8585858585859	0.157136931001444\\
33.0505050505051	0.152180880672376\\
33.2424242424242	0.147207688263512\\
33.4343434343434	0.142217318375471\\
33.6262626262626	0.137209735403184\\
33.8181818181818	0.132184903549283\\
34.010101010101	0.127142786841826\\
34.2020202020202	0.122083349136204\\
34.3939393939394	0.117006554138929\\
34.5858585858586	0.11191236541673\\
34.7777777777778	0.106800746407275\\
34.969696969697	0.101671660431637\\
35.1616161616162	0.0965250707037444\\
35.3535353535354	0.0913609403466217\\
35.5454545454545	0.086179232395266\\
35.7373737373737	0.0809799098147524\\
35.9292929292929	0.0757629355039984\\
36.1212121212121	0.0705282723092684\\
36.3131313131313	0.0652758830318016\\
36.5050505050505	0.060747841470375\\
36.6969696969697	0.0562204024087161\\
36.8888888888889	0.0516681967781532\\
37.0808080808081	0.047090152590221\\
37.2727272727273	0.0424852164554963\\
37.4646464646465	0.0378522116611427\\
37.6565656565657	0.0331900356156376\\
37.8484848484849	0.0284974089220164\\
38.040404040404	0.0237729782335787\\
38.2323232323232	0.0190154064479984\\
38.4242424242424	0.0142232010746713\\
38.6161616161616	0.00939485552023511\\
38.8080808080808	0.00452872553618929\\
39	-0.000376934262356323\\
};

\end{axis}
\end{tikzpicture}%

%% file: figure3_wage.tex
%
%

\pgfplotsset{every axis/.append style={
font=\small,
thin,
tick style={ultra thin}}}
\pgfplotsset{every axis y label/.style={
at={(-0.4,0.5)},
xshift=32pt,
rotate=90}}

\begin{tikzpicture}

\begin{axis}[%
width=1.794in,
height=1.03in,
at={(1.358in,0in)},
scale only axis,
xmin=20,
xmax=40,
xtick={20,25,30,35,40},
xticklabels={{20},{25},{30},{35},{40}},
xlabel style={font=\color{white!15!black}},
xlabel={Wage Floor (\$/hour)},
ymin=22,
ymax=40,
ylabel style={font=\color{white!15!black}},
ylabel={Driver Wage},
axis background/.style={fill=white},
legend style={legend cell align=left, align=left, draw=white!15!black}
]
\addplot [color=black, line width=1.0pt]
  table[row sep=crcr]{%
20	22.6013813813814\\
20.1919191919192	22.6013813813814\\
20.3838383838384	22.6013813813814\\
20.5757575757576	22.6013813813814\\
20.7676767676768	22.6013813813814\\
20.959595959596	22.6013813813814\\
21.1515151515152	22.6013813813814\\
21.3434343434343	22.6013813813814\\
21.5353535353535	22.6013813813814\\
21.7272727272727	22.6013813813814\\
21.9191919191919	22.6013813813814\\
22.1111111111111	22.6013813813814\\
22.3030303030303	22.6013813813814\\
22.4949494949495	22.6013813813814\\
22.6868686868687	22.6868686868687\\
22.8787878787879	22.8787878787879\\
23.0707070707071	23.0707070707071\\
23.2626262626263	23.2626262626263\\
23.4545454545455	23.4545454545455\\
23.6464646464646	23.6464646464646\\
23.8383838383838	23.8383838383838\\
24.030303030303	24.030303030303\\
24.2222222222222	24.2222222222222\\
24.4141414141414	24.4141414141414\\
24.6060606060606	24.6060606060606\\
24.7979797979798	24.7979797979798\\
24.989898989899	24.989898989899\\
25.1818181818182	25.1818181818182\\
25.3737373737374	25.3737373737374\\
25.5656565656566	25.5656565656566\\
25.7575757575758	25.7575757575758\\
25.9494949494949	25.949494949495\\
26.1414141414141	26.1414141414141\\
26.3333333333333	26.3333333333333\\
26.5252525252525	26.5252525252525\\
26.7171717171717	26.7171717171717\\
26.9090909090909	26.9090909090909\\
27.1010101010101	27.1010101010101\\
27.2929292929293	27.2929292929293\\
27.4848484848485	27.4848484848485\\
27.6767676767677	27.6767676767677\\
27.8686868686869	27.8686868686869\\
28.0606060606061	28.0606060606061\\
28.2525252525253	28.2525252525253\\
28.4444444444444	28.4444444444444\\
28.6363636363636	28.6363636363636\\
28.8282828282828	28.8282828282828\\
29.020202020202	29.020202020202\\
29.2121212121212	29.2121212121212\\
29.4040404040404	29.4040404040404\\
29.5959595959596	29.5959595959596\\
29.7878787878788	29.7878787878788\\
29.979797979798	29.979797979798\\
30.1717171717172	30.1717171717172\\
30.3636363636364	30.3636363636364\\
30.5555555555556	30.5555555555556\\
30.7474747474747	30.7474747474748\\
30.9393939393939	30.9393939393939\\
31.1313131313131	31.1313131313131\\
31.3232323232323	31.3232323232323\\
31.5151515151515	31.5151515151515\\
31.7070707070707	31.7070707070707\\
31.8989898989899	31.8989898989899\\
32.0909090909091	32.0909090909091\\
32.2828282828283	32.2828282828283\\
32.4747474747475	32.4747474747475\\
32.6666666666667	32.6666666666667\\
32.8585858585859	32.8585858585859\\
33.0505050505051	33.0505050505051\\
33.2424242424242	33.2424242424242\\
33.4343434343434	33.4343434343434\\
33.6262626262626	33.6262626262626\\
33.8181818181818	33.8181818181818\\
34.010101010101	34.010101010101\\
34.2020202020202	34.2020202020202\\
34.3939393939394	34.3939393939394\\
34.5858585858586	34.5858585858586\\
34.7777777777778	34.7777777777778\\
34.969696969697	34.969696969697\\
35.1616161616162	35.1616161616162\\
35.3535353535354	35.3535353535354\\
35.5454545454545	35.5454545454545\\
35.7373737373737	35.7373737373737\\
35.9292929292929	35.9292929292929\\
36.1212121212121	36.1212121212121\\
36.3131313131313	36.3131313131313\\
36.5050505050505	36.5050505050505\\
36.6969696969697	36.6969696969697\\
36.8888888888889	36.8888888888889\\
37.0808080808081	37.0808080808081\\
37.2727272727273	37.2727272727273\\
37.4646464646465	37.4646464646465\\
37.6565656565657	37.6565656565657\\
37.8484848484849	37.8484848484849\\
38.040404040404	38.040404040404\\
38.2323232323232	38.2323232323232\\
38.4242424242424	38.4242424242424\\
38.6161616161616	38.6161616161616\\
38.8080808080808	38.8080808080808\\
39	39\\
};

\end{axis}
\end{tikzpicture}%

%% file: figure4_wage.tex
%
%
\begin{tikzpicture}

\pgfplotsset{every axis/.append style={
font=\small,
thin,
tick style={ultra thin}}}
\pgfplotsset{every axis y label/.style={
at={(-0.41,0.5)},
xshift=32pt,
rotate=90}}

\begin{axis}[%
width=1.794in,
height=1.03in,
at={(1.358in,0.0in)},
scale only axis,
xmin=20,
xmax=40,
xtick={20,25,30,35,40},
xticklabels={{20},{25},{30},{35},{40}},
xlabel style={font=\color{white!15!black}},
xlabel={Wage Floor (\$/hour)},
ymin=3.8,
ymax=5,
ylabel style={font=\color{white!15!black}},
ylabel={Pickup Time (min)},
axis background/.style={fill=white},
legend style={legend cell align=left, align=left, draw=white!15!black}
]
\addplot [color=black, line width=1.0pt]
  table[row sep=crcr]{%
20	4.99981625054523\\
20.1919191919192	4.99981625054523\\
20.3838383838384	4.99981625054523\\
20.5757575757576	4.99981625054523\\
20.7676767676768	4.99981625054523\\
20.959595959596	4.99981625054523\\
21.1515151515152	4.99981625054523\\
21.3434343434343	4.99981625054523\\
21.5353535353535	4.99981625054523\\
21.7272727272727	4.99981625054523\\
21.9191919191919	4.99981625054523\\
22.1111111111111	4.99981625054523\\
22.3030303030303	4.99981625054523\\
22.4949494949495	4.99981625054523\\
22.6868686868687	4.9915459727123\\
22.8787878787879	4.97300440193703\\
23.0707070707071	4.95461481914931\\
23.2626262626263	4.93638166456434\\
23.4545454545455	4.91831314268935\\
23.6464646464646	4.90038509226182\\
23.8383838383838	4.88261780753582\\
24.030303030303	4.86498599229998\\
24.2222222222222	4.84751123059\\
24.4141414141414	4.83017427619997\\
24.6060606060606	4.81297333023868\\
24.7979797979798	4.79591205631946\\
24.989898989899	4.7789780165963\\
25.1818181818182	4.7621911988671\\
25.3737373737374	4.74552770291737\\
25.5656565656566	4.72899688705993\\
25.7575757575758	4.71258738421015\\
25.9494949494949	4.6963073435322\\
26.1414141414141	4.6801496551525\\
26.3333333333333	4.66410394439205\\
26.5252525252525	4.64818318641835\\
26.7171717171717	4.63238579415439\\
26.9090909090909	4.61670066547504\\
27.1010101010101	4.60112655776781\\
27.2929292929293	4.58566271492028\\
27.4848484848485	4.57031625033679\\
27.6767676767677	4.55506784001948\\
27.8686868686869	4.53993519569022\\
28.0606060606061	4.52490297792758\\
28.2525252525253	4.50997859553207\\
28.4444444444444	4.49515731363358\\
28.6363636363636	4.48043758533487\\
28.8282828282828	4.46581832910635\\
29.020202020202	4.45129847922159\\
29.2121212121212	4.4368689380357\\
29.4040404040404	4.42253645394412\\
29.5959595959596	4.40830873975566\\
29.7878787878788	4.39416027316754\\
29.979797979798	4.3801145659861\\
30.1717171717172	4.36615425771603\\
30.3636363636364	4.35228620844711\\
30.5555555555556	4.33850234648818\\
30.7474747474747	4.32481373054279\\
30.9393939393939	4.31120294830184\\
31.1313131313131	4.29768137455521\\
31.3232323232323	4.28424848256528\\
31.5151515151515	4.27089505529111\\
31.7070707070707	4.25762117955052\\
31.8989898989899	4.24443322465137\\
32.0909090909091	4.23132261207246\\
32.2828282828283	4.21828943100314\\
32.4747474747475	4.20533302909478\\
32.6666666666667	4.19245919296169\\
32.8585858585859	4.17966038034259\\
33.0505050505051	4.1669363219486\\
33.2424242424242	4.15428569912989\\
33.4343434343434	4.14170860302071\\
33.6262626262626	4.12920408165379\\
33.8181818181818	4.1167715400686\\
34.010101010101	4.10441709526924\\
34.2020202020202	4.09212669256934\\
34.3939393939394	4.07990652168897\\
34.5858585858586	4.06775633796066\\
34.7777777777778	4.0556749204519\\
34.969696969697	4.04366235959593\\
35.1616161616162	4.03171777525775\\
35.3535353535354	4.01983431947891\\
35.5454545454545	4.00801786022855\\
35.7373737373737	3.99626785847455\\
35.9292929292929	3.98458378063037\\
36.1212121212121	3.97295932060056\\
36.3131313131313	3.96139953248702\\
36.5050505050505	3.97661920094652\\
36.6969696969697	3.99284265491926\\
36.8888888888889	4.00920330689983\\
37.0808080808081	4.02570191840131\\
37.2727272727273	4.04234537593879\\
37.4646464646465	4.05913580198721\\
37.6565656565657	4.07608349537896\\
37.8484848484849	4.09319147274425\\
38.040404040404	4.11046190481837\\
38.2323232323232	4.12790709743456\\
38.4242424242424	4.14552923786869\\
38.6161616161616	4.16333751487023\\
38.8080808080808	4.18133918764976\\
39	4.19954131187932\\
};

\end{axis}
\end{tikzpicture}%

%% file: GC_wage.tex
%
%
\begin{tikzpicture}

\pgfplotsset{every axis/.append style={
font=\small,
thin,
tick style={ultra thin}}}
\pgfplotsset{every axis y label/.style={
at={(-0.45,0.5)},
xshift=32pt,
rotate=90}}

\begin{axis}[%
width=1.794in,
height=1.03in,
at={(1.358in,0in)},
scale only axis,
xmin=20,
xmax=40,
xtick={20,25,30,35,40},
xticklabels={{20},{25},{30},{35},{40}},
xlabel style={font=\color{white!15!black}},
xlabel={Wage Floor (\$/hour)},
ymin=30,
ymax=34,
ylabel style={font=\color{white!15!black}},
ylabel={Total Cost (\$/trip)},
axis background/.style={fill=white},
legend style={legend cell align=left, align=left, draw=white!15!black}
]
\addplot [color=black, line width=1.0pt]
  table[row sep=crcr]{%
20	33.3848260602649\\
20.1919191919192	33.3848260602649\\
20.3838383838384	33.3848260602649\\
20.5757575757576	33.3848260602649\\
20.7676767676768	33.3848260602649\\
20.959595959596	33.3848260602649\\
21.1515151515152	33.3848260602649\\
21.3434343434343	33.3848260602649\\
21.5353535353535	33.3848260602649\\
21.7272727272727	33.3848260602649\\
21.9191919191919	33.3848260602649\\
22.1111111111111	33.3848260602649\\
22.3030303030303	33.3848260602649\\
22.4949494949495	33.3848260602649\\
22.6868686868687	33.3655021155341\\
22.8787878787879	33.3222776470103\\
23.0707070707071	33.2791219443914\\
23.2626262626263	33.2360265694404\\
23.4545454545455	33.1929777766125\\
23.6464646464646	33.1500045842798\\
23.8383838383838	33.1070772104289\\
24.030303030303	33.0642269643559\\
24.2222222222222	33.0214217731234\\
24.4141414141414	32.9786852714316\\
24.6060606060606	32.9360175356447\\
24.7979797979798	32.8934108911669\\
24.989898989899	32.8508806108257\\
25.1818181818182	32.8083953853249\\
25.3737373737374	32.765987287602\\
25.5656565656566	32.7236402811882\\
25.7575757575758	32.6813688752695\\
25.9494949494949	32.63915856066\\
26.1414141414141	32.5970176992326\\
26.3333333333333	32.554960112896\\
26.5252525252525	32.5129636178685\\
26.7171717171717	32.4710282141503\\
26.9090909090909	32.4291691745685\\
27.1010101010101	32.3873864991231\\
27.2929292929293	32.3456794241729\\
27.4848484848485	32.3040342041731\\
27.6767676767677	32.2624806211372\\
27.8686868686869	32.220987365769\\
28.0606060606061	32.1795780727689\\
28.2525252525253	32.138238309315\\
28.4444444444444	32.0969741463562\\
28.6363636363636	32.0557863475338\\
28.8282828282828	32.014674912848\\
29.020202020202	31.9736398422985\\
29.2121212121212	31.9326956450715\\
29.4040404040404	31.8918285756223\\
29.5959595959596	31.8510225974822\\
29.7878787878788	31.8103227654919\\
29.979797979798	31.7696840248107\\
30.1717171717172	31.7291369210932\\
30.3636363636364	31.6886669451536\\
30.5555555555556	31.6482878425364\\
30.7474747474747	31.6079766658185\\
30.9393939393939	31.5677655643015\\
31.1313131313131	31.5276308269209\\
31.3232323232323	31.4875716900354\\
31.5151515151515	31.447604953755\\
31.7070707070707	31.4077290907971\\
31.8989898989899	31.3679295919755\\
32.0909090909091	31.3282224937591\\
32.2828282828283	31.2886062688651\\
32.4747474747475	31.2490809172935\\
32.6666666666667	31.2096326934997\\
32.8585858585859	31.1702761066696\\
33.0505050505051	31.1310103931619\\
33.2424242424242	31.0918370802594\\
33.4343434343434	31.0527546406792\\
33.6262626262626	31.0137638380628\\
33.8181818181818	30.9748646724101\\
34.010101010101	30.9360418708939\\
34.2020202020202	30.8973259791687\\
34.3939393939394	30.8587017244073\\
34.5858585858586	30.8201683429682\\
34.7777777777778	30.7817273621343\\
34.969696969697	30.7433772546228\\
35.1616161616162	30.705118784075\\
35.3535353535354	30.6669672233182\\
35.5454545454545	30.6289072995252\\
35.7373737373737	30.590939012696\\
35.9292929292929	30.5530623628305\\
36.1212121212121	30.5152918591147\\
36.3131313131313	30.477613756004\\
36.5050505050505	30.5271932097817\\
36.6969696969697	30.579846358234\\
36.8888888888889	30.6327216881413\\
37.0808080808081	30.6858253468167\\
37.2727272727273	30.7391657724971\\
37.4646464646465	30.7927428888185\\
37.6565656565657	30.8465796432038\\
37.8484848484849	30.9006829466074\\
38.040404040404	30.9550527226652\\
38.2323232323232	31.00970508421\\
38.4242424242424	31.0646399548776\\
38.6161616161616	31.1198879566867\\
38.8080808080808	31.1754635224592\\
39	31.2313522193733\\
};

\end{axis}
\end{tikzpicture}%

%% file: figure6_wage.tex
%
%
\begin{tikzpicture}

\pgfplotsset{every axis/.append style={
font=\small,
thin,
tick style={ultra thin}}}
\pgfplotsset{every axis y label/.style={
at={(-0.4,0.5)},
xshift=32pt,
rotate=90}}

\begin{axis}[%
width=1.794in,
height=1.03in,
at={(1.358in,0.0in)},
scale only axis,
xmin=20,
xmax=40,
xtick={20,25,30,35,40},
xticklabels={{20},{25},{30},{35},{40}},
xlabel style={font=\color{white!15!black}},
xlabel={Wage Floor (\$/hour)},
ymin=-10000,
ymax=80000,
ylabel style={font=\color{white!15!black}},
ylabel={Platform Rent},
axis background/.style={fill=white},
legend style={legend cell align=left, align=left, draw=white!15!black}
]
\addplot [color=black, line width=1.0pt]
  table[row sep=crcr]{%
20	76006.2882956179\\
20.1919191919192	76006.2882956179\\
20.3838383838384	76006.2882956179\\
20.5757575757576	76006.2882956179\\
20.7676767676768	76006.2882956179\\
20.959595959596	76006.2882956179\\
21.1515151515152	76006.2882956179\\
21.3434343434343	76006.2882956179\\
21.5353535353535	76006.2882956179\\
21.7272727272727	76006.2882956179\\
21.9191919191919	76006.2882956179\\
22.1111111111111	76006.2882956179\\
22.3030303030303	76006.2882956179\\
22.4949494949495	76006.2882956179\\
22.6868686868687	76004.1770515142\\
22.8787878787879	75984.5141648079\\
23.0707070707071	75944.1489745042\\
23.2626262626263	75883.0169563291\\
23.4545454545455	75801.0552934974\\
23.6464646464646	75698.20283232\\
23.8383838383838	75574.4000459654\\
24.030303030303	75429.5889926188\\
24.2222222222222	75263.7132816259\\
24.4141414141414	75076.718032856\\
24.6060606060606	74868.549846833\\
24.7979797979798	74639.1567658806\\
24.989898989899	74388.4882445146\\
25.1818181818182	74116.4951188443\\
25.3737373737374	73823.129569937\\
25.5656565656566	73508.3451000363\\
25.7575757575758	73172.0965021716\\
25.9494949494949	72814.3398290353\\
26.1414141414141	72435.0323697441\\
26.3333333333333	72034.1326226391\\
26.5252525252525	71611.6002655818\\
26.7171717171717	71167.396137533\\
26.9090909090909	70701.482208856\\
27.1010101010101	70213.8215622429\\
27.2929292929293	69704.3783663705\\
27.4848484848485	69173.1178548966\\
27.6767676767677	68620.0063078242\\
27.8686868686869	68045.0110252177\\
28.0606060606061	67448.1003103464\\
28.2525252525253	66829.243450844\\
28.4444444444444	66188.4106961864\\
28.6363636363636	65525.5732409998\\
28.8282828282828	64840.7032064072\\
29.020202020202	64133.7736220374\\
29.2121212121212	63404.7584104729\\
29.4040404040404	62653.6323659317\\
29.5959595959596	61880.3711441466\\
29.7878787878788	61084.9512404529\\
29.979797979798	60267.3499783642\\
30.1717171717172	59427.5454908791\\
30.3636363636364	58565.516709234\\
30.5555555555556	57681.2433452031\\
30.7474747474747	56774.7058810056\\
30.9393939393939	55845.8855489207\\
31.1313131313131	54894.7643255248\\
31.3232323232323	53921.324913389\\
31.5151515151515	52925.5507286155\\
31.7070707070707	51907.4258912351\\
31.8989898989899	50866.9352094025\\
32.0909090909091	49804.0641709197\\
32.2828282828283	48718.798928813\\
32.4747474747475	47611.1262900028\\
32.6666666666667	46481.0337054518\\
32.8585858585859	45328.509258772\\
33.0505050505051	44153.5416554243\\
33.2424242424242	42956.1202119956\\
33.4343434343434	41736.2348455454\\
33.6262626262626	40493.8760642402\\
33.8181818181818	39229.0349572607\\
34.010101010101	37941.7031865448\\
34.2020202020202	36631.8729734642\\
34.3939393939394	35299.5370933131\\
34.5858585858586	33944.6888650918\\
34.7777777777778	32567.32214215\\
34.969696969697	31167.4313037262\\
35.1616161616162	29745.0112457441\\
35.3535353535354	28300.0573741977\\
35.5454545454545	26832.5655942483\\
35.7373737373737	25342.532304772\\
35.9292929292929	23829.9543882662\\
36.1212121212121	22294.8292043183\\
36.3131313131313	20737.1545813048\\
36.5050505050505	19172.9079566887\\
36.6969696969697	17620.187837085\\
36.8888888888889	16079.0179649776\\
37.0808080808081	14549.4066450961\\
37.2727272727273	13031.364034919\\
37.4646464646465	11524.9022031647\\
37.6565656565657	10030.0351935236\\
37.8484848484849	8546.7790934305\\
38.040404040404	7075.15210926564\\
38.2323232323232	5615.17464726123\\
38.4242424242424	4166.86940194755\\
38.6161616161616	2730.26145211173\\
38.8080808080808	1305.37836497883\\
39	-107.749689914581\\
};

\end{axis}
\end{tikzpicture}%

%% file: Monopsony-Fig1.tex
\begin{tikzpicture}[dot/.style={circle,inner sep=1pt,fill,label={#1},name=#1},
extended line/.style={shorten >=-#1,shorten <=-#1},
extended line/.default=1cm,
one end extended/.style={shorten >=-#1},
one end extended/.default=1cm,
]
	\node(origin) at (0,0) {};
	\node [left of=origin,node distance=3pt] (originol){};
	\node [below of=origin,node distance=1.5pt] (origino){};
	\node [left of=origin,node distance=8pt] (originleft){};
	\node [below of=origin,node distance=8pt] (originbelow){};
	\node [right of=origin,node distance=1.6in] (xaxis) {};
	\node [above of=origin,node distance=1.2in] (yaxis) {Wage};
	\node [below of=xaxis, node distance=10pt] (xlabel) {\textbf{\#} Drivers};	
	\draw [->,line width=0.3mm, black] (originleft) -- (xaxis);
	\draw [->,line width=0.3mm, black] (originbelow) -- (yaxis);

	\node [right of=yaxis,node distance=15pt] (MRP0) {};
	\node [above of=xaxis, node distance=10pt] (MRP1) {};
	\draw [name path=MRP,line width=0.3mm, red] (MRP0.center) node[above right]{\footnotesize MRP} -- (MRP1.center) ;

	\node [below of=MRP0,node distance=1.1in] (WL0) {};
	\node [above of=MRP1, node distance=1 in] (WL1) {};
	\draw [name path=WL, line width=0.3mm,black!60!green] (WL0.center) -- (WL1.center) node[above left]{\footnotesize $W(L)$};
	
	\path [name intersections={of=MRP and WL,by=E}];
	\node [fill=black,inner sep=2pt,label=-0:$E$] at (E) {};
	\node [below of=E, node distance=0.645 in] (Lc) {};
	\draw [line width=0.3mm,dotted] (E.center) -- (Lc.center) node[below]{\footnotesize $L^*_c$};
	\node [left of=E, node distance=1 in] (Wc) {};
	\draw [line width=0.3mm,dotted] (E.center) -- (Wc.east) node[left]{\footnotesize $\omega^*_c$};

	\node [above of=E, node distance=0.8 in] (MCL1) {};
	\draw [name path=MCL, line width=0.3mm,black!60!blue,dotted,one end extended] (WL0.center) -- (MCL1.east);
	\node [left of=MCL1, node distance=10pt,black!60!blue] (MCLlabel) {\footnotesize MCL};
	
	\path [name intersections={of=MRP and MCL,by=A}];
	\node [fill=black,inner sep=2pt,label=-0:$A$] at (A) {};
	\node [below of=A, node distance=0.86 in] (L1) {};
	\draw [name path=L1path,line width=0.3mm,dotted] (A.center) -- (L1.center) node[below]{\footnotesize $L^*_1$};
	\path [name intersections={of=WL and L1path,by=AA}];
	\node [fill=black,inner sep=2pt,label=-0:$A'$] at (AA) {};
	\node [left of=AA, node distance=0.66 in] (W1) {};
	\draw [line width=0.3mm,dotted] (AA.center) -- (W1.center) node[left]{\footnotesize $\omega^*_1$};

	\node [left of=AA, node distance=0.66 in] (W1) {};
	\node [below of=W1, node distance=0.15 in] (Wm) {};
	\draw [name path=MLCm1,opacity=0] (Wm.center) --+(1.5in,0);
	\path [name intersections={of=MCL and MLCm1,by=M1}];
	\draw [line width=0.3mm,black!60!blue] (Wm.center) node[left]{\footnotesize $\omega_m$} -- (M1.west) ;
	\draw [line width=0.3mm,black!60!blue,one end extended] (M1.center) -- (MCL1.east) ;

\end{tikzpicture}						

%% file: Monopsony-Fig2.tex
\begin{tikzpicture}[dot/.style={circle,inner sep=1pt,fill,label={#1},name=#1},
extended line/.style={shorten >=-#1,shorten <=-#1},
extended line/.default=1cm,
one end extended/.style={shorten >=-#1},
one end extended/.default=1cm,
]
	\node(origin) at (0,0) {};
	\node [left of=origin,node distance=3pt] (originol){};
	\node [below of=origin,node distance=1.5pt] (origino){};
	\node [left of=origin,node distance=8pt] (originleft){};
	\node [below of=origin,node distance=8pt] (originbelow){};
	\node [right of=origin,node distance=1.6in] (xaxis) {};
	\node [above of=origin,node distance=1.2in] (yaxis) {Wage};
	\node [below of=xaxis, node distance=10pt] (xlabel) {\textbf{\#} Drivers};	
	\draw [->,line width=0.3mm, black] (originleft) -- (xaxis);
	\draw [->,line width=0.3mm, black] (originbelow) -- (yaxis);

	\node [right of=yaxis,node distance=15pt] (MRP0) {};
	\node [above of=xaxis, node distance=10pt] (MRP1) {};
	\draw [name path=MRP,line width=0.3mm, red] (MRP0.center) node[above right]{\footnotesize MRP} -- (MRP1.center) ;

	\node [below of=MRP0,node distance=1.1in] (WL0) {};
	\node [above of=MRP1, node distance=1 in] (WL1) {};
	\draw [name path=WL, line width=0.3mm,black!60!green] (WL0.center) -- (WL1.center) node[above left]{\footnotesize $W(L)$};
	
	\path [name intersections={of=MRP and WL,by=E}];
	\node [fill=black,inner sep=2pt,label=-0:$E$] at (E) {};
	\node [below of=E, node distance=0.645 in] (Lc) {};
	\node [left of=E, node distance=1 in] (Wc) {};

	\node [above of=E, node distance=0.8 in] (MCL1) {};
	\draw [name path=MCL, line width=0.3mm,black!60!blue,dotted] (WL0.center) -- (MCL1.east);
	\node [left of=MCL1, node distance=10pt,black!60!blue] (MCLlabel) {\footnotesize MCL};
	
	\path [name intersections={of=MRP and MCL,by=A}];
	\node [fill=black,inner sep=2pt,label=-180:$A$] at (A) {};
	\node [below of=A, node distance=0.86 in] (L1) {};
	\draw [name path=L1path,line width=0.3mm,dotted] (A.center) -- (L1.center) node[below]{};
	\path [name intersections={of=WL and L1path,by=AA}];
	\node [fill=black,inner sep=2pt,label=-180:$\hspace*{-5pt}A'$] at (AA) {};
	\node [left of=AA, node distance=0.66 in] (W1) {};
	\node [above of=W1, node distance=0.1 in] (Wm) {};
	\draw [name path=MLCm1,opacity=0] (Wm.center) --+(1.5in,0);
	\path [name intersections={of=WL and MLCm1,by=M1}];
	\draw [line width=0.3mm,black!60!blue] (Wm.center) node[left]{\footnotesize $\omega_m$} -- (M1.west) ;
	\draw [name path=MLCm2,opacity=0] (M1.center) --+(0,1.5in);
	\path [name intersections={of=MCL and MLCm2,by=M2}];
	\draw [line width=0.3mm,black!60!blue] (M1.center) -- (M2.center) ;
	\draw [line width=0.3mm,black!60!blue,one end extended] (M2.center) -- (MCL1.east) ;
	\path [name intersections={of=MRP and MLCm2,by=B}];
	\node [fill=black,inner sep=2pt,label=5:$\hspace*{-5pt}B$] at (B.east) {};
	\node [below of=B, node distance=0.77 in] (L2) {};
	\draw [name path=L2path,line width=0.3mm,dotted] (B.center) -- (L2.center) node[below]{\footnotesize $L^*_2$};
	\node [left of=B, node distance=0.84 in] (W2) {};
	\draw [name path=W2path,line width=0.3mm,dotted] (B.center) -- (W2.east) node[left]{\footnotesize $\omega^*_2$};
	

\end{tikzpicture}						

%% file: Monopsony-Fig3.tex
\hspace*{-25pt}\begin{tikzpicture}[dot/.style={circle,inner sep=1pt,fill,label={#1},name=#1},
extended line/.style={shorten >=-#1,shorten <=-#1},
extended line/.default=1cm,
one end extended/.style={shorten >=-#1},
one end extended/.default=1cm,
]
	\node(origin) at (0,0) {};
	\node [left of=origin,node distance=3pt] (originol){};
	\node [below of=origin,node distance=1.5pt] (origino){};
	\node [left of=origin,node distance=8pt] (originleft){};
	\node [below of=origin,node distance=8pt] (originbelow){};
	\node [right of=origin,node distance=1.6in] (xaxis) {};
	\node [above of=origin,node distance=1.2in] (yaxis) {Wage};
	\node [below of=xaxis, node distance=10pt] (xlabel) {\textbf{\#} Drivers};	
	\draw [->,line width=0.3mm, black] (originleft) -- (xaxis);
	\draw [->,line width=0.3mm, black] (originbelow) -- (yaxis);

	\node [right of=yaxis,node distance=15pt] (MRP0) {};
	\node [above of=xaxis, node distance=10pt] (MRP1) {};
	\draw [name path=MRP,line width=0.3mm, red] (MRP0.center) node[above right]{\footnotesize MRP} -- (MRP1.center) ;

	\node [below of=MRP0,node distance=1.1in] (WL0) {};
	\node [above of=MRP1, node distance=1 in] (WL1) {};
	\draw [name path=WL, line width=0.3mm,black!60!green] (WL0.center) -- (WL1.center) node[above left]{\footnotesize $W(L)$};
	
	\path [name intersections={of=MRP and WL,by=E}];
	\node [fill=black,inner sep=2pt,label=-0:$E$] at (E) {};
	\node [below of=E, node distance=0.645 in] (Lc) {};
	\node [left of=E, node distance=1 in] (Wc) {};

	\node [above of=E, node distance=0.8 in] (MCL1) {};
	\draw [name path=MCL, line width=0.3mm,blue!60!green,dotted,one end extended] (WL0.center) -- (MCL1.east);
	\node [left of=MCL1, node distance=10pt,black!60!blue] (MCLlabel) {};
	
	\path [name intersections={of=MRP and MCL,by=A}];
	\node [fill=black,inner sep=2pt,label=-180:$A$] at (A) {};
	\node [below of=A, node distance=0.86 in] (L1) {};
	\draw [name path=L1path,line width=0.3mm,dotted,opacity=0] (A.center) -- (L1.center) node[below]{};
	\path [name intersections={of=WL and L1path,by=AA}];
	\node [fill=black,inner sep=2pt,label=-180:$\hspace*{-5pt}A'$] at (AA) {};
	
	\node [left of=AA, node distance=0.66 in] (W1) {};
	\node [above of=W1, node distance=0.32 in] (Wm) {};
	\draw [name path=MLCm1,opacity=0] (Wm.center) --+(1.5in,0);
	\path [name intersections={of=WL and MLCm1,by=M1}];
	\draw [line width=0.3mm,black!60!blue] (Wm.center) node[left]{\footnotesize $\omega_3^*\hspace*{-1pt}=\hspace*{-1pt}\omega_m$} -- (M1.west) ;
	
	\draw [name path=MLCm2,opacity=0] (M1.center) --+(0,1.5in);
	\path [name intersections={of=MCL and MLCm2,by=M2}];
	\draw [line width=0.3mm,black!60!blue] (M1.center) --+ (0,0.865in) node[left](MCLend) {\footnotesize MCL};
	\draw [line width=0.3mm,black!60!blue] (MCLend.east) --+ (0.1in,0.15in) ;
	\path [name intersections={of=MRP and MLCm1,by=C}];
	\node [fill=black,inner sep=2pt,label=5:$\hspace*{-5pt}C$] at (C.east) {};
	\node [below of=B, node distance=0.77 in] (L2) {};
	\draw [name path=L2path,line width=0.3mm,dotted] (C.center) -- (L2.center) node[below]{\footnotesize $L^*_3$};
	\node [left of=B, node distance=0.84 in] (W2) {};
	

\end{tikzpicture}						

%% file: figure1_cap.tex
%
%
\begin{tikzpicture}

\begin{axis}[%
width=1.794in,
height=1.03in,
at={(1.358in,0.0in)},
scale only axis,
xmin=3000,
xmax=6000,
xtick={3000,4000,5000,6000},
xticklabels={{3k},{4k},{5k},{6k}},
xlabel style={font=\color{white!15!black}},
xlabel={Cap},
ymin=3000,
ymax=6000,
ytick={3000,4000,5000,6000},
yticklabels={{3k},{4k},{5k},{6k}},
ylabel style={font=\color{white!15!black}},
ylabel={Number of Drivers},
axis background/.style={fill=white},
legend style={legend cell align=left, align=left, draw=white!15!black}
]
\addplot [color=black, line width=1.0pt]
  table[row sep=crcr]{%
3000	3000\\
3030.30303030303	3030.30303030303\\
3060.60606060606	3060.60606060606\\
3090.90909090909	3090.90909090909\\
3121.21212121212	3121.21212121212\\
3151.51515151515	3151.51515151515\\
3181.81818181818	3181.81818181818\\
3212.12121212121	3212.12121212121\\
3242.42424242424	3242.42424242424\\
3272.72727272727	3272.72727272727\\
3303.0303030303	3303.0303030303\\
3333.33333333333	3333.33333333333\\
3363.63636363636	3363.63636363636\\
3393.93939393939	3393.93939393939\\
3424.24242424242	3424.24242424242\\
3454.54545454545	3454.54545454545\\
3484.84848484848	3484.84848484848\\
3515.15151515152	3515.15151515152\\
3545.45454545455	3545.45454545455\\
3575.75757575758	3575.75757575758\\
3606.06060606061	3606.06060606061\\
3636.36363636364	3636.36363636364\\
3666.66666666667	3666.66666666667\\
3696.9696969697	3696.9696969697\\
3727.27272727273	3727.27272727273\\
3757.57575757576	3757.57575757576\\
3787.87878787879	3787.87878787879\\
3818.18181818182	3818.18181818182\\
3848.48484848485	3848.48484848485\\
3878.78787878788	3878.78787878788\\
3909.09090909091	3909.09090909091\\
3939.39393939394	3939.39393939394\\
3969.69696969697	3969.69696969697\\
4000	4000\\
4030.30303030303	4030.30303030303\\
4060.60606060606	4060.60606060606\\
4090.90909090909	4090.90909090909\\
4121.21212121212	4121.21212121212\\
4151.51515151515	4151.51515151515\\
4181.81818181818	4181.81818181818\\
4212.12121212121	4212.12121212121\\
4242.42424242424	4242.42424242424\\
4272.72727272727	4272.72727272727\\
4303.0303030303	4303.0303030303\\
4333.33333333333	4333.33333333333\\
4363.63636363636	4363.63636363636\\
4393.93939393939	4393.93939393939\\
4424.24242424242	4424.24242424242\\
4454.54545454545	4454.54545454545\\
4484.84848484848	4484.84848484848\\
4515.15151515152	4515.15151515152\\
4545.45454545455	4545.45454545455\\
4575.75757575758	4575.75757575758\\
4606.06060606061	4606.06060606061\\
4636.36363636364	4636.36363636364\\
4666.66666666667	4666.66666666667\\
4696.9696969697	4696.9696969697\\
4727.27272727273	4727.27272727273\\
4757.57575757576	4757.57575757576\\
4787.87878787879	4787.87878787879\\
4818.18181818182	4818.18181818182\\
4848.48484848485	4848.48484848485\\
4878.78787878788	4878.78787878788\\
4909.09090909091	4909.09090909091\\
4939.39393939394	4939.39393939394\\
4969.69696969697	4969.69696969697\\
5000	5000\\
5030.30303030303	5030.30303030303\\
5060.60606060606	5060.60606060606\\
5090.90909090909	5089.33333333333\\
5121.21212121212	5089.33333333333\\
5151.51515151515	5089.33333333333\\
5181.81818181818	5089.33333333333\\
5212.12121212121	5089.33333333333\\
5242.42424242424	5089.33333333333\\
5272.72727272727	5089.33333333333\\
5303.0303030303	5089.33333333333\\
5333.33333333333	5089.33333333333\\
5363.63636363636	5089.33333333333\\
5393.93939393939	5089.33333333333\\
5424.24242424242	5089.33333333333\\
5454.54545454545	5089.33333333333\\
5484.84848484848	5089.33333333333\\
5515.15151515152	5089.33333333333\\
5545.45454545455	5089.33333333333\\
5575.75757575758	5089.33333333333\\
5606.06060606061	5089.33333333333\\
5636.36363636364	5089.33333333333\\
5666.66666666667	5089.33333333333\\
5696.9696969697	5089.33333333333\\
5727.27272727273	5089.33333333333\\
5757.57575757576	5089.33333333333\\
5787.87878787879	5089.33333333333\\
5818.18181818182	5089.33333333333\\
5848.48484848485	5089.33333333333\\
5878.78787878788	5089.33333333333\\
5909.09090909091	5089.33333333333\\
5939.39393939394	5089.33333333333\\
5969.69696969697	5089.33333333333\\
6000	5089.33333333333\\
};

\end{axis}
\end{tikzpicture}%

%% file: figure2_cap.tex
%
%
\begin{tikzpicture}

\pgfplotsset{every axis/.append style={
font=\small,
thin,
tick style={ultra thin}}}
\pgfplotsset{every axis y label/.style={
at={(-0.42,0.5)},
xshift=32pt,
rotate=90}}

\begin{axis}[%
width=1.794in,
height=1.03in,
at={(1.358in,0.0in)},
scale only axis,
xmin=3000,
xmax=6000,
xtick={3000,4000,5000,6000},
xticklabels={{3k},{4k},{5k},{6k}},
xlabel style={font=\color{white!15!black}},
xlabel={Cap},
ymin=100,
ymax=190,
ylabel style={font=\color{white!15!black}},
ylabel={ Ride Arrival/min},
axis background/.style={fill=white},
legend style={legend cell align=left, align=left, draw=white!15!black}
]
\addplot [color=black, line width=1.0pt]
  table[row sep=crcr]{%
3000	102.331752015664\\
3030.30303030303	103.580517526981\\
3060.60606060606	104.829198494105\\
3090.90909090909	106.077771619468\\
3121.21212121212	107.326214046798\\
3151.51515151515	108.574503345954\\
3181.81818181818	109.822617498422\\
3212.12121212121	111.070534883427\\
3242.42424242424	112.318234264632\\
3272.72727272727	113.565694777395\\
3303.0303030303	114.812895916556\\
3333.33333333333	116.059817524722\\
3363.63636363636	117.306439781043\\
3393.93939393939	118.552743190426\\
3424.24242424242	119.798708573196\\
3454.54545454545	121.044317055166\\
3484.84848484848	122.289550058091\\
3515.15151515152	123.534389290513\\
3545.45454545455	124.778816738947\\
3575.75757575758	126.022814659414\\
3606.06060606061	127.266365569304\\
3636.36363636364	128.50945223954\\
3666.66666666667	129.752057687041\\
3696.9696969697	130.994165167475\\
3727.27272727273	132.235758168268\\
3757.57575757576	133.47682040189\\
3787.87878787879	134.71733579937\\
3818.18181818182	135.957288504062\\
3848.48484848485	137.196662865628\\
3878.78787878788	138.435443434244\\
3909.09090909091	139.67361495501\\
3939.39393939394	140.911162362562\\
3969.69696969697	142.148070775875\\
4000	143.384325493242\\
4030.30303030303	144.619911987442\\
4060.60606060606	145.854815901062\\
4090.90909090909	147.089023041994\\
4121.21212121212	148.322519379075\\
4151.51515151515	149.555291037887\\
4181.81818181818	150.787324296692\\
4212.12121212121	152.018605582513\\
4242.42424242424	153.249121467338\\
4272.72727272727	154.47885866446\\
4303.0303030303	155.707804024936\\
4333.33333333333	156.935944534163\\
4363.63636363636	158.16326730857\\
4393.93939393939	159.389759592415\\
4424.24242424242	160.615408754698\\
4454.54545454545	161.84020228616\\
4484.84848484848	163.064127796394\\
4515.15151515152	164.287173011042\\
4545.45454545455	165.509325769086\\
4575.75757575758	166.730574020228\\
4606.06060606061	167.950905822351\\
4636.36363636364	169.170309339064\\
4666.66666666667	170.388772837328\\
4696.9696969697	171.606284685154\\
4727.27272727273	172.822833349379\\
4757.57575757576	174.038407393508\\
4787.87878787879	175.252995475632\\
4818.18181818182	176.466586346404\\
4848.48484848485	177.679168847085\\
4878.78787878788	178.890731907653\\
4909.09090909091	180.101264544968\\
4939.39393939394	181.310755860997\\
4969.69696969697	182.519195041095\\
5000	183.726571352344\\
5030.30303030303	184.93287414194\\
5060.60606060606	186.138092835633\\
5090.90909090909	187.27962962963\\
5121.21212121212	187.27962962963\\
5151.51515151515	187.27962962963\\
5181.81818181818	187.27962962963\\
5212.12121212121	187.27962962963\\
5242.42424242424	187.27962962963\\
5272.72727272727	187.27962962963\\
5303.0303030303	187.27962962963\\
5333.33333333333	187.27962962963\\
5363.63636363636	187.27962962963\\
5393.93939393939	187.27962962963\\
5424.24242424242	187.27962962963\\
5454.54545454545	187.27962962963\\
5484.84848484848	187.27962962963\\
5515.15151515152	187.27962962963\\
5545.45454545455	187.27962962963\\
5575.75757575758	187.27962962963\\
5606.06060606061	187.27962962963\\
5636.36363636364	187.27962962963\\
5666.66666666667	187.27962962963\\
5696.9696969697	187.27962962963\\
5727.27272727273	187.27962962963\\
5757.57575757576	187.27962962963\\
5787.87878787879	187.27962962963\\
5818.18181818182	187.27962962963\\
5848.48484848485	187.27962962963\\
5878.78787878788	187.27962962963\\
5909.09090909091	187.27962962963\\
5939.39393939394	187.27962962963\\
5969.69696969697	187.27962962963\\
6000	187.27962962963\\
};

\end{axis}
\end{tikzpicture}%

%% file: figure3_cap.tex
%
%
\begin{tikzpicture}

\pgfplotsset{every axis/.append style={
font=\small,
thin,
tick style={ultra thin}}}
\pgfplotsset{every axis y label/.style={
at={(-0.45,0.5)},
xshift=32pt,
rotate=90}}

\begin{axis}[%
width=1.794in,
height=1.03in,
at={(1.358in,0.0in)},
scale only axis,
xmin=3000,
xmax=6000,
xtick={3000,4000,5000,6000},
xticklabels={{3k},{4k},{5k},{6k}},
xlabel style={font=\color{white!15!black}},
xlabel={Cap},
ymin=0.555,
ymax=0.6,
ylabel style={font=\color{white!15!black}},
ylabel={Occupancy Rate},
axis background/.style={fill=white},
legend style={legend cell align=left, align=left, draw=white!15!black}
]
\addplot [color=black, line width=1.0pt]
  table[row sep=crcr]{%
3000	0.556002519285108\\
3030.30303030303	0.557159603777633\\
3060.60606060606	0.558293325445335\\
3090.90909090909	0.559404248569724\\
3121.21212121212	0.560492917823035\\
3151.51515151515	0.561559859132583\\
3181.81818181818	0.562605580499059\\
3212.12121212121	0.563630572771655\\
3242.42424242424	0.56463531038267\\
3272.72727272727	0.565620252044083\\
3303.0303030303	0.566585841408397\\
3333.33333333333	0.567532507695892\\
3363.63636363636	0.568460666290297\\
3393.93939393939	0.569370719304732\\
3424.24242424242	0.570263056119666\\
3454.54545454545	0.571138053894505\\
3484.84848484848	0.571996078054325\\
3515.15151515152	0.572837482753165\\
3545.45454545455	0.573662611315209\\
3575.75757575758	0.574471796655075\\
3606.06060606061	0.575265361678393\\
3636.36363636364	0.576043619663737\\
3666.66666666667	0.576806874626938\\
3696.9696969697	0.577555421668726\\
3727.27272727273	0.578289547306597\\
3757.57575757576	0.579009529791748\\
3787.87878787879	0.579715639411851\\
3818.18181818182	0.580408138780437\\
3848.48484848485	0.581087283113554\\
3878.78787878788	0.581753320494372\\
3909.09090909091	0.582406492126354\\
3939.39393939394	0.583047032575555\\
3969.69696969697	0.583675170002618\\
4000	0.584291126384961\\
4030.30303030303	0.584895117729661\\
4060.60606060606	0.585487354277473\\
4090.90909090909	0.586068040698434\\
4121.21212121212	0.586637376279445\\
4151.51515151515	0.587195555104228\\
4181.81818181818	0.587742766226019\\
4212.12121212121	0.588279193833335\\
4242.42424242424	0.588805017409149\\
4272.72727272727	0.58932041188378\\
4303.0303030303	0.589825547781783\\
4333.33333333333	0.590320591363122\\
4363.63636363636	0.590805704758886\\
4393.93939393939	0.591281046101794\\
4424.24242424242	0.591746769651726\\
4454.54545454545	0.5922030259165\\
4484.84848484848	0.59264996176811\\
4515.15151515152	0.593087720554627\\
4545.45454545455	0.593516442207944\\
4575.75757575758	0.593936263347555\\
4606.06060606061	0.594347317380543\\
4636.36363636364	0.594749734597924\\
4666.66666666667	0.595143642267523\\
4696.9696969697	0.595529164723511\\
4727.27272727273	0.595906423452762\\
4757.57575757576	0.596275537178141\\
4787.87878787879	0.596636621938877\\
4818.18181818182	0.596989791168117\\
4848.48484848485	0.597335155767794\\
4878.78787878788	0.597672824180912\\
4909.09090909091	0.598002902461348\\
4939.39393939394	0.598325494341289\\
4969.69696969697	0.598640701296371\\
5000	0.598948622608641\\
5030.30303030303	0.599249355427406\\
5060.60606060606	0.599542994828067\\
5090.90909090909	0.599814899717638\\
5121.21212121212	0.599814899717638\\
5151.51515151515	0.599814899717638\\
5181.81818181818	0.599814899717638\\
5212.12121212121	0.599814899717638\\
5242.42424242424	0.599814899717638\\
5272.72727272727	0.599814899717638\\
5303.0303030303	0.599814899717638\\
5333.33333333333	0.599814899717638\\
5363.63636363636	0.599814899717638\\
5393.93939393939	0.599814899717638\\
5424.24242424242	0.599814899717638\\
5454.54545454545	0.599814899717638\\
5484.84848484848	0.599814899717638\\
5515.15151515152	0.599814899717638\\
5545.45454545455	0.599814899717638\\
5575.75757575758	0.599814899717638\\
5606.06060606061	0.599814899717638\\
5636.36363636364	0.599814899717638\\
5666.66666666667	0.599814899717638\\
5696.9696969697	0.599814899717638\\
5727.27272727273	0.599814899717638\\
5757.57575757576	0.599814899717638\\
5787.87878787879	0.599814899717638\\
5818.18181818182	0.599814899717638\\
5848.48484848485	0.599814899717638\\
5878.78787878788	0.599814899717638\\
5909.09090909091	0.599814899717638\\
5939.39393939394	0.599814899717638\\
5969.69696969697	0.599814899717638\\
6000	0.599814899717638\\
};

\end{axis}
\end{tikzpicture}%

%% file: figure4_cap.tex
%
%
\definecolor{mycolor1}{rgb}{0.60000,0.20000,0.00000}%
\begin{tikzpicture}

\pgfplotsset{every axis y label/.style={
at={(-0.42,0.5)},
xshift=32pt,
rotate=90}}

\begin{axis}[%
width=1.794in,
height=1.03in,
at={(1.358in,0.0in)},
scale only axis,
xmin=3000,
xmax=6000,
xtick={3000,4000,5000,6000},
xticklabels={{3k},{4k},{5k},{6k}},
xlabel style={font=\color{white!15!black}},
xlabel={Cap},
ymin=6,
ymax=18,
ylabel style={font=\color{white!15!black}},
ylabel={Price (\$/trip)},
axis background/.style={fill=white},
legend style={at={(0.0,0.3)}, anchor=south west, legend cell align=left, align=left, draw=white!15!black}
]
\addplot [color=black, line width=1.0pt]
  table[row sep=crcr]{%
3000	15.265427531637\\
3030.30303030303	15.309189806976\\
3060.60606060606	15.3521324324423\\
3090.90909090909	15.3942736754464\\
3121.21212121212	15.4356312435769\\
3151.51515151515	15.4762223070404\\
3181.81818181818	15.5160635199973\\
3212.12121212121	15.5551710408581\\
3242.42424242424	15.5935605516008\\
3272.72727272727	15.6312472761643\\
3303.0303030303	15.66824599797\\
3333.33333333333	15.7045710766211\\
3363.63636363636	15.7402364638245\\
3393.93939393939	15.7752557185788\\
3424.24242424242	15.8096420216669\\
3454.54545454545	15.8434081894936\\
3484.84848484848	15.8765666872993\\
3515.15151515152	15.9091296417869\\
3545.45454545455	15.9411088531903\\
3575.75757575758	15.9725158068138\\
3606.06060606061	16.0033616840719\\
3636.36363636364	16.0336573730519\\
3666.66666666667	16.0634134786261\\
3696.9696969697	16.0926403321358\\
3727.27272727273	16.1213480006679\\
3757.57575757576	16.1495462959445\\
3787.87878787879	16.1772447828455\\
3818.18181818182	16.2044527875805\\
3848.48484848485	16.2311794055294\\
3878.78787878788	16.2574335087646\\
3909.09090909091	16.2832237532728\\
3939.39393939394	16.3085585858889\\
3969.69696969697	16.3334462509565\\
4000	16.3578947967272\\
4030.30303030303	16.3819120815105\\
4060.60606060606	16.4055057795877\\
4090.90909090909	16.4286833868973\\
4121.21212121212	16.4514522265056\\
4151.51515151515	16.4738194538696\\
4181.81818181818	16.4957920619035\\
4212.12121212121	16.5173768858558\\
4242.42424242424	16.5385806080063\\
4272.72727272727	16.5594097621909\\
4303.0303030303	16.5798707381611\\
4333.33333333333	16.5999697857862\\
4363.63636363636	16.6197130191034\\
4393.93939393939	16.6391064202247\\
4424.24242424242	16.6581558431039\\
4454.54545454545	16.6768670171715\\
4484.84848484848	16.6952455508428\\
4515.15151515152	16.7132969349029\\
4545.45454545455	16.7310265457759\\
4575.75757575758	16.7484396486814\\
4606.06060606061	16.7655414006837\\
4636.36363636364	16.7823368536371\\
4666.66666666667	16.7988309570329\\
4696.9696969697	16.8150285607507\\
4727.27272727273	16.8309344177186\\
4757.57575757576	16.8465531864851\\
4787.87878787879	16.8618894337068\\
4818.18181818182	16.8769476365555\\
4848.48484848485	16.8917321850466\\
4878.78787878788	16.9062473842926\\
4909.09090909091	16.9204974566847\\
4939.39393939394	16.9344865440049\\
4969.69696969697	16.9482187094712\\
5000	16.9616979397185\\
5030.30303030303	16.9749281467175\\
5060.60606060606	16.9879131696342\\
5090.90909090909	17\\
5121.21212121212	17\\
5151.51515151515	17\\
5181.81818181818	17\\
5212.12121212121	17\\
5242.42424242424	17\\
5272.72727272727	17\\
5303.0303030303	17\\
5333.33333333333	17\\
5363.63636363636	17\\
5393.93939393939	17\\
5424.24242424242	17\\
5454.54545454545	17\\
5484.84848484848	17\\
5515.15151515152	17\\
5545.45454545455	17\\
5575.75757575758	17\\
5606.06060606061	17\\
5636.36363636364	17\\
5666.66666666667	17\\
5696.9696969697	17\\
5727.27272727273	17\\
5757.57575757576	17\\
5787.87878787879	17\\
5818.18181818182	17\\
5848.48484848485	17\\
5878.78787878788	17\\
5909.09090909091	17\\
5939.39393939394	17\\
5969.69696969697	17\\
6000	17\\
};
\addlegendentry{$p_f$}

\addplot [color=mycolor1, dotted, line width=1.0pt]
  table[row sep=crcr]{%
3000	6.5092116330656\\
3030.30303030303	6.56130665598388\\
3060.60606060606	6.61346249112585\\
3090.90909090909	6.66567857963095\\
3121.21212121212	6.71795440669905\\
3151.51515151515	6.77028949906186\\
3181.81818181818	6.82268342262474\\
3212.12121212121	6.87513578026598\\
3242.42424242424	6.92764620978155\\
3272.72727272727	6.98021438196441\\
3303.0303030303	7.03283999880852\\
3333.33333333333	7.08552279182823\\
3363.63636363636	7.13826252048471\\
3393.93939393939	7.19105897071177\\
3424.24242424242	7.2439119535339\\
3454.54545454545	7.29682130376995\\
3484.84848484848	7.34978687881663\\
3515.15151515152	7.40280855750606\\
3545.45454545455	7.45588623903244\\
3575.75757575758	7.50901984194299\\
3606.06060606061	7.56220930318897\\
3636.36363636364	7.61545457723257\\
3666.66666666667	7.66875563520609\\
3696.9696969697	7.72211246411984\\
3727.27272727273	7.77552506611574\\
3757.57575757576	7.82899345776346\\
3787.87878787879	7.8825176693965\\
3818.18181818182	7.93609774448561\\
3848.48484848485	7.9897337390472\\
3878.78787878788	8.0434257210845\\
3909.09090909091	8.09717377005947\\
3939.39393939394	8.15097797639357\\
3969.69696969697	8.20483844099555\\
4000	8.25875527481468\\
4030.30303030303	8.31272859841784\\
4060.60606060606	8.36675854158908\\
4090.90909090909	8.4208452429502\\
4121.21212121212	8.47498884960132\\
4151.51515151515	8.52918951677998\\
4181.81818181818	8.58344740753794\\
4212.12121212121	8.63776269243449\\
4242.42424242424	8.69213554924541\\
4272.72727272727	8.74656616268667\\
4303.0303030303	8.80105472415193\\
4333.33333333333	8.85560143146321\\
4363.63636363636	8.91020648863397\\
4393.93939393939	8.96487010564369\\
4424.24242424242	9.01959249822371\\
4454.54545454545	9.07437388765327\\
4484.84848484848	9.12921450056556\\
4515.15151515152	9.18411456876296\\
4545.45454545455	9.23907432904121\\
4575.75757575758	9.29409402302177\\
4606.06060606061	9.3491738969922\\
4636.36363636364	9.40431420175391\\
4666.66666666667	9.45951519247702\\
4696.9696969697	9.51477712856194\\
4727.27272727273	9.57010027350721\\
4757.57575757576	9.62548489478347\\
4787.87878787879	9.68093126371304\\
4818.18181818182	9.73643965535495\\
4848.48484848485	9.79201034839509\\
4878.78787878788	9.84764362504125\\
4909.09090909091	9.90333977092268\\
4939.39393939394	9.95909907499421\\
4969.69696969697	10.0149218294444\\
5000	10.0708083296075\\
5030.30303030303	10.1267588738796\\
5060.60606060606	10.1827737636381\\
5090.90909090909	10.2359355687178\\
5121.21212121212	10.2359355687178\\
5151.51515151515	10.2359355687178\\
5181.81818181818	10.2359355687178\\
5212.12121212121	10.2359355687178\\
5242.42424242424	10.2359355687178\\
5272.72727272727	10.2359355687178\\
5303.0303030303	10.2359355687178\\
5333.33333333333	10.2359355687178\\
5363.63636363636	10.2359355687178\\
5393.93939393939	10.2359355687178\\
5424.24242424242	10.2359355687178\\
5454.54545454545	10.2359355687178\\
5484.84848484848	10.2359355687178\\
5515.15151515152	10.2359355687178\\
5545.45454545455	10.2359355687178\\
5575.75757575758	10.2359355687178\\
5606.06060606061	10.2359355687178\\
5636.36363636364	10.2359355687178\\
5666.66666666667	10.2359355687178\\
5696.9696969697	10.2359355687178\\
5727.27272727273	10.2359355687178\\
5757.57575757576	10.2359355687178\\
5787.87878787879	10.2359355687178\\
5818.18181818182	10.2359355687178\\
5848.48484848485	10.2359355687178\\
5878.78787878788	10.2359355687178\\
5909.09090909091	10.2359355687178\\
5939.39393939394	10.2359355687178\\
5969.69696969697	10.2359355687178\\
6000	10.2359355687178\\
};
\addlegendentry{$p_d$}

\end{axis}
\end{tikzpicture}%

%% file: commission_cap.tex
%
%
\begin{tikzpicture}

\pgfplotsset{every axis/.append style={
font=\small,
thin,
tick style={ultra thin}}}
\pgfplotsset{every axis y label/.style={
at={(-0.45,0.5)},
xshift=32pt,
rotate=90}}

\begin{axis}[%
width=1.794in,
height=1.03in,
at={(1.358in,0in)},
scale only axis,
xmin=3000,
xmax=6000,
xtick={3000,4000,5000,6000},
xticklabels={{3k},{4k},{5k},{6k}},
xlabel style={font=\color{white!15!black}},
xlabel={Cap},
ymin=0.38,
ymax=0.58,
ylabel style={font=\color{white!15!black}},
ylabel={Commisison Rate},
axis background/.style={fill=white},
legend style={legend cell align=left, align=left, draw=white!15!black}
]
\addplot [color=black, line width=1.0pt]
  table[row sep=crcr]{%
3000	0.573597816400786\\
3030.30303030303	0.571413854115646\\
3060.60606060606	0.569215382929462\\
3090.90909090909	0.567002723209825\\
3121.21212121212	0.564776179173459\\
3151.51515151515	0.562536039820135\\
3181.81818181818	0.56028257980308\\
3212.12121212121	0.558016060240844\\
3242.42424242424	0.555736729475145\\
3272.72727272727	0.55344482377882\\
3303.0303030303	0.551140568017652\\
3333.33333333333	0.548824176269531\\
3363.63636363636	0.546495852404096\\
3393.93939393939	0.544155790625774\\
3424.24242424242	0.541804175982847\\
3454.54545454545	0.539441184845016\\
3484.84848484848	0.537066985351675\\
3515.15151515152	0.534681737832984\\
3545.45454545455	0.532285595205612\\
3575.75757575758	0.529878703344926\\
3606.06060606061	0.527461201435221\\
3636.36363636364	0.525033222299484\\
3666.66666666667	0.522594892710064\\
3696.9696969697	0.520146333681531\\
3727.27272727273	0.517687660746881\\
3757.57575757576	0.515218984218184\\
3787.87878787879	0.512740409432686\\
3818.18181818182	0.510252036985289\\
3848.48484848485	0.507753962948289\\
3878.78787878788	0.505246279079156\\
3909.09090909091	0.502729073017129\\
3939.39393939394	0.500202428469291\\
3969.69696969697	0.497666425386801\\
4000	0.495121140131854\\
4030.30303030303	0.49256664563595\\
4060.60606060606	0.490003011549983\\
4090.90909090909	0.487430304386641\\
4121.21212121212	0.484848587655566\\
4151.51515151515	0.4822579219917\\
4181.81818181818	0.479658365277219\\
4212.12121212121	0.477049972757405\\
4242.42424242424	0.474432797150828\\
4272.72727272727	0.471806888754139\\
4303.0303030303	0.469172295541789\\
4333.33333333333	0.466529063260953\\
4363.63636363636	0.463877235521925\\
4393.93939393939	0.461216853884234\\
4424.24242424242	0.458547957938715\\
4454.54545454545	0.455870585385749\\
4484.84848484848	0.453184772109884\\
4515.15151515152	0.450490552251035\\
4545.45454545455	0.447787958272422\\
4575.75757575758	0.44507702102545\\
4606.06060606061	0.44235776981166\\
4636.36363636364	0.439630232441926\\
4666.66666666667	0.436894435293024\\
4696.9696969697	0.434150403361722\\
4727.27272727273	0.431398160316497\\
4757.57575757576	0.428637728547025\\
4787.87878787879	0.425869129211527\\
4818.18181818182	0.423092382282103\\
4848.48484848485	0.420307506588137\\
4878.78787878788	0.417514519857874\\
4909.09090909091	0.414713438758255\\
4939.39393939394	0.411904278933104\\
4969.69696969697	0.409087055039731\\
5000	0.40626178078404\\
5030.30303030303	0.403428468954203\\
5060.60606060606	0.400587131452982\\
5090.90909090909	0.3978861430166\\
5121.21212121212	0.3978861430166\\
5151.51515151515	0.3978861430166\\
5181.81818181818	0.3978861430166\\
5212.12121212121	0.3978861430166\\
5242.42424242424	0.3978861430166\\
5272.72727272727	0.3978861430166\\
5303.0303030303	0.3978861430166\\
5333.33333333333	0.3978861430166\\
5363.63636363636	0.3978861430166\\
5393.93939393939	0.3978861430166\\
5424.24242424242	0.3978861430166\\
5454.54545454545	0.3978861430166\\
5484.84848484848	0.3978861430166\\
5515.15151515152	0.3978861430166\\
5545.45454545455	0.3978861430166\\
5575.75757575758	0.3978861430166\\
5606.06060606061	0.3978861430166\\
5636.36363636364	0.3978861430166\\
5666.66666666667	0.3978861430166\\
5696.9696969697	0.3978861430166\\
5727.27272727273	0.3978861430166\\
5757.57575757576	0.3978861430166\\
5787.87878787879	0.3978861430166\\
5818.18181818182	0.3978861430166\\
5848.48484848485	0.3978861430166\\
5878.78787878788	0.3978861430166\\
5909.09090909091	0.3978861430166\\
5939.39393939394	0.3978861430166\\
5969.69696969697	0.3978861430166\\
6000	0.3978861430166\\
};

\end{axis}
\end{tikzpicture}%

%% file: figure6_cap.tex
%
%
\begin{tikzpicture}

\pgfplotsset{every axis/.append style={
font=\small,
thin,
tick style={ultra thin}}}
\pgfplotsset{every axis y label/.style={
at={(-0.4,0.5)},
xshift=32pt,
rotate=90}}

\begin{axis}[%
width=1.794in,
height=1.03in,
at={(1.358in,0.0in)},
scale only axis,
xmin=3000,
xmax=6000,
xtick={3000,4000,5000,6000},
xticklabels={{3k},{4k},{5k},{6k}},
xlabel style={font=\color{white!15!black}},
xlabel={Cap},
ymin=13,
ymax=23,
ylabel style={font=\color{white!15!black}},
ylabel={Driver Wage/hour},
axis background/.style={fill=white},
legend style={legend cell align=left, align=left, draw=white!15!black}
]
\addplot [color=black, line width=1.0pt]
  table[row sep=crcr]{%
3000	13.3219806130469\\
3030.30303030303	13.4565460737847\\
3060.60606060606	13.5911115345226\\
3090.90909090909	13.7256769952604\\
3121.21212121212	13.8602424559983\\
3151.51515151515	13.9948079167361\\
3181.81818181818	14.129373377474\\
3212.12121212121	14.2639388382118\\
3242.42424242424	14.3985042989497\\
3272.72727272727	14.5330697596875\\
3303.0303030303	14.6676352204254\\
3333.33333333333	14.8022006811632\\
3363.63636363636	14.9367661419011\\
3393.93939393939	15.0713316026389\\
3424.24242424242	15.2058970633768\\
3454.54545454545	15.3404625241146\\
3484.84848484848	15.4750279848525\\
3515.15151515152	15.6095934455903\\
3545.45454545455	15.7441589063282\\
3575.75757575758	15.878724367066\\
3606.06060606061	16.0132898278038\\
3636.36363636364	16.1478552885417\\
3666.66666666667	16.2824207492795\\
3696.9696969697	16.4169862100174\\
3727.27272727273	16.5515516707552\\
3757.57575757576	16.6861171314931\\
3787.87878787879	16.8206825922309\\
3818.18181818182	16.9552480529688\\
3848.48484848485	17.0898135137066\\
3878.78787878788	17.2243789744445\\
3909.09090909091	17.3589444351823\\
3939.39393939394	17.4935098959202\\
3969.69696969697	17.628075356658\\
4000	17.7626408173959\\
4030.30303030303	17.8972062781337\\
4060.60606060606	18.0317717388716\\
4090.90909090909	18.1663371996094\\
4121.21212121212	18.3009026603473\\
4151.51515151515	18.4354681210851\\
4181.81818181818	18.5700335818229\\
4212.12121212121	18.7045990425608\\
4242.42424242424	18.8391645032986\\
4272.72727272727	18.9737299640365\\
4303.0303030303	19.1082954247743\\
4333.33333333333	19.2428608855122\\
4363.63636363636	19.37742634625\\
4393.93939393939	19.5119918069879\\
4424.24242424242	19.6465572677257\\
4454.54545454545	19.7811227284636\\
4484.84848484848	19.9156881892014\\
4515.15151515152	20.0502536499393\\
4545.45454545455	20.1848191106771\\
4575.75757575758	20.319384571415\\
4606.06060606061	20.4539500321528\\
4636.36363636364	20.5885154928907\\
4666.66666666667	20.7230809536285\\
4696.9696969697	20.8576464143664\\
4727.27272727273	20.9922118751042\\
4757.57575757576	21.1267773358421\\
4787.87878787879	21.2613427965799\\
4818.18181818182	21.3959082573177\\
4848.48484848485	21.5304737180556\\
4878.78787878788	21.6650391787934\\
4909.09090909091	21.7996046395313\\
4939.39393939394	21.9341701002691\\
4969.69696969697	22.068735561007\\
5000	22.2033010217448\\
5030.30303030303	22.3378664824827\\
5060.60606060606	22.4724319432205\\
5090.90909090909	22.6\\
5121.21212121212	22.6\\
5151.51515151515	22.6\\
5181.81818181818	22.6\\
5212.12121212121	22.6\\
5242.42424242424	22.6\\
5272.72727272727	22.6\\
5303.0303030303	22.6\\
5333.33333333333	22.6\\
5363.63636363636	22.6\\
5393.93939393939	22.6\\
5424.24242424242	22.6\\
5454.54545454545	22.6\\
5484.84848484848	22.6\\
5515.15151515152	22.6\\
5545.45454545455	22.6\\
5575.75757575758	22.6\\
5606.06060606061	22.6\\
5636.36363636364	22.6\\
5666.66666666667	22.6\\
5696.9696969697	22.6\\
5727.27272727273	22.6\\
5757.57575757576	22.6\\
5787.87878787879	22.6\\
5818.18181818182	22.6\\
5848.48484848485	22.6\\
5878.78787878788	22.6\\
5909.09090909091	22.6\\
5939.39393939394	22.6\\
5969.69696969697	22.6\\
6000	22.6\\
};

\end{axis}
\end{tikzpicture}%

%% file: figure5_cap.tex
%
%
\begin{tikzpicture}

\pgfplotsset{every axis/.append style={
font=\small,
thin,
tick style={ultra thin}}}
\pgfplotsset{every axis y label/.style={
at={(-0.4,0.5)},
xshift=32pt,
rotate=90}}

\begin{axis}[%
width=1.794in,
height=1.03in,
at={(1.358in,0.0in)},
scale only axis,
xmin=3000,
xmax=6000,
xtick={3000,4000,5000,6000},
xticklabels={{3k},{4k},{5k},{6k}},
xlabel style={font=\color{white!15!black}},
xlabel={Cap},
ymin=5,
ymax=6.2,
ylabel style={font=\color{white!15!black}},
ylabel={Pickup Time (min)},
axis background/.style={fill=white},
legend style={legend cell align=left, align=left, draw=white!15!black}
]
\addplot [color=black, line width=1.0pt]
  table[row sep=crcr]{%
3000	6.18272641970833\\
3030.30303030303	6.15976671452138\\
3060.60606060606	6.13705778051081\\
3090.90909090909	6.1145942224742\\
3121.21212121212	6.09237081264752\\
3151.51515151515	6.07038248397423\\
3181.81818181818	6.04862432370609\\
3212.12121212121	6.02709156731633\\
3242.42424242424	6.00577959270727\\
3272.72727272727	5.9846839146955\\
3303.0303030303	5.96380017975889\\
3333.33333333333	5.94312416103088\\
3363.63636363636	5.92265175352818\\
3393.93939393939	5.90237896959901\\
3424.24242424242	5.88230193458001\\
3454.54545454545	5.8624168826503\\
3484.84848484848	5.84272015287229\\
3515.15151515152	5.82320818540919\\
3545.45454545455	5.80387751791003\\
3575.75757575758	5.7847247820532\\
3606.06060606061	5.76574670024059\\
3636.36363636364	5.74694008243435\\
3666.66666666667	5.72830182312906\\
3696.9696969697	5.70982889845251\\
3727.27272727273	5.69151836338854\\
3757.57575757576	5.67336734911603\\
3787.87878787879	5.65537306045804\\
3818.18181818182	5.63753277343597\\
3848.48484848485	5.61984383292347\\
3878.78787878788	5.60230365039532\\
3909.09090909091	5.5849097017668\\
3939.39393939394	5.56765952531913\\
3969.69696969697	5.55055071970704\\
4000	5.53358094204456\\
4030.30303030303	5.51674790606543\\
4060.60606060606	5.50004938035464\\
4090.90909090909	5.48348318664788\\
4121.21212121212	5.46704719819591\\
4151.51515151515	5.45073933819064\\
4181.81818181818	5.43455757825056\\
4212.12121212121	5.41849993696247\\
4242.42424242424	5.40256447847734\\
4272.72727272727	5.38674931115781\\
4303.0303030303	5.37105258627498\\
4333.33333333333	5.35547249675254\\
4363.63636363636	5.34000727595616\\
4393.93939393939	5.32465519652604\\
4424.24242424242	5.30941456925112\\
4454.54545454545	5.2942837419828\\
4484.84848484848	5.27926109858694\\
4515.15151515152	5.26434505793213\\
4545.45454545455	5.24953407291317\\
4575.75757575758	5.23482662950785\\
4606.06060606061	5.22022124586613\\
4636.36363636364	5.20571647143008\\
4666.66666666667	5.1913108860835\\
4696.9696969697	5.17700309932997\\
4727.27272727273	5.1627917494983\\
4757.57575757576	5.14867550297414\\
4787.87878787879	5.13465305345694\\
4818.18181818182	5.12072312124098\\
4848.48484848485	5.10688445251991\\
4878.78787878788	5.09313581871352\\
4909.09090909091	5.07947601581624\\
4939.39393939394	5.06590386376619\\
4969.69696969697	5.05241820583435\\
5000	5.03901790803284\\
5030.30303030303	5.02570185854175\\
5060.60606060606	5.01246896715371\\
5090.90909090909	5\\
5121.21212121212	5\\
5151.51515151515	5\\
5181.81818181818	5\\
5212.12121212121	5\\
5242.42424242424	5\\
5272.72727272727	5\\
5303.0303030303	5\\
5333.33333333333	5\\
5363.63636363636	5\\
5393.93939393939	5\\
5424.24242424242	5\\
5454.54545454545	5\\
5484.84848484848	5\\
5515.15151515152	5\\
5545.45454545455	5\\
5575.75757575758	5\\
5606.06060606061	5\\
5636.36363636364	5\\
5666.66666666667	5\\
5696.9696969697	5\\
5727.27272727273	5\\
5757.57575757576	5\\
5787.87878787879	5\\
5818.18181818182	5\\
5848.48484848485	5\\
5878.78787878788	5\\
5909.09090909091	5\\
5939.39393939394	5\\
5969.69696969697	5\\
6000	5\\
};

\end{axis}
\end{tikzpicture}%

%% file: GC_cap.tex
%
%
\begin{tikzpicture}

\pgfplotsset{every axis/.append style={
font=\small,
thin,
tick style={ultra thin}}}
\pgfplotsset{every axis y label/.style={
at={(-0.45,0.5)},
xshift=32pt,
rotate=90}}

\begin{axis}[%
width=1.794in,
height=1.03in,
at={(1.358in,0in)},
scale only axis,
xmin=3000,
xmax=6000,
xtick={3000,4000,5000,6000},
xticklabels={{3k},{4k},{5k},{6k}},
xlabel style={font=\color{white!15!black}},
xlabel={Cap},
ymin=33,
ymax=36,
ylabel style={font=\color{white!15!black}},
ylabel={Total Cost (\$/trip)},
axis background/.style={fill=white},
legend style={legend cell align=left, align=left, draw=white!15!black}
]
\addplot [color=black, line width=1.0pt]
  table[row sep=crcr]{%
3000	35.52631550199\\
3030.30303030303	35.4948384762429\\
3060.60606060606	35.4633635815602\\
3090.90909090909	35.4318914051924\\
3121.21212121212	35.4004225232665\\
3151.51515151515	35.3689575011681\\
3181.81818181818	35.3374968939071\\
3212.12121212121	35.3060412464681\\
3242.42424242424	35.2745910941449\\
3272.72727272727	35.2431469628626\\
3303.0303030303	35.2117093694847\\
3333.33333333333	35.1802788221087\\
3363.63636363636	35.1488558203494\\
3393.93939393939	35.11744085561\\
3424.24242424242	35.0860344113432\\
3454.54545454545	35.0546369633016\\
3484.84848484848	35.0232489797779\\
3515.15151515152	34.991870921836\\
3545.45454545455	34.9605032435331\\
3575.75757575758	34.929146392133\\
3606.06060606061	34.8978008083115\\
3636.36363636364	34.8664669263535\\
3666.66666666667	34.8351451743432\\
3696.9696969697	34.8038359743469\\
3727.27272727273	34.7725397425891\\
3757.57575757576	34.7412568896216\\
3787.87878787879	34.7099878204871\\
3818.18181818182	34.6787329348764\\
3848.48484848485	34.64749262728\\
3878.78787878788	34.6162672871341\\
3909.09090909091	34.5850572989616\\
3939.39393939394	34.5538630425078\\
3969.69696969697	34.5226848928713\\
4000	34.4915232206311\\
4030.30303030303	34.4603783919678\\
4060.60606060606	34.4292507687816\\
4090.90909090909	34.3981407088061\\
4121.21212121212	34.3670485657181\\
4151.51515151515	34.3359746892433\\
4181.81818181818	34.3049194252588\\
4212.12121212121	34.2738831158921\\
4242.42424242424	34.2428660996166\\
4272.72727272727	34.2118687113438\\
4303.0303030303	34.1808912825127\\
4333.33333333333	34.1499341411762\\
4363.63636363636	34.1189976120842\\
4393.93939393939	34.0880820167644\\
4424.24242424242	34.0571876736003\\
4454.54545454545	34.0263148979068\\
4484.84848484848	33.9954640020028\\
4515.15151515152	33.9646352952821\\
4545.45454545455	33.9338290842815\\
4575.75757575758	33.9030456727472\\
4606.06060606061	33.8722853616982\\
4636.36363636364	33.8415484494886\\
4666.66666666667	33.8108352318676\\
4696.9696969697	33.7801460020369\\
4727.27272727273	33.7494810507075\\
4757.57575757576	33.7188406661535\\
4787.87878787879	33.6882251342648\\
4818.18181818182	33.6576347385985\\
4848.48484848485	33.6270697604272\\
4878.78787878788	33.5965304787878\\
4909.09090909091	33.5660171705268\\
4939.39393939394	33.5355301103458\\
4969.69696969697	33.5050695708442\\
5000	33.4746358225614\\
5030.30303030303	33.4442291340176\\
5060.60606060606	33.4138497717527\\
5090.90909090909	33.3850756082045\\
5121.21212121212	33.3850756082045\\
5151.51515151515	33.3850756082045\\
5181.81818181818	33.3850756082045\\
5212.12121212121	33.3850756082045\\
5242.42424242424	33.3850756082045\\
5272.72727272727	33.3850756082045\\
5303.0303030303	33.3850756082045\\
5333.33333333333	33.3850756082045\\
5363.63636363636	33.3850756082045\\
5393.93939393939	33.3850756082045\\
5424.24242424242	33.3850756082045\\
5454.54545454545	33.3850756082045\\
5484.84848484848	33.3850756082045\\
5515.15151515152	33.3850756082045\\
5545.45454545455	33.3850756082045\\
5575.75757575758	33.3850756082045\\
5606.06060606061	33.3850756082045\\
5636.36363636364	33.3850756082045\\
5666.66666666667	33.3850756082045\\
5696.9696969697	33.3850756082045\\
5727.27272727273	33.3850756082045\\
5757.57575757576	33.3850756082045\\
5787.87878787879	33.3850756082045\\
5818.18181818182	33.3850756082045\\
5848.48484848485	33.3850756082045\\
5878.78787878788	33.3850756082045\\
5909.09090909091	33.3850756082045\\
5939.39393939394	33.3850756082045\\
5969.69696969697	33.3850756082045\\
6000	33.3850756082045\\
};

\end{axis}
\end{tikzpicture}%

%% file: figure7_cap.tex
%
%
\begin{tikzpicture}

\pgfplotsset{every axis/.append style={
font=\small,
thin,
tick style={ultra thin}}}
\pgfplotsset{every axis y label/.style={
at={(-0.4,0.5)},
xshift=32pt,
rotate=90}}

\begin{axis}[%
width=1.794in,
height=1.03in,
at={(1.358in,0.0in)},
scale only axis,
xmin=3000,
xmax=6000,
xtick={3000,4000,5000,6000},
xticklabels={{3k},{4k},{5k},{6k}},
xlabel style={font=\color{white!15!black}},
xlabel={Cap},
ymin=50000,
ymax=80000,
ylabel style={font=\color{white!15!black}},
ylabel={Platform Rent},
axis background/.style={fill=white},
legend style={legend cell align=left, align=left, draw=white!15!black}
]
\addplot [color=black, line width=1.0pt]
  table[row sep=crcr]{%
3000	53762.3348356937\\
3030.30303030303	54366.6158427194\\
3060.60606060606	54964.0659511637\\
3090.90909090909	55554.5950279629\\
3121.21212121212	56138.1150111337\\
3151.51515151515	56714.5398473711\\
3181.81818181818	57283.7854321283\\
3212.12121212121	57845.7695520582\\
3242.42424242424	58400.4118297016\\
3272.72727272727	58947.633670315\\
3303.0303030303	59487.3582107358\\
3333.33333333333	60019.5102701899\\
3363.63636363636	60544.0163029496\\
3393.93939393939	61060.8043527577\\
3424.24242424242	61569.8040089361\\
3454.54545454545	62070.9463641033\\
3484.84848484848	62564.1639734259\\
3515.15151515152	63049.3908153394\\
3545.45454545455	63526.56225367\\
3575.75757575758	63995.615001097\\
3606.06060606061	64456.4870838994\\
3636.36363636364	64909.117807926\\
3666.66666666667	65353.4477257442\\
3696.9696969697	65789.4186049101\\
3727.27272727273	66216.9733973178\\
3757.57575757576	66636.0562095801\\
3787.87878787879	67046.6122743988\\
3818.18181818182	67448.5879228843\\
3848.48484848485	67841.930557785\\
3878.78787878788	68226.5886275896\\
3909.09090909091	68602.5116014681\\
3939.39393939394	68969.6499450173\\
3969.69696969697	69327.9550967797\\
4000	69677.3794455047\\
4030.30303030303	70017.8763081231\\
4060.60606060606	70349.3999084079\\
4090.90909090909	70671.9053562958\\
4121.21212121212	70985.3486278417\\
4151.51515151515	71289.6865457844\\
4181.81818181818	71584.8767607009\\
4212.12121212121	71870.8777327247\\
4242.42424242424	72147.64871381\\
4272.72727272727	72415.1497305203\\
4303.0303030303	72673.3415673227\\
4333.33333333333	72922.1857503694\\
4363.63636363636	73161.644531748\\
4393.93939393939	73391.6808741869\\
4424.24242424242	73612.2584361947\\
4454.54545454545	73823.3415576227\\
4484.84848484848	74024.8952456337\\
4515.15151515152	74216.885161063\\
4545.45454545455	74399.2776051583\\
4575.75757575758	74572.039506685\\
4606.06060606061	74735.1384093861\\
4636.36363636364	74888.5424597817\\
4666.66666666667	75032.2203952995\\
4696.9696969697	75166.1415327236\\
4727.27272727273	75290.275756952\\
4757.57575757576	75404.5935100505\\
4787.87878787879	75509.0657805971\\
4818.18181818182	75603.6640933032\\
4848.48484848485	75688.3604989048\\
4878.78787878788	75763.1275643154\\
4909.09090909091	75827.9383630309\\
4939.39393939394	75882.766465779\\
4969.69696969697	75927.5859314059\\
5000	75962.3712979909\\
5030.30303030303	75987.0975741852\\
5060.60606060606	76001.7402307635\\
5090.90909090909	76006.288888889\\
5121.21212121212	76006.288888889\\
5151.51515151515	76006.288888889\\
5181.81818181818	76006.288888889\\
5212.12121212121	76006.288888889\\
5242.42424242424	76006.288888889\\
5272.72727272727	76006.288888889\\
5303.0303030303	76006.288888889\\
5333.33333333333	76006.288888889\\
5363.63636363636	76006.288888889\\
5393.93939393939	76006.288888889\\
5424.24242424242	76006.288888889\\
5454.54545454545	76006.288888889\\
5484.84848484848	76006.288888889\\
5515.15151515152	76006.288888889\\
5545.45454545455	76006.288888889\\
5575.75757575758	76006.288888889\\
5606.06060606061	76006.288888889\\
5636.36363636364	76006.288888889\\
5666.66666666667	76006.288888889\\
5696.9696969697	76006.288888889\\
5727.27272727273	76006.288888889\\
5757.57575757576	76006.288888889\\
5787.87878787879	76006.288888889\\
5818.18181818182	76006.288888889\\
5848.48484848485	76006.288888889\\
5878.78787878788	76006.288888889\\
5909.09090909091	76006.288888889\\
5939.39393939394	76006.288888889\\
5969.69696969697	76006.288888889\\
6000	76006.288888889\\
};

\end{axis}
\end{tikzpicture}%

%% file: figure1_tax.tex
%
%
\begin{tikzpicture}
\pgfplotsset{every axis y label/.style={
at={(-0.49,0.5)},
xshift=32pt,
rotate=90}}

\begin{axis}[%
width=1.794in,
height=1.03in,
at={(1.358in,0.0in)},
scale only axis,
xmin=0,
xmax=3,
xlabel style={font=\color{white!15!black}},
xlabel={Congestion Surcharge (\$/trip)},
ymin=4400,
ymax=5100,
ytick={4400,4600,4800,5000},
yticklabels={{4.4K},{4.6K},{4.8K},{5K}},
ylabel style={font=\color{white!15!black}},
ylabel={Number of Drivers},
axis background/.style={fill=white},
legend style={legend cell align=left, align=left, draw=white!15!black}
]
\addplot [color=black, line width=1.0pt]
  table[row sep=crcr]{%
0	5089.33333333333\\
0.0303030303030303	5082.78555643999\\
0.0606060606060606	5076.23497278673\\
0.0909090909090909	5069.68156914627\\
0.121212121212121	5063.1253321978\\
0.151515151515152	5056.56624852608\\
0.181818181818182	5050.00430462052\\
0.212121212121212	5043.43948687424\\
0.242424242424242	5036.87178158319\\
0.272727272727273	5030.30117494513\\
0.303030303030303	5023.72765305875\\
0.333333333333333	5017.15120192264\\
0.363636363636364	5010.57180743436\\
0.393939393939394	5003.98945538939\\
0.424242424242424	4997.40413148019\\
0.454545454545455	4990.81582129511\\
0.484848484848485	4984.22451031742\\
0.515151515151515	4977.63018392424\\
0.545454545454545	4971.03282738548\\
0.575757575757576	4964.43242586278\\
0.606060606060606	4957.82896440842\\
0.636363636363636	4951.2224279642\\
0.666666666666667	4944.61280136038\\
0.696969696969697	4938.0000693145\\
0.727272727272727	4931.38421643027\\
0.757575757575758	4924.76522719641\\
0.787878787878788	4918.14308598546\\
0.818181818181818	4911.51777705263\\
0.848484848484849	4904.88928453455\\
0.878787878787879	4898.25759244811\\
0.909090909090909	4891.62268468916\\
0.939393939393939	4884.98454503133\\
0.96969696969697	4878.34315712472\\
1	4871.69850449462\\
1.03030303030303	4865.05057054023\\
1.06060606060606	4858.39933853332\\
1.09090909090909	4851.74479161692\\
1.12121212121212	4845.08691280395\\
1.15151515151515	4838.42568497586\\
1.18181818181818	4831.76109088121\\
1.21212121212121	4825.09311313433\\
1.24242424242424	4818.42173421383\\
1.27272727272727	4811.74693646116\\
1.3030303030303	4805.06870207918\\
1.33333333333333	4798.38701313067\\
1.36363636363636	4791.70185153677\\
1.39393939393939	4785.01319907551\\
1.42424242424242	4778.32103738025\\
1.45454545454545	4771.6253479381\\
1.48484848484848	4764.92611208833\\
1.51515151515152	4758.22331102077\\
1.54545454545455	4751.51692577415\\
1.57575757575758	4744.80693723447\\
1.60606060606061	4738.0933261333\\
1.63636363636364	4731.37607304609\\
1.66666666666667	4724.65515839042\\
1.6969696969697	4717.93056242425\\
1.72727272727273	4711.20226524418\\
1.75757575757576	4704.47024678359\\
1.78787878787879	4697.73448681083\\
1.81818181818182	4690.99496492738\\
1.84848484848485	4684.25166056597\\
1.87878787878788	4677.50455298862\\
1.90909090909091	4670.75362128479\\
1.93939393939394	4663.9988443693\\
1.96969696969697	4657.24020098047\\
2	4650.47766967797\\
2.03030303030303	4643.71122884085\\
2.06060606060606	4636.94085666542\\
2.09090909090909	4630.16653116317\\
2.12121212121212	4623.38823015856\\
2.15151515151515	4616.60593128692\\
2.18181818181818	4609.81961199219\\
2.21212121212121	4603.02924952469\\
2.24242424242424	4596.23482093887\\
2.27272727272727	4589.43630309095\\
2.3030303030303	4582.63367263661\\
2.33333333333333	4575.8269060286\\
2.36363636363636	4569.01597951435\\
2.39393939393939	4562.20086913346\\
2.42424242424242	4555.38155071524\\
2.45454545454545	4548.55799987621\\
2.48484848484848	4541.73019201749\\
2.51515151515152	4534.89810232219\\
2.54545454545455	4528.0617057528\\
2.57575757575758	4521.22097704848\\
2.60606060606061	4514.3758907223\\
2.63636363636364	4507.52642105853\\
2.66666666666667	4500.67254210976\\
2.6969696969697	4493.81422769408\\
2.72727272727273	4486.95145139215\\
2.75757575757576	4480.08418654426\\
2.78787878787879	4473.21240624732\\
2.81818181818182	4466.33608335183\\
2.84848484848485	4459.45519045873\\
2.87878787878788	4452.56969991635\\
2.90909090909091	4445.67958381711\\
2.93939393939394	4438.78481399433\\
2.96969696969697	4431.88536201892\\
3	4424.981199196\\
};

\end{axis}
\end{tikzpicture}%

%% file: figure2_tax.tex
%
%
\begin{tikzpicture}

\pgfplotsset{every axis/.append style={
font=\small,
thin,
tick style={ultra thin}}}
\pgfplotsset{every axis y label/.style={
at={(-0.44,0.5)},
xshift=32pt,
rotate=90}}

\begin{axis}[%
width=1.794in,
height=1.03in,
at={(1.358in,0.0in)},
scale only axis,
xmin=0,
xmax=3,
xlabel style={font=\color{white!15!black}},
xlabel={Congestion Surcharge (\$/trip)},
ymin=150,
ymax=190,
ylabel style={font=\color{white!15!black}},
ylabel={ Ride Arrival/min},
axis background/.style={fill=white},
legend style={legend cell align=left, align=left, draw=white!15!black}
]
\addplot [color=black, line width=1.0pt]
  table[row sep=crcr]{%
0	187.27962962963\\
0.0303030303030303	186.927462828644\\
0.0606060606060606	186.575240321977\\
0.0909090909090909	186.222961790944\\
0.121212121212121	185.870626914352\\
0.151515151515152	185.518235368473\\
0.181818181818182	185.165786827023\\
0.212121212121212	184.813280961129\\
0.242424242424242	184.460717439309\\
0.272727272727273	184.108095927439\\
0.303030303030303	183.755416088731\\
0.333333333333333	183.4026775837\\
0.363636363636364	183.049880070141\\
0.393939393939394	182.697023203099\\
0.424242424242424	182.344106634837\\
0.454545454545455	181.991130014812\\
0.484848484848485	181.638092989642\\
0.515151515151515	181.284995203077\\
0.545454545454545	180.93183629597\\
0.575757575757576	180.578615906243\\
0.606060606060606	180.22533366886\\
0.636363636363636	179.871989215794\\
0.666666666666667	179.518582175992\\
0.696969696969697	179.165112175347\\
0.727272727272727	178.811578836665\\
0.757575757575758	178.457981779626\\
0.787878787878788	178.10432062076\\
0.818181818181818	177.750594973403\\
0.848484848484849	177.39680444767\\
0.878787878787879	177.042948650415\\
0.909090909090909	176.689027185199\\
0.939393939393939	176.33503965225\\
0.96969696969697	175.980985648432\\
1	175.626864767203\\
1.03030303030303	175.272676598578\\
1.06060606060606	174.918420729097\\
1.09090909090909	174.564096741777\\
1.12121212121212	174.20970421608\\
1.15151515151515	173.855242727872\\
1.18181818181818	173.500711849381\\
1.21212121212121	173.146111149156\\
1.24242424242424	172.791440192031\\
1.27272727272727	172.436698539075\\
1.3030303030303	172.081885747555\\
1.33333333333333	171.727001370894\\
1.36363636363636	171.372044958621\\
1.39393939393939	171.017016056334\\
1.42424242424242	170.661914205651\\
1.45454545454545	170.306738944163\\
1.48484848484848	169.951489805395\\
1.51515151515152	169.596166318749\\
1.54545454545455	169.240768009465\\
1.57575757575758	168.885294398568\\
1.60606060606061	168.52974500282\\
1.63636363636364	168.174119334674\\
1.66666666666667	167.818416902217\\
1.6969696969697	167.462637209125\\
1.72727272727273	167.106779754606\\
1.75757575757576	166.750844033352\\
1.78787878787879	166.394829535484\\
1.81818181818182	166.038735746497\\
1.84848484848485	165.682562147204\\
1.87878787878788	165.326308213684\\
1.90909090909091	164.969973417222\\
1.93939393939394	164.613557224254\\
1.96969696969697	164.257059096305\\
2	163.900478489934\\
2.03030303030303	163.543814856672\\
2.06060606060606	163.187067642961\\
2.09090909090909	162.830236290093\\
2.12121212121212	162.473320234143\\
2.15151515151515	162.116318905913\\
2.18181818181818	161.759231730859\\
2.21212121212121	161.402058129033\\
2.24242424242424	161.044797515007\\
2.27272727272727	160.687449297813\\
2.3030303030303	160.330012880874\\
2.33333333333333	159.972487661925\\
2.36363636363636	159.614873032955\\
2.39393939393939	159.257168380123\\
2.42424242424242	158.899373083691\\
2.45454545454545	158.541486517951\\
2.48484848484848	158.183508051145\\
2.51515151515152	157.825437045388\\
2.54545454545455	157.467272856594\\
2.57575757575758	157.109014834393\\
2.60606060606061	156.750662322054\\
2.63636363636364	156.392214656398\\
2.66666666666667	156.033671167721\\
2.6969696969697	155.675031179702\\
2.72727272727273	155.316294009326\\
2.75757575757576	154.957458966788\\
2.78787878787879	154.598525355409\\
2.81818181818182	154.239492471548\\
2.84848484848485	153.880359604502\\
2.87878787878788	153.521126036423\\
2.90909090909091	153.161791042214\\
2.93939393939394	152.80235388944\\
2.96969696969697	152.442813838225\\
3	152.083170141155\\
};

\end{axis}
\end{tikzpicture}%

%% file: figure3_tax.tex
%
%
\definecolor{mycolor1}{rgb}{0.60000,0.20000,0.00000}%
\begin{tikzpicture}

\pgfplotsset{every axis y label/.style={
at={(-0.40,0.5)},
xshift=32pt,
rotate=90}}

\begin{axis}[%
width=1.794in,
height=1.03in,
at={(1.358in,0.0in)},
scale only axis,
xmin=0,
xmax=3,
xlabel style={font=\color{white!15!black}},
xlabel={Congestion Surcharge (\$/trip)},
ymin=9,
ymax=17,
ylabel style={font=\color{white!15!black}},
ylabel={Price (\$/trip)},
axis background/.style={fill=white},
legend style={at={(0.0,0.3)}, anchor=south west, legend cell align=left, align=left, draw=white!15!black}
]
\addplot [color=black, line width=1.0pt]
  table[row sep=crcr]{%
0	17\\
0.0303030303030303	16.9753248967496\\
0.0606060606060606	16.9506416173187\\
0.0909090909090909	16.9259501217601\\
0.121212121212121	16.9012503698378\\
0.151515151515152	16.876542321024\\
0.181818181818182	16.8518259344963\\
0.212121212121212	16.8271011691351\\
0.242424242424242	16.8023679835201\\
0.272727272727273	16.7776263359278\\
0.303030303030303	16.7528761843285\\
0.333333333333333	16.7281174863828\\
0.363636363636364	16.7033501994393\\
0.393939393939394	16.6785742805307\\
0.424242424242424	16.6537896863711\\
0.454545454545455	16.6289963733527\\
0.484848484848485	16.6041942975427\\
0.515151515151515	16.5793834146797\\
0.545454545454545	16.5545636801707\\
0.575757575757576	16.5297350490877\\
0.606060606060606	16.5048974761641\\
0.636363636363636	16.4800509157914\\
0.666666666666667	16.4551953220158\\
0.696969696969697	16.4303306485344\\
0.727272727272727	16.4054568486919\\
0.757575757575758	16.3805738754769\\
0.787878787878788	16.355681681518\\
0.818181818181818	16.3307802190805\\
0.848484848484849	16.3058694400623\\
0.878787878787879	16.2809492959901\\
0.909090909090909	16.2560197380157\\
0.939393939393939	16.2310807169118\\
0.96969696969697	16.2061321830684\\
1	16.1811740864886\\
1.03030303030303	16.1562063767843\\
1.06060606060606	16.1312290031724\\
1.09090909090909	16.1062419144706\\
1.12121212121212	16.0812450590927\\
1.15151515151515	16.0562383850448\\
1.18181818181818	16.0312218399208\\
1.21212121212121	16.006195370898\\
1.24242424242424	15.9811589247322\\
1.27272727272727	15.9561124477539\\
1.3030303030303	15.931055885863\\
1.33333333333333	15.9059891845248\\
1.36363636363636	15.8809122887645\\
1.39393939393939	15.8558251431629\\
1.42424242424242	15.8307276918517\\
1.45454545454545	15.805619878508\\
1.48484848484848	15.7805016463496\\
1.51515151515152	15.7553729381301\\
1.54545454545455	15.7302336961334\\
1.57575757575758	15.7050838621689\\
1.60606060606061	15.6799233775656\\
1.63636363636364	15.6547521831675\\
1.66666666666667	15.6295702193275\\
1.6969696969697	15.6043774259022\\
1.72727272727273	15.5791737422463\\
1.75757575757576	15.5539591072069\\
1.78787878787879	15.5287334591175\\
1.81818181818182	15.5034967357926\\
1.84848484848485	15.4782488745214\\
1.87878787878788	15.452989812062\\
1.90909090909091	15.4277194846352\\
1.93939393939394	15.4024378279183\\
1.96969696969697	15.3771447770389\\
2	15.3518402665684\\
2.03030303030303	15.3265242305156\\
2.06060606060606	15.3011966023204\\
2.09090909090909	15.2758573148465\\
2.12121212121212	15.2505063003752\\
2.15151515151515	15.2251434905985\\
2.18181818181818	15.1997688166117\\
2.21212121212121	15.174382208907\\
2.24242424242424	15.1489835973656\\
2.27272727272727	15.1235729112509\\
2.3030303030303	15.0981500792012\\
2.33333333333333	15.0727150292215\\
2.36363636363636	15.0472676886768\\
2.39393939393939	15.0218079842837\\
2.42424242424242	14.9963358421026\\
2.45454545454545	14.9708511875303\\
2.48484848484848	14.945353945291\\
2.51515151515152	14.9198440394289\\
2.54545454545455	14.8943213932992\\
2.57575757575758	14.8687859295601\\
2.60606060606061	14.8432375701639\\
2.63636363636364	14.8176762363483\\
2.66666666666667	14.7921018486276\\
2.6969696969697	14.7665143267834\\
2.72727272727273	14.7409135898557\\
2.75757575757576	14.7152995561337\\
2.78787878787879	14.6896721431458\\
2.81818181818182	14.6640312676503\\
2.84848484848485	14.6383768456258\\
2.87878787878788	14.6127087922609\\
2.90909090909091	14.5870270219441\\
2.93939393939394	14.5613314482541\\
2.96969696969697	14.5356219839484\\
3	14.5098985409537\\
};
\addlegendentry{$p_f$}

\addplot [color=mycolor1, dotted, line width=1.0pt]
  table[row sep=crcr]{%
0	10.2359355687178\\
0.0303030303030303	10.2288487084643\\
0.0606060606060606	10.2217608572752\\
0.0909090909090909	10.2146720111884\\
0.121212121212121	10.2075821662171\\
0.151515151515152	10.2004913183492\\
0.181818181818182	10.1933994635477\\
0.212121212121212	10.1863065977498\\
0.242424242424242	10.1792127168672\\
0.272727272727273	10.1721178167856\\
0.303030303030303	10.1650218933643\\
0.333333333333333	10.1579249424367\\
0.363636363636364	10.150826959809\\
0.393939393939394	10.1437279412609\\
0.424242424242424	10.1366278825448\\
0.454545454545455	10.1295267793857\\
0.484848484848485	10.1224246274811\\
0.515151515151515	10.1153214225004\\
0.545454545454545	10.1082171600851\\
0.575757575757576	10.1011118358482\\
0.606060606060606	10.0940054453741\\
0.636363636363636	10.086897984218\\
0.666666666666667	10.0797894479063\\
0.696969696969697	10.0726798319357\\
0.727272727272727	10.0655691317732\\
0.757575757575758	10.0584573428558\\
0.787878787878788	10.0513444605902\\
0.818181818181818	10.0442304803524\\
0.848484848484849	10.0371153974878\\
0.878787878787879	10.0299992073103\\
0.909090909090909	10.0228819051026\\
0.939393939393939	10.0157634861154\\
0.96969696969697	10.0086439455676\\
1	10.0015232786455\\
1.03030303030303	9.99440148050277\\
1.06060606060606	9.98727854626012\\
1.09090909090909	9.98015447100492\\
1.12121212121212	9.97302924979089\\
1.15151515151515	9.96590287763777\\
1.18181818181818	9.958775349531\\
1.21212121212121	9.95164666042138\\
1.24242424242424	9.94451680522474\\
1.27272727272727	9.93738577882156\\
1.3030303030303	9.93025357605666\\
1.33333333333333	9.92312019173884\\
1.36363636363636	9.9159856206405\\
1.39393939393939	9.90884985749729\\
1.42424242424242	9.90171289700776\\
1.45454545454545	9.89457473383295\\
1.48484848484848	9.88743536259606\\
1.51515151515152	9.880294777882\\
1.54545454545455	9.87315297423709\\
1.57575757575758	9.86600994616859\\
1.60606060606061	9.85886568814435\\
1.63636363636364	9.85172019459238\\
1.66666666666667	9.84457345990046\\
1.6969696969697	9.83742547841571\\
1.72727272727273	9.83027624444419\\
1.75757575757576	9.82312575225046\\
1.78787878787879	9.81597399605713\\
1.81818181818182	9.80882097004448\\
1.84848484848485	9.80166666834997\\
1.87878787878788	9.79451108506779\\
1.90909090909091	9.78735421424843\\
1.93939393939394	9.78019604989821\\
1.96969696969697	9.7730365859788\\
2	9.76587581640675\\
2.03030303030303	9.75871373505304\\
2.06060606060606	9.75155033574254\\
2.09090909090909	9.74438561225355\\
2.12121212121212	9.73721955831728\\
2.15151515151515	9.73005216761736\\
2.18181818181818	9.72288343378933\\
2.21212121212121	9.71571335042007\\
2.24242424242424	9.70854191104733\\
2.27272727272727	9.70136910915915\\
2.3030303030303	9.69419493819332\\
2.33333333333333	9.68701939153683\\
2.36363636363636	9.67984246252535\\
2.39393939393939	9.67266414444257\\
2.42424242424242	9.66548443051969\\
2.45454545454545	9.65830331393486\\
2.48484848484848	9.65112078781245\\
2.51515151515152	9.64393684522261\\
2.54545454545455	9.63675147918056\\
2.57575757575758	9.62956468264597\\
2.60606060606061	9.62237644852239\\
2.63636363636364	9.61518676965653\\
2.66666666666667	9.60799563883766\\
2.6969696969697	9.60080304879695\\
2.72727272727273	9.59360899220675\\
2.75757575757576	9.58641346167999\\
2.78787878787879	9.57921644976941\\
2.81818181818182	9.57201794896691\\
2.84848484848485	9.5648179517028\\
2.87878787878788	9.55761645034514\\
2.90909090909091	9.55041343719891\\
2.93939393939394	9.54320890450535\\
2.96969696969697	9.53600284444116\\
3	9.52879524911774\\
};
\addlegendentry{$p_d$}

\end{axis}
\end{tikzpicture}%

%% file: figure4_tax.tex
%
%
\begin{tikzpicture}

\pgfplotsset{every axis/.append style={
font=\small,
thin,
tick style={ultra thin}}}
\pgfplotsset{every axis y label/.style={
at={(-0.46,0.5)},
xshift=32pt,
rotate=90}}

\begin{axis}[%
width=1.794in,
height=1.03in,
at={(1.358in,0.0in)},
scale only axis,
xmin=0,
xmax=3,
xlabel style={font=\color{white!15!black}},
xlabel={tax},
ymin=19.5,
ymax=23,
ylabel style={font=\color{white!15!black}},
ylabel={Driver Wage/hour},
axis background/.style={fill=white},
legend style={legend cell align=left, align=left, draw=white!15!black}
]
\addplot [color=black, line width=1.0pt]
  table[row sep=crcr]{%
0	22.6\\
0.0303030303030303	22.5709235477228\\
0.0606060606060606	22.5418346315785\\
0.0909090909090909	22.5127331928293\\
0.121212121212121	22.4836191723219\\
0.151515151515152	22.4544925104839\\
0.181818181818182	22.4253531473193\\
0.212121212121212	22.3962010224046\\
0.242424242424242	22.3670360748847\\
0.272727272727273	22.3378582434687\\
0.303030303030303	22.3086674664254\\
0.333333333333333	22.2794636815795\\
0.363636363636364	22.2502468263066\\
0.393939393939394	22.2210168375295\\
0.424242424242424	22.1917736517132\\
0.454545454545455	22.1625172048604\\
0.484848484848485	22.1332474325073\\
0.515151515151515	22.1039642697186\\
0.545454545454545	22.074667651083\\
0.575757575757576	22.0453575107085\\
0.606060606060606	22.0160337822171\\
0.636363636363636	21.9866963987407\\
0.666666666666667	21.9573452929155\\
0.696969696969697	21.9279803968773\\
0.727272727272727	21.8986016422565\\
0.757575757575758	21.8692089601727\\
0.787878787878788	21.8398022812296\\
0.818181818181818	21.8103815355101\\
0.848484848484849	21.7809466525702\\
0.878787878787879	21.7514975614345\\
0.909090909090909	21.7220341905898\\
0.939393939393939	21.6925564679804\\
0.96969696969697	21.6630643210019\\
1	21.6335576764956\\
1.03030303030303	21.6040364607432\\
1.06060606060606	21.5745005994603\\
1.09090909090909	21.5449500177906\\
1.12121212121212	21.5153846403005\\
1.15151515151515	21.4858043909722\\
1.18181818181818	21.4562091931979\\
1.21212121212121	21.4265989697739\\
1.24242424242424	21.3969736428935\\
1.27272727272727	21.3673331341411\\
1.3030303030303	21.3376773644858\\
1.33333333333333	21.3080062542743\\
1.36363636363636	21.2783197232246\\
1.39393939393939	21.2486176904192\\
1.42424242424242	21.2189000742979\\
1.45454545454545	21.1891667926515\\
1.48484848484848	21.1594177626139\\
1.51515151515152	21.1296529006555\\
1.54545454545455	21.0998721225758\\
1.57575757575758	21.070075343496\\
1.60606060606061	21.0402624778516\\
1.63636363636364	21.0104334393847\\
1.66666666666667	20.9805881411364\\
1.6969696969697	20.9507264954391\\
1.72727272727273	20.9208484139085\\
1.75757575757576	20.8909538074356\\
1.78787878787879	20.8610425861785\\
1.81818181818182	20.8311146595544\\
1.84848484848485	20.8011699362308\\
1.87878787878788	20.7712083241177\\
1.90909090909091	20.7412297303582\\
1.93939393939394	20.7112340613203\\
1.96969696969697	20.6812212225882\\
2	20.6511911189525\\
2.03030303030303	20.621143654402\\
2.06060606060606	20.5910787321139\\
2.09090909090909	20.5609962544448\\
2.12121212121212	20.5308961229205\\
2.15151515151515	20.5007782382272\\
2.18181818181818	20.4706425002011\\
2.21212121212121	20.4404888078186\\
2.24242424242424	20.4103170591862\\
2.27272727272727	20.3801271515304\\
2.3030303030303	20.3499189811869\\
2.33333333333333	20.3196924435905\\
2.36363636363636	20.2894474332639\\
2.39393939393939	20.2591838438072\\
2.42424242424242	20.2289015678867\\
2.45454545454545	20.1986004972234\\
2.48484848484848	20.1682805225822\\
2.51515151515152	20.1379415337598\\
2.54545454545455	20.107583419573\\
2.57575757575758	20.0772060678469\\
2.60606060606061	20.0468093654029\\
2.63636363636364	20.0163931980461\\
2.66666666666667	19.9859574505529\\
2.6969696969697	19.9555020066583\\
2.72727272727273	19.9250267490429\\
2.75757575757576	19.8945315593202\\
2.78787878787879	19.8640163180226\\
2.81818181818182	19.8334809045883\\
2.84848484848485	19.8029251973475\\
2.87878787878788	19.7723490735086\\
2.90909090909091	19.7417524091433\\
2.93939393939394	19.7111350791731\\
2.96969696969697	19.6804969573541\\
3	19.649837916262\\
};

\end{axis}
\end{tikzpicture}%

%% file: figure5_tax.tex
%
%
\begin{tikzpicture}

\pgfplotsset{every axis/.append style={
font=\small,
thin,
tick style={ultra thin}}}
\pgfplotsset{every axis y label/.style={
at={(-0.45,0.5)},
xshift=32pt,
rotate=90}}

\begin{axis}[%
width=1.794in,
height=1.03in,
at={(1.358in,0in)},
scale only axis,
xmin=0,
xmax=3,
xlabel style={font=\color{white!15!black}},
xlabel={tax},
ymin=33.3,
ymax=34.3,
ylabel style={font=\color{white!15!black}},
ylabel={Total Cost (\$/trip)},
axis background/.style={fill=white},
legend style={legend cell align=left, align=left, draw=white!15!black}
]
\addplot [color=black, line width=1.0pt]
  table[row sep=crcr]{%
0	33.3850756082045\\
0.0303030303030303	33.3939525057177\\
0.0606060606060606	33.4028308073769\\
0.0909090909090909	33.411710521215\\
0.121212121212121	33.4205916553284\\
0.151515151515152	33.4294742178769\\
0.181818181818182	33.438358217085\\
0.212121212121212	33.4472436612426\\
0.242424242424242	33.4561305587052\\
0.272727272727273	33.4650189178947\\
0.303030303030303	33.4739087473005\\
0.333333333333333	33.4828000554799\\
0.363636363636364	33.4916928510585\\
0.393939393939394	33.5005871427316\\
0.424242424242424	33.5094829392643\\
0.454545454545455	33.5183802494927\\
0.484848484848485	33.5272790823242\\
0.515151515151515	33.5361794467386\\
0.545454545454545	33.5450813517889\\
0.575757575757576	33.5539848066017\\
0.606060606060606	33.5628898203784\\
0.636363636363636	33.5717964023956\\
0.666666666666667	33.5807045620062\\
0.696969696969697	33.5896143086402\\
0.727272727272727	33.5985256518054\\
0.757575757575758	33.6074386010882\\
0.787878787878788	33.6163531661546\\
0.818181818181818	33.6252693567511\\
0.848484848484849	33.6341871827054\\
0.878787878787879	33.6431066539272\\
0.909090909090909	33.6520277804094\\
0.939393939393939	33.6609505722288\\
0.96969696969697	33.6698750395472\\
1	33.6788011926119\\
1.03030303030303	33.6877290417574\\
1.06060606060606	33.6966585974055\\
1.09090909090909	33.7055898700669\\
1.12121212121212	33.7145228703418\\
1.15151515151515	33.7234576089213\\
1.18181818181818	33.7323940965879\\
1.21212121212121	33.7413323442171\\
1.24242424242424	33.7502723627778\\
1.27272727272727	33.759214163334\\
1.3030303030303	33.7681577570454\\
1.33333333333333	33.7771031551688\\
1.36363636363636	33.7860503690588\\
1.39393939393939	33.7949994101693\\
1.42424242424242	33.8039502900547\\
1.45454545454545	33.8129030203706\\
1.48484848484848	33.8218576128752\\
1.51515151515152	33.8308140794306\\
1.54545454545455	33.839772432004\\
1.57575757575758	33.8487326826685\\
1.60606060606061	33.8576948436049\\
1.63636363636364	33.8666589271027\\
1.66666666666667	33.8756249455611\\
1.6969696969697	33.8845929114909\\
1.72727272727273	33.8935628375152\\
1.75757575757576	33.9025347363711\\
1.78787878787879	33.9115086209108\\
1.81818181818182	33.9204845041032\\
1.84848484848485	33.9294623990353\\
1.87878787878788	33.9384423189131\\
1.90909090909091	33.947424277064\\
1.93939393939394	33.956408286937\\
1.96969696969697	33.9653943621052\\
2	33.974382516267\\
2.03030303030303	33.9833727632473\\
2.06060606060606	33.9923651169994\\
2.09090909090909	34.0013595916064\\
2.12121212121212	34.0103562012828\\
2.15151515151515	34.0193549603762\\
2.18181818181818	34.0283558833687\\
2.21212121212121	34.037358984879\\
2.24242424242424	34.0463642796635\\
2.27272727272727	34.0553717826184\\
2.3030303030303	34.0643815087816\\
2.33333333333333	34.0733934733339\\
2.36363636363636	34.0824076916011\\
2.39393939393939	34.0914241790562\\
2.42424242424242	34.1004429513203\\
2.45454545454545	34.1094640241657\\
2.48484848484848	34.1184874135166\\
2.51515151515152	34.1275131354518\\
2.54545454545455	34.1365412062067\\
2.57575757575758	34.1455716421746\\
2.60606060606061	34.1546044599096\\
2.63636363636364	34.1636396761279\\
2.66666666666667	34.1726773077107\\
2.6969696969697	34.1817173717055\\
2.72727272727273	34.1907598853288\\
2.75757575757576	34.1998048659681\\
2.78787878787879	34.2088523311844\\
2.81818181818182	34.2179022987138\\
2.84848484848485	34.2269547864709\\
2.87878787878788	34.23600981255\\
2.90909090909091	34.2450673952283\\
2.93939393939394	34.2541275529679\\
2.96969696969697	34.2631903044185\\
3	34.2722556684199\\
};

\end{axis}
\end{tikzpicture}%

%% file: figure6_tax.tex
%
%
\begin{tikzpicture}

\pgfplotsset{every axis/.append style={
font=\small,
thin,
tick style={ultra thin}}}
\pgfplotsset{every axis y label/.style={
at={(-0.4,0.5)},
xshift=32pt,
rotate=90}}

\begin{axis}[%
width=1.794in,
height=1.03in,
at={(1.358in,0.0in)},
scale only axis,
xmin=0,
xmax=3,
xlabel style={font=\color{white!15!black}},
xlabel={tax},
ymin=44897.4509229417,
ymax=79897.4509229417,
ylabel style={font=\color{white!15!black}},
ylabel={Platform Rent},
axis background/.style={fill=white},
legend style={legend cell align=left, align=left, draw=white!15!black}
]
\addplot [color=black, line width=1.0pt]
  table[row sep=crcr]{%
0	76006.2888888889\\
0.0303030303030303	75666.1006146016\\
0.0606060606060606	75326.5526941821\\
0.0909090909090909	74987.645229203\\
0.121212121212121	74649.3783218183\\
0.151515151515152	74311.7520747685\\
0.181818181818182	73974.7665913848\\
0.212121212121212	73638.4219755941\\
0.242424242424242	73302.7183319236\\
0.272727272727273	72967.6557655057\\
0.303030303030303	72633.2343820825\\
0.333333333333333	72299.4542880109\\
0.363636363636364	71966.315590268\\
0.393939393939394	71633.818396455\\
0.424242424242424	71301.962814803\\
0.454545454545455	70970.7489541778\\
0.484848484848485	70640.1769240854\\
0.515151515151515	70310.2468346769\\
0.545454545454545	69980.9587967533\\
0.575757575757576	69652.3129217717\\
0.606060606060606	69324.3093218504\\
0.636363636363636	68996.9481097733\\
0.666666666666667	68670.229398997\\
0.696969696969697	68344.153303655\\
0.727272727272727	68018.7199385642\\
0.757575757575758	67693.9294192297\\
0.787878787878788	67369.7818618513\\
0.818181818181818	67046.2773833289\\
0.848484848484849	66723.4161012682\\
0.878787878787879	66401.1981339869\\
0.909090909090909	66079.62360052\\
0.939393939393939	65758.692620627\\
0.96969696969697	65438.4053147971\\
1	65118.7618042554\\
1.03030303030303	64799.7622109692\\
1.06060606060606	64481.4066576549\\
1.09090909090909	64163.6952677833\\
1.12121212121212	63846.6281655867\\
1.15151515151515	63530.2054760657\\
1.18181818181818	63214.4273249946\\
1.21212121212121	62899.2938389295\\
1.24242424242424	62584.805145214\\
1.27272727272727	62270.9613719861\\
1.3030303030303	61957.762648186\\
1.33333333333333	61645.2091035618\\
1.36363636363636	61333.3008686771\\
1.39393939393939	61022.0380749183\\
1.42424242424242	60711.4208545014\\
1.45454545454545	60401.4493404799\\
1.48484848484848	60092.1236667512\\
1.51515151515152	59783.4439680649\\
1.54545454545455	59475.41038003\\
1.57575757575758	59168.0230391223\\
1.60606060606061	58861.2820826925\\
1.63636363636364	58555.1876489741\\
1.66666666666667	58249.7398770903\\
1.6969696969697	57944.9389070632\\
1.72727272727273	57640.7848798215\\
1.75757575757576	57337.2779372081\\
1.78787878787879	57034.4182219889\\
1.81818181818182	56732.2058778614\\
1.84848484848485	56430.6410494624\\
1.87878787878788	56129.7238823771\\
1.90909090909091	55829.4545231477\\
1.93939393939394	55529.8331192819\\
1.96969696969697	55230.8598192624\\
2	54932.5347725549\\
2.03030303030303	54634.8581296182\\
2.06060606060606	54337.8300419127\\
2.09090909090909	54041.4506619098\\
2.12121212121212	53745.720143102\\
2.15151515151515	53450.6386400112\\
2.18181818181818	53156.2063081997\\
2.21212121212121	52862.4233042788\\
2.24242424242424	52569.2897859195\\
2.27272727272727	52276.8059118623\\
2.3030303030303	51984.9718419272\\
2.33333333333333	51693.7877370239\\
2.36363636363636	51403.2537591626\\
2.39393939393939	51113.370071464\\
2.42424242424242	50824.1368381701\\
2.45454545454545	50535.5542246559\\
2.48484848484848	50247.6223974387\\
2.51515151515152	49960.3415241904\\
2.54545454545455	49673.7117737482\\
2.57575757575758	49387.733316126\\
2.60606060606061	49102.4063225261\\
2.63636363636364	48817.7309653507\\
2.66666666666667	48533.7074182136\\
2.6969696969697	48250.3358559524\\
2.72727272727273	47967.6164546402\\
2.75757575757576	47685.5493915985\\
2.78787878787879	47404.1348454086\\
2.81818181818182	47123.3729959256\\
2.84848484848485	46843.2640242896\\
2.87878787878788	46563.80811294\\
2.90909090909091	46285.0054456275\\
2.93939393939394	46006.8562074279\\
2.96969696969697	45729.3605847557\\
3	45452.5187653771\\
};

\end{axis}
\end{tikzpicture}%

%% file: figure_subsidy1.tex
%
%
\definecolor{mycolor1}{rgb}{0.60000,0.20000,0.00000}%
\begin{tikzpicture}

\pgfplotsset{every axis y label/.style={
at={(-0.42,0.5)},
xshift=32pt,
rotate=90}}

\begin{axis}[%
width=1.794in,
height=1.03in,
at={(1.358in,0.0in)},
scale only axis,
xmin=0,
xmax=900,
xlabel style={font=\color{white!15!black}},
xlabel={Subsidy Budget (\$/min)},
ymin=0,
ymax=3.5,
ylabel style={font=\color{white!15!black}},
ylabel={Subsidy (\$/trip)},
axis background/.style={fill=white},
legend style={at={(0.0,0.59)}, anchor=south west, legend cell align=left, align=left, draw=white!15!black}
]
\addplot [color=black, line width=1.0pt]
  table[row sep=crcr]{%
868.738666666668	1.34784760989307\\
859.96352861953	1.33929472224596\\
851.188390572392	1.33074183459886\\
842.413252525254	1.32200281229132\\
833.638114478116	1.3134452712777\\
824.862976430977	1.30470624897016\\
816.087838383839	1.29596722666262\\
807.312700336701	1.28722820435508\\
798.537562289563	1.27848918204754\\
789.762424242425	1.26975015974\\
780.987286195287	1.26082500277202\\
772.212148148149	1.25189984580404\\
763.437010101011	1.24297468883607\\
754.661872053873	1.23404953186809\\
745.886734006735	1.2251197215336\\
737.111595959597	1.21600842990518\\
728.336457912459	1.20689713827677\\
719.561319865321	1.19778584664835\\
710.786181818183	1.18867455501994\\
702.011043771045	1.17937712873108\\
693.235905723907	1.17007970244223\\
684.460767676769	1.16078227615338\\
675.685629629631	1.15148484986452\\
666.910491582492	1.14218742357567\\
658.135353535354	1.13270386262638\\
649.360215488216	1.12322030167709\\
640.585077441078	1.1137367407278\\
631.80993939394	1.10406704511807\\
623.034801346802	1.09440200287485\\
614.259663299664	1.08473230726512\\
605.484525252526	1.07506261165539\\
596.709387205388	1.06520678138523\\
587.93424915825	1.05535095111506\\
579.159111111112	1.0454951208449\\
570.383973063974	1.03563929057473\\
561.608835016836	1.02559732564413\\
552.833696969698	1.01556001408004\\
544.05855892256	1.005331914489\\
535.283420875422	0.995289949558392\\
526.508282828283	0.985061849967351\\
517.733144781145	0.974647615715873\\
508.958006734007	0.964238034830906\\
500.182868686869	0.953823800579428\\
491.407730639731	0.943409566327949\\
482.632592592593	0.932809197416034\\
473.857454545455	0.922213481870629\\
465.082316498317	0.911426978298275\\
456.307178451179	0.900640474725922\\
447.532040404041	0.88985862452008\\
438.756902356903	0.878885986287289\\
429.981764309765	0.867727213394061\\
421.206626262627	0.856759228527781\\
412.431488215489	0.845414320974115\\
403.656350168351	0.834260201447398\\
394.881212121213	0.822915293893732\\
386.106074074075	0.811384251679629\\
377.330936026937	0.799857862832036\\
368.555797979798	0.788140685957496\\
359.78065993266	0.776428162449466\\
351.005521885522	0.764529504280998\\
342.230383838384	0.75262619274602\\
333.455245791246	0.740541399917115\\
324.680107744108	0.728265819061262\\
315.90496969697	0.71599489157192\\
307.129831649832	0.70353782942214\\
298.354693602694	0.691076113905849\\
289.579555555556	0.678246782435194\\
280.804417508418	0.665417450964539\\
272.02927946128	0.652401984833447\\
263.254141414142	0.639381865335843\\
254.479003367004	0.625994129883876\\
245.703865319866	0.612606394431909\\
236.928727272728	0.598846389659066\\
228.153589225589	0.585086384886224\\
219.378451178452	0.571140245452944\\
210.603313131313	0.556821836698789\\
201.828175084175	0.542503427944634\\
193.053037037037	0.527812749869604\\
184.277898989899	0.512940590500648\\
175.502760942761	0.497691508444305\\
166.727622895623	0.482256291727526\\
157.952484848485	0.466639593716819\\
149.177346801347	0.450645973018727\\
140.402208754209	0.434098601705833\\
131.627070707071	0.41736044236599\\
122.851932659933	0.400254667071784\\
114.076794612795	0.382590487796265\\
105.301656565657	0.364554039199871\\
96.5265185185186	0.345959186622164\\
87.7513804713806	0.326619795402708\\
78.9762424242425	0.30616100802532\\
70.2011043771045	0.284148526947444\\
61.4259663299664	0.260148054626522\\
52.6508282828283	0.233725293519998\\
43.8756902356903	0.204445946085314\\
35.1005521885522	0.171875714779913\\
26.3254141414142	0.135580302061239\\
17.5502760942761	0.0951254103867329\\
8.77513804713806	0.0500767422138391\\
0	0\\
};
\addlegendentry{rider}

\addplot [color=mycolor1, dotted, line width=1.0pt]
  table[row sep=crcr]{%
868.738666666668	3.29336873702595\\
859.96352861953	3.27686378061178\\
851.188390572392	3.26023997444558\\
842.413252525254	3.2436816701666\\
833.638114478116	3.226819423091\\
824.862976430977	3.21001431966383\\
816.087838383839	3.19308295333558\\
807.312700336701	3.17602334788765\\
798.537562289563	3.15883347419992\\
789.762424242425	3.14151124822146\\
780.987286195287	3.12424066355929\\
772.212148148149	3.10683338523608\\
763.437010101011	3.08928715138\\
754.661872053873	3.07159963605004\\
745.886734006735	3.05377310000299\\
737.111595959597	3.03598190927958\\
728.336457912459	3.01804204840709\\
719.561319865321	2.9999509080256\\
710.786181818183	2.98170580029735\\
702.011043771045	2.96349009028474\\
693.235905723907	2.94511478834534\\
684.460767676769	2.92657695049987\\
675.685629629631	2.90787353963799\\
666.910491582492	2.88900142129148\\
658.135353535354	2.87014349411023\\
649.360215488216	2.85111028086388\\
640.585077441078	2.83189832760569\\
631.80993939394	2.81269019838059\\
623.034801346802	2.79329141203149\\
614.259663299664	2.7737074539863\\
605.484525252526	2.75392971471192\\
596.709387205388	2.73414023106684\\
587.93424915825	2.71414862206629\\
579.159111111112	2.69395048484676\\
570.383973063974	2.67354125158777\\
561.608835016836	2.65310231574661\\
552.833696969698	2.63243796500333\\
544.05855892256	2.61173852334269\\
535.283420875422	2.59062174697928\\
526.508282828283	2.56945436809416\\
517.733144781145	2.54823062183306\\
508.958006734007	2.52675371334454\\
500.182868686869	2.50502668691761\\
491.407730639731	2.48303835461928\\
482.632592592593	2.46096802807887\\
473.857454545455	2.43861778916777\\
465.082316498317	2.41617562276454\\
456.307178451179	2.39344293564906\\
447.532040404041	2.37040689780601\\
438.756902356903	2.34725437741416\\
429.981764309765	2.32397172247307\\
421.206626262627	2.30017255734169\\
412.431488215489	2.27641465126221\\
403.656350168351	2.25211982844411\\
394.881212121213	2.22765863495401\\
386.106074074075	2.20301478048669\\
377.330936026937	2.17798518336142\\
368.555797979798	2.15275228233114\\
359.78065993266	2.1271068649575\\
351.005521885522	2.10122511411646\\
342.230383838384	2.07491003682011\\
333.455245791246	2.04832196108678\\
324.680107744108	2.02145261461153\\
315.90496969697	1.9940877961384\\
307.129831649832	1.96639823737762\\
298.354693602694	1.93818097522243\\
289.579555555556	1.90977591866197\\
280.804417508418	1.88079065745151\\
272.02927946128	1.85138460234829\\
263.254141414142	1.82134740005303\\
254.479003367004	1.79101074115313\\
245.703865319866	1.75997304718706\\
236.928727272728	1.72856959765304\\
228.153589225589	1.69638771234833\\
219.378451178452	1.66356922925864\\
210.603313131313	1.63025155675379\\
201.828175084175	1.59600861798595\\
193.053037037037	1.56115301594356\\
184.277898989899	1.52542749812304\\
175.502760942761	1.48895310228418\\
166.727622895623	1.45145529521656\\
157.952484848485	1.41283459937213\\
149.177346801347	1.37317837407576\\
140.402208754209	1.33253934761319\\
131.627070707071	1.29041088568218\\
122.851932659933	1.24680193928428\\
114.076794612795	1.20170505998613\\
105.301656565657	1.15469719452929\\
96.5265185185186	1.10567734858767\\
87.7513804713806	1.05447796597341\\
78.9762424242425	0.999463472295689\\
70.2011043771045	0.938617788576696\\
61.4259663299664	0.869924835838602\\
52.6508282828283	0.791368535103585\\
43.8756902356903	0.700932807393822\\
35.1005521885522	0.59660157373149\\
26.3254141414142	0.476358755138764\\
17.5502760942761	0.338188272637823\\
8.77513804713806	0.180074047250843\\
0	0\\
};
\addlegendentry{driver}

\end{axis}
\end{tikzpicture}%

%% file: figure_subsidy2.tex
%
%
\begin{tikzpicture}

\pgfplotsset{every axis/.append style={
font=\small,
thin,
tick style={ultra thin}}}
\pgfplotsset{every axis y label/.style={
at={(-0.45,0.5)},
xshift=32pt,
rotate=90}}

\begin{axis}[%
width=1.794in,
height=1.03in,
at={(1.358in,0.0in)},
scale only axis,
xmin=0,
xmax=900,
xlabel style={font=\color{white!15!black}},
xlabel={Subsidy Budget (\$/min)},
ymin=180,
ymax=360,
ylabel style={font=\color{white!15!black}},
ylabel={Ride Arrival/min},
axis background/.style={fill=white},
legend style={legend cell align=left, align=left, draw=white!15!black}
]
\addplot [color=black, line width=1.0pt]
  table[row sep=crcr]{%
868.738666666668	351.113194863944\\
859.96352861953	350.288938552904\\
851.188390572392	349.460431912881\\
842.413252525254	348.627608862993\\
833.638114478116	347.790401613788\\
824.862976430977	346.948740488315\\
816.087838383839	346.10255404541\\
807.312700336701	345.25176890933\\
798.537562289563	344.396309672489\\
789.762424242425	343.536098816281\\
780.987286195287	342.671056706428\\
772.212148148149	341.801101442035\\
763.437010101011	340.926148824436\\
754.661872053873	340.046112190253\\
745.886734006735	339.160902311685\\
737.111595959597	338.270427402051\\
728.336457912459	337.374592861182\\
719.561319865321	336.473301280751\\
710.786181818183	335.566452228083\\
702.011043771045	334.6539421821\\
693.235905723907	333.7356643444\\
684.460767676769	332.811508585385\\
675.685629629631	331.881361181749\\
666.910491582492	330.945104651006\\
658.135353535354	330.002617727028\\
649.360215488216	329.053775007269\\
640.585077441078	328.098446860998\\
631.80993939394	327.136499225054\\
623.034801346802	326.167793290049\\
614.259663299664	325.192185419783\\
605.484525252526	324.20952676541\\
596.709387205388	323.219663019858\\
587.93424915825	322.222434274763\\
579.159111111112	321.217674599977\\
570.383973063974	320.205211670265\\
561.608835016836	319.184866566841\\
552.833696969698	318.156453314332\\
544.05855892256	317.119778519718\\
535.283420875422	316.074641031849\\
526.508282828283	315.02083127079\\
517.733144781145	313.958130998237\\
508.958006734007	312.886312561493\\
500.182868686869	311.805138572337\\
491.407730639731	310.714361061724\\
482.632592592593	309.613720977265\\
473.857454545455	308.502947387985\\
465.082316498317	307.381756826869\\
456.307178451179	306.249852266554\\
447.532040404041	305.106922349812\\
438.756902356903	303.952640457753\\
429.981764309765	302.786663374446\\
421.206626262627	301.60863039864\\
412.431488215489	300.418161811718\\
403.656350168351	299.214857696591\\
394.881212121213	297.998296080203\\
386.106074074075	296.768031557258\\
377.330936026937	295.52359311095\\
368.555797979798	294.264482218348\\
359.78065993266	292.990170461127\\
351.005521885522	291.700096924993\\
342.230383838384	290.393665416231\\
333.455245791246	289.070241208394\\
324.680107744108	287.729147479078\\
315.90496969697	286.369661353069\\
307.129831649832	284.991009267189\\
298.354693602694	283.592361964066\\
289.579555555556	282.172828787331\\
280.804417508418	280.731450833406\\
272.02927946128	279.267193674379\\
263.254141414142	277.778938720285\\
254.479003367004	276.26547321973\\
245.703865319866	274.725478788799\\
236.928727272728	273.157518155616\\
228.153589225589	271.560019608605\\
219.378451178452	269.931258453471\\
210.603313131313	268.269335698351\\
201.828175084175	266.572152107992\\
193.053037037037	264.837377426837\\
184.277898989899	263.062413088938\\
175.502760942761	261.244346999433\\
166.727622895623	259.379897801728\\
157.952484848485	257.465345852975\\
149.177346801347	255.496446527834\\
140.402208754209	253.468320558048\\
131.627070707071	251.375312344208\\
122.851932659933	249.210805455652\\
114.076794612795	246.966977188662\\
105.301656565657	244.63446669376\\
96.5265185185186	242.201914442431\\
87.7513804713806	239.655309202957\\
78.9762424242425	236.914098776969\\
70.2011043771045	233.880059856274\\
61.4259663299664	230.45496913268\\
52.6508282828283	226.540603297996\\
43.8756902356903	222.038739044028\\
35.1005521885522	216.851153062586\\
26.3254141414142	210.879622045476\\
17.5502760942761	204.025922684507\\
8.77513804713806	196.191831671487\\
0	187.279125698224\\
};

\end{axis}
\end{tikzpicture}%

%% file: figure_subsidy3.tex
%
%
\begin{tikzpicture}

\pgfplotsset{every axis/.append style={
font=\small,
thin,
tick style={ultra thin}}}
\pgfplotsset{every axis y label/.style={
at={(-0.45,0.5)},
xshift=32pt,
rotate=90}}

\begin{axis}[%
width=1.794in,
height=1.03in,
at={(1.358in,0.0in)},
scale only axis,
xmin=0,
xmax=900,
xlabel style={font=\color{white!15!black}},
xlabel={Subsidy Budget (\$/min)},
ymin=5000,
ymax=8500,
ytick={5000,6000,7000, 8000},
yticklabels={{5K},{6K},{7K}, {8K}},
ylabel style={font=\color{white!15!black}},
ylabel={Number of Drivers},
axis background/.style={fill=white},
legend style={legend cell align=left, align=left, draw=white!15!black}
]
\addplot [color=black, line width=1.0pt]
  table[row sep=crcr]{%
868.738666666668	8011.48754384047\\
859.96352861953	7997.19580632433\\
851.188390572392	7982.81781370725\\
842.413252525254	7968.40726350946\\
833.638114478116	7953.85411270354\\
824.862976430977	7939.26422205928\\
816.087838383839	7924.58252445793\\
807.312700336701	7909.80755548948\\
798.537562289563	7894.93781066696\\
789.762424242425	7879.97174387124\\
780.987286195287	7864.96255435749\\
772.212148148149	7849.85375370522\\
763.437010101011	7834.64366208932\\
754.661872053873	7819.33055106548\\
745.886734006735	7803.91400785932\\
737.111595959597	7788.44408416208\\
728.336457912459	7772.86557542313\\
719.561319865321	7757.17653914468\\
710.786181818183	7741.37497316919\\
702.011043771045	7725.51328561756\\
693.235905723907	7709.53480302024\\
684.460767676769	7693.43732965178\\
675.685629629631	7677.21859888679\\
666.910491582492	7660.87626994758\\
658.135353535354	7644.46221042969\\
649.360215488216	7627.91956122373\\
640.585077441078	7611.24573942696\\
631.80993939394	7594.49224080255\\
623.034801346802	7577.60070123804\\
614.259663299664	7560.57096624254\\
605.484525252526	7543.39871692037\\
596.709387205388	7526.13488493767\\
587.93424915825	7508.72219848265\\
579.159111111112	7491.15735295636\\
570.383973063974	7473.4369175848\\
561.608835016836	7455.61116507216\\
552.833696969698	7437.6211275451\\
544.05855892256	7419.51939247977\\
535.283420875422	7401.19307258355\\
526.508282828283	7382.7454451433\\
517.733144781145	7364.17203060515\\
508.958006734007	7345.41325720658\\
500.182868686869	7326.4670798438\\
491.407730639731	7307.3272249338\\
482.632592592593	7288.04196026606\\
473.857454545455	7268.55104489052\\
465.082316498317	7248.90483229601\\
456.307178451179	7229.04265797491\\
447.532040404041	7208.95699936239\\
438.756902356903	7188.69719285732\\
429.981764309765	7168.25492665724\\
421.206626262627	7147.51541757801\\
412.431488215489	7126.6329586132\\
403.656350168351	7105.43776055526\\
394.881212121213	7084.03004720189\\
386.106074074075	7062.39947314138\\
377.330936026937	7040.48240449363\\
368.555797979798	7018.32429994244\\
359.78065993266	6995.85918150439\\
351.005521885522	6973.12982912146\\
342.230383838384	6950.0729887209\\
333.455245791246	6926.7259035233\\
324.680107744108	6903.07756742048\\
315.90496969697	6879.05846034999\\
307.129831649832	6854.70681625512\\
298.354693602694	6829.955062472\\
289.579555555556	6804.88719342773\\
280.804417508418	6779.38091752816\\
272.02927946128	6753.46770315577\\
263.254141414142	6727.07450565873\\
254.479003367004	6700.2783622632\\
245.703865319866	6672.95041457368\\
236.928727272728	6645.16553221178\\
228.153589225589	6616.78883723959\\
219.378451178452	6587.83764879181\\
210.603313131313	6558.32549385395\\
201.828175084175	6528.10752579972\\
193.053037037037	6497.24006150766\\
184.277898989899	6465.61859329455\\
175.502760942761	6433.23874065659\\
166.727622895623	6399.98263384478\\
157.952484848485	6365.77498739624\\
149.177346801347	6330.58468746322\\
140.402208754209	6294.36251900725\\
131.627070707071	6256.89691944256\\
122.851932659933	6218.10436317533\\
114.076794612795	6177.88074169747\\
105.301656565657	6135.98941236935\\
96.5265185185186	6092.25467936458\\
87.7513804713806	6046.44777035602\\
78.9762424242425	5997.06943291365\\
70.2011043771045	5942.28556560209\\
61.4259663299664	5880.26206698597\\
52.6508282828283	5809.16483562989\\
43.8756902356903	5727.15977009849\\
35.1005521885522	5632.41276895637\\
26.3254141414142	5523.08973076816\\
17.5502760942761	5397.35655409848\\
8.77513804713806	5253.37913751194\\
0	5089.32337957316\\
};

\end{axis}
\end{tikzpicture}%